\batchmode
\documentclass[11pt]{article}

\usepackage{epsfig}
\usepackage{graphicx}
\usepackage{color}
\usepackage{mathtools}

\newtheorem{theorem}{Theorem}
\newtheorem{lemma}{Lemma}
\newtheorem{corollary}{Corollary}

\newtheorem{definition}{Definition}

\newcommand{\be}{\begin{equation}}
\newcommand{\ee}{\end{equation}}
\newcommand{\bea}{\begin{eqnarray}}
\newcommand{\eea}{\end{eqnarray}}
\newcommand{\beas}{\begin{eqnarray*}}
\newcommand{\eeas}{\end{eqnarray*}}
\newcommand{\ba}{\begin{array}}
\newcommand{\ea}{\end{array}}

\DeclarePairedDelimiter{\floor}{\lfloor}{\rfloor}
\DeclarePairedDelimiter{\ceil}{\lceil}{\rceil}

\definecolor{armygreen}{rgb}{0.29, 0.33, 0.13}

\newcommand{\real}{\mbox{$\mathbb{R}$}}
\newcommand{\Natural}{\mbox{$\mathrm{I\!N}$}}

\newcommand{\eps}{\ensuremath{\epsilon}}

\def\XXint#1#2#3{{\setbox0=\hbox{$#1{#2#3}{\int}$}
     \vcenter{\hbox{$#2#3$}}\kern-.5\wd0}}

\newcommand{\mcF}{\ensuremath{\mathcal{F}}}

\newcommand{\mcH}{\ensuremath{\mathcal{H}}}

\newcommand{\mcL}{\ensuremath{\mathcal{L}}}

\newcommand{\mrI}{\ensuremath{\mathrm{I}}}
\newcommand{\mrJ}{\ensuremath{\mathrm{J}}}

\newcommand{\wtilde}{\ensuremath{\widetilde}}

\def\qed{\hbox{\vrule width 6pt height 6pt depth 0pt}}

\usepackage{amsmath}
\usepackage{amssymb}

\title{Regularity of the solution to fractional diffusion, advection, reaction equations} 
\author{
	V.J. Ervin\thanks{School of Mathematical and Statistical Sciences,
	  Clemson University, Clemson, South Carolina 29634-0975, USA.
	  email: {\tt vjervin@clemson.edu}. }
            } 

\date{\today}

\begin{document}
\maketitle

\begin{abstract}
In this report we investigate the regularity of the solution to the fractional diffusion, advection, reaction equation
on a bounded domain in $\mathbb{R}^{1}$. The analysis is performed in the weighted Sobolev spaces,
$H_{(a , b)}^{s}(\mrI)$. Three different characterizations of $H_{(a , b)}^{s}(\mrI)$ are presented, together 
with needed embedding theorems for these spaces. The analysis shows that the regularity of the solution is
bounded by the endpoint behavior of the solution, which is determined by the parameters $\alpha$ and $r$
defining the fractional diffusion operator. Additionally, the analysis shows that for a sufficiently smooth right
hand side function, the regularity of the solution to fractional diffusion reaction equation is lower than that of
the fractional diffusion equation. Also, the regularity of the solution to fractional diffusion 
advection reaction equation is two orders lower than that of the fractional diffusion reaction equation.
\end{abstract}



\setcounter{equation}{0}
\setcounter{figure}{0}
\setcounter{table}{0}
\setcounter{theorem}{0}
\setcounter{lemma}{0}
\setcounter{corollary}{0}
\setcounter{definition}{0}
%
%%\section{Introduction}
%% \label{sec_intro}
% \input{RegFDeq_Intro}
\section{Introduction}
 \label{sec_intro}
Of interest in this report is the regularity of the solution of the fractional diffusion equation
\begin{align}
 \mcL_{r}^{\alpha}u(x) \ := \  
   -  \left( r D^{\alpha} \, + \, (1 - r) D^{\alpha *}  \right) u(x) &= \ f(x) \, , \  \ x \in \mrI \, , \
   \label{DefProb1}  \\
   \mbox{subject to } u(0) \, = \, u(1) &= \,  0 \, ,
\label{DefBC1}
\end{align}  
and  the regularity of the solution of the fractional diffusion, advection, reaction equation
\begin{align}
 \mcL_{r}^{\alpha}u(x) \ + \ b(x) D u(x) \ + \ c(x) u(x)   &= \ f(x) \, , \ \  x \in \mrI \, , \
   \label{DefProb2}  \\
   \mbox{subject to } u(0) \, = \, u(1) &= \,  0 \, ,
\label{DefBC2}
\end{align}  
where $\mrI \, := \, (0 , 1)$, $1 < \alpha < 2$, $0 \le r \le 1$, $c(x) \, - \, \frac{1}{2}D b(x) \ge 0$, 
$D$ denotes the usual derivative
operator, $D^{\alpha}$ the $\alpha$-order left fractional derivative operator, and
$D^{\alpha *}$ the $\alpha$-order right fractional derivative operator, defined by:
\begin{align}
D^{\alpha} u(x) &:= \ 
D \, \frac{1}{\Gamma(2 - \alpha)}   \int_{0}^{x} \frac{1}{(x - s)^{\alpha - 1}} \, D u(s) \, ds \, ,   \label{defDalpha}  \\
D^{\alpha *} u(x) &:= \ 
 D \,  \frac{1}{\Gamma(2 - \alpha)} \int_{x}^{1} \frac{1}{(s - x)^{\alpha - 1}} \, D u(s) \, ds \, .   \label{defDalpha*} 
\end{align}
The regularity of the solution to a differential equation plays a fundamental role in designing optimal approximation
schemes for the solution.

In recent years fractional order differential equations have received increased attention due to their application
in the modeling of physical phenomena such as  in contaminant transport in ground water flow \cite{ben001, cus931},
viscoelasticity \cite{mai101}, image processing \cite{ant191, bua101, gat151, gil081},
turbulent flow \cite{mai971, shl871}, and chaotic dynamics \cite{zas931}.

The diffusion operator, $\mcL_{r}^{\alpha}$, arises in a random walk process in which the jumps have an
unbounded variance (L\'{e}vy process) \cite{ben001, roo041}.

For the case $r \, = \, 1/2$,  $\mcL_{1/2}^{\alpha}$ represents the \textit{integral} fractional Laplacian 
operator \cite{aco171}. The existence, uniqueness and regularity of the solution to the fractional Laplacian equation
has been investigated by a number of authors, in $\mathbb{R}^{1}$ see \cite{aco181}, 
in $\mathbb{R}^{n \ge 2}$ \cite{alb151, coz171, gru151, ros141}. Recently the regularity results for the fractional Laplacian 
equation was extended by Hao and Zhang in \cite{hao182, zha191} to the fractional Laplacian equation with a
constant advection and reaction term (i.e. \eqref{DefProb2},\eqref{DefBC2}, for $r = 1/2$, 
$b(x) = b$, $c(x) = c$, $b, \, c \in \mathbb{R}$).

Fewer results on the regularity of the solution to the general fractional diffusion, advection, reaction equation
have been established. In \cite{erv061} Ervin and Roop established existence and uniqueness of solution,
$u \in H^{\alpha/2}_{0}(\mrI)$, for $f \in H^{- \alpha/2}(\mrI)$. More recently in \cite{erv162, jia181} precise 
regularity results were obtained for the solution of \eqref{DefProb1},\eqref{DefBC1} for $f \in H^{s}_{(a , b)}(\mrI)$,
where $H^{s}_{(a , b)}(\mrI)$ denotes an appropriated weighted Sobolev space.
In \cite{hao181} Hao, Guang and Zhang obtained regularity estimates for the solution of 
\eqref{DefProb2} for $b(x) = 0$, $c(x) = c$. Their
numerical experiments indicated that their regularity estimates were not optimal.

In this article we present the general regularity results for \eqref{DefProb2},\eqref{DefBC2}, in appropriately
weighted Sobolev spaces. 
The analysis establishes that the presence of a reaction term (i.e. $c(x) \neq 0)$ limits the 
regularity of the solution, regardless of the smoothness of the right hand side function, $f(x)$.
This reduction in regularity is greater (by a factor of 2) when an advective term (i.e. $b(x) \neq 0)$
appears in \eqref{DefProb2}. This behavior of the solution is in sharp contrast to that for the
integer order ($\alpha = 2$) diffusion, advection, reaction equation. In that case, assuming 
$b(x)$ and $c(x)$ are sufficiently regular, for the right hand side function $f \in H^{s}(\mrI)$
the solution lies in $H^{s + 2}(\mrI)$.

The results we present herein extend those in \cite{erv162} for the fractional diffusion equation,
and those in \cite{hao182} for the fractional Laplacian equation with a constant advection and reaction term. 
The proofs given are significantly different that those used in \cite{hao181, hao182}.

The analysis of \eqref{DefProb2},\eqref{DefBC2} is most appropriately performed in  weighted Sobolev spaces
(due to the singular behavior of the solution at the endpoints). There are different ways to define the weighted 
Sobolev spaces: (i) using interpolation (Section \ref{sec_Wspace1}), (ii) using an appropriate basis (Section \ref{sec_Wspace2}),
(iii) using an explicit definition for the fractional order norms (Section \ref{sec_Wspace3}). Each of these representations
have their advantages and are used in the analysis.

This paper is organized as follows. In the next section we present some preliminary definitions and results.
Sections \ref{sec_Wspace1}, \ref{sec_Wspace2}, and \ref{sec_Wspace3} present the three different (but equivalent) definitions
of the weighted Sobolev spaces. Section \ref{sec_Wspace2} also establishes some useful properties of the 
space $H^{s}_{(a , b)}(\mrI)$. For example, which $H^{s}_{(a , b)}(\mrI)$ space the function $f(x) = x^{\mu}$ lies in,
and for which $H^{s}_{(a , b)}(\mrI)$ we have the embedding $H^{s}_{(a , b)}(\mrI) \subset C^{k}(\mrI)$. The
regularity of the solution to the fractional diffusion problem \eqref{DefProb1},\eqref{DefBC1}, is discussed in 
Section \ref{secRFD}. The theorems needed for determining which $H^{t}_{(b , a)}(\mrI)$ space
$(1 - x)^{a} x^{b} \phi(x)$ lies in, when $\phi(x) \in H^{s}_{(a , b)}(\mrI)$, are established in Section \ref{sec_Wspace3}.
Section \ref{secRFDAR2} then discusses the regularity of solutions to \eqref{DefProb2},\eqref{DefBC2}.

In the last section we relate the regularity results obtained in weighted Sobolev spaces to
the usual (unweighted) Sobolev spaces.

%%%%%
%%%%% 
 \setcounter{equation}{0}
\setcounter{figure}{0}
\setcounter{table}{0}
\setcounter{theorem}{0}
\setcounter{lemma}{0}
\setcounter{corollary}{0}
\setcounter{definition}{0}
\section{Notation and Properties}
\label{sec_not}
Jacobi polynomial have an important connection with fractional order diffusion equations
\cite{aco181, erv162, mao181, mao161}. We briefly review their definition and some of their important
properties \cite{abr641, sze751}. 

\textbf{Usual Jacobi Polynomials, $P_{n}^{(a , b)}(t)$, on $(-1 \, , \, 1)$}.   \\    
\underline{Definition}: $ P_{n}^{(a , b)}(t) \ := \ 
\sum_{m = 0}^{n} \, p_{n , m} \, (t - 1)^{(n - m)} (t + 1)^{m}$, where
\begin{equation}
       p_{n , m} \ := \ \frac{1}{2^{n}} \, \left( \begin{array}{c}
                                                              n + a \\
                                                              m  \end{array} \right) \,
                                                    \left( \begin{array}{c}
                                                              n + b \\
                                                              n - m  \end{array} \right) \, .
\label{spm21}
\end{equation}
\underline{Orthogonality}:    
\begin{align}
 & \int_{-1}^{1} (1 - t)^{a} (1 + t)^{b} \, P_{j}^{(a , b)}(t) \, P_{k}^{(a , b)}(t)  \, dt 
 \ = \
   \left\{ \begin{array}{ll} 
   0 , & k \ne j  \\
   |\| P_{j}^{(a , b)} |\|^{2}
   \, , & k = j  
    \end{array} \right.  \, ,  \nonumber \\
& \quad \quad \mbox{where } \  \ |\| P_{j}^{(a , b)} |\| \ = \
 \left( \frac{2^{(a + b + 1)}}{(2j \, + \, a \, + \, b \, + 1)} 
   \frac{\Gamma(j + a + 1) \, \Gamma(j + b + 1)}{\Gamma(j + 1) \, \Gamma(j + a + b + 1)}
   \right)^{1/2} \, .
  \label{spm22}
\end{align}                                                    

In order to transform the domain of the family of Jacobi polynomials to $[0 , 1]$, let $t \rightarrow 2x - 1$ and 
introduce $G_{n}^{(a , b)}(x) \, = \, P_{n}^{(a , b)}( t(x) )$. From \eqref{spm22},
\begin{align}
 \int_{-1}^{1} (1 - t)^{a} (1 + t)^{b} \, P_{j}^{(a , b)}(t) \, P_{k}^{(a , b)}(t)  \, dt 
 &= \
 \int_{0}^{1} 2^{a} \, (1 - x)^{a} \, 2^{b} \, x^{b} \, P_{j}^{(a , b)}(2x - 1) \, 
 P_{k}^{(a , b)}(2x - 1)  \, 2 \,  dx
 \nonumber \\
  &= \
2^{a + b + 1} \int_{0}^{1}   (1 - x)^{a}  \, x^{b} \, G_{j}^{(a , b)}(x) \, G_{k}^{(a , b)}(x)  \,  dx
\nonumber \\
&= \
   \left\{ \begin{array}{ll} 
   0 , & k \ne j \, , \\
  2^{a + b + 1} \, |\| G_{j}^{(a , b)} |\|^{2}
   \, , & k = j  
   \, . \end{array} \right.    \nonumber \\
 \quad \quad \mbox{where } \  \ |\| G_{j}^{(a , b)} |\| &= \
 \left( \frac{1}{(2j \, + \, a \, + \, b \, + 1)} 
   \frac{\Gamma(j + a + 1) \, \Gamma(j + b + 1)}{\Gamma(j + 1) \, \Gamma(j + a + b + 1)}
   \right)^{1/2} \, .  \label{spm22g} 
\end{align}                                                    

\begin{equation}
\mbox{Note that } \quad  |\| G_{j}^{(a , b)} |\| \ = \ |\| G_{j}^{(b , a)} |\| \, .
\label{nmeqG}
\end{equation}
%and from \cite{abr641, sze751}
%\begin{equation}
%G_{j}^{(a , b)}(0) \ = \ (-1)^{j} \, \frac{\Gamma(j + b + 1)}{\Gamma(j + 1) \, \Gamma(b + 1)} \, .
%\label{eqG0}
%\end{equation}

 From \cite[equation (2.19)]{mao161} we have that
\begin{equation}
   \frac{d^{k}}{dt^{k}} P_{n}^{(a , b)}(t) \ = \ 
   \frac{\Gamma(n + k + a + b + 1)}{2^{k} \, \Gamma(n + a + b + 1)} P_{n - k}^{(a + k \, , \, b + k)}(t) \, .
   \label{derP}
\end{equation}   
Hence,
\begin{align}
\frac{d^{k}}{dx^{k}} G_{n}^{(a , b)}(x) 
  &= \ \frac{\Gamma(n + k + a + b + 1)}{  \Gamma(n + a + b + 1)} 
  G_{n - k}^{(a + k \, , \, b + k)}(x)  \, .  \label{eqC4}
\end{align}   

Also, from \cite[equation (2.15)]{mao161}, 
\begin{equation}
\frac{d^{k}}{dt^{k}} \left\{ (1 - t)^{a + k} \, (1 + t)^{b + k} \, P_{n - k}^{(a + k \, , \, b + k)}(t) \right\}
 \ = \ 
 \frac{(-1)^{k} \, 2^{k} \, n!}{(n - k)!} \, (1 - t)^{a} \, (1 + t)^{b} \, P_{n}^{(a \, , \, b)}(t) \, , \
 n \ge k \ge 0 \, ,
 \label{eqB0}
\end{equation}
from which it follows that
\begin{equation}
 \frac{d^{k}}{dx^{k}} \left\{  \ (1 \, - \, x)^{a + k} \,  x^{b + k} \,
 G_{n - k}^{(a + k \, , \, b + k)}(x) \right\} 
 \ = \ 
 \frac{(-1)^{k} \,  n!}{(n - k)!} \,  (1 \, - \, x)^{a} \,  x^{b} \,
 G_{n}^{(a \, , \, b)}(x) \, . 
  \label{eqC2}
 \end{equation}
 
For compactness of notation we introduce
\begin{equation}
 \rho^{(a , b)} \, = \, \rho^{(a , b)}(x) \, := \, (1 - x)^{a} \, x^{b} \, .
 \label{defrho}
\end{equation} 

 We let $\mathbb{N}_{0}  := \mathbb{N} \cup {0}$ and
 use $y_{n} \sim n^{p}$ to denote that there exists constants $c$ and $C > 0$ such that, as 
 $n \rightarrow \infty$,  
 $c \, n^{p} \le | y_{n} | \le C \, n^{p}$. Additionally, we use $a \, \lesssim \, b$ to denote that there exists a constant $C$ such that
 $a \, \le \, C \,  b$. 
 
 For $s \in \mathbb{R}$, $\floor{s}$ is used to denote the largest integer that is less than or equal to $s$, and
 $\ceil{s}$ is used to denote the smallest integer that is greater than or equal to $s$.
 
Note, 
from Stirling's formula we have that
\begin{equation}
\lim_{n \rightarrow \infty} \, \frac{\Gamma(n + \sigma)}{\Gamma(n) \, n^{\sigma}}
\ = \ 1 \, , \mbox{ for } \sigma \in \mathbb{R}.  
 \label{eqStrf}
\end{equation} 

\underline{Definition}: \textbf{Condition A} \\
For  $\alpha$ and  $r$ given, satisfying $1 < \alpha < 2$, $0 \le r \le 1$, let $\beta$ be determined by
$\alpha - 1 \, \le  \, \beta \,  \le  \, 1$ and
\begin{equation}
  r \ = \ \frac{\sin( \pi \, \beta)}{\sin( \pi ( \alpha - \beta)) \, + \,  \sin( \pi \, \beta)} \, . \label{propker0} 
\end{equation}
Furthermore, introduce the constant 
$c_{*}^{*}$ defined by
\begin{equation}
  c_{*}^{*} \ = \ \frac{\sin(\pi \alpha)}{\sin(\pi (\alpha - \beta)) \, + \, \sin(\pi \beta)} \, .  \label{defcss}
\end{equation}

\textbf{Function space $L_{\omega}^{2}(\mrI)$}. \\
For $\omega(x) > 0, \ x \in (0 , 1)$, let 
\begin{equation}
L_{\omega}^{2}(\mrI) \, := \, \{ f(x) \, : \, \int_{0}^{1} \omega(x) \, f(x)^{2} \, dx \ < \ \infty \} \, .
\label{defLw}
\end{equation}
Associated with $L_{\omega}^{2}(0 , 1)$ is the inner product, $\langle \cdot , \cdot \rangle_{\omega}$, and
norm, $\| \cdot \|_{\omega}$, defined by
\[
\langle f \,  , \, g \rangle_{\omega} \, := \, \int_{0}^{1} \omega(x) \, f(x) \, g(x) \, dx \, , \quad \mbox{and} \quad
 \| f \|_{\omega} \, := \, \left( \langle f \,  , \, f \rangle_{\omega} \right)^{1/2} \, .
\]

The set of orthogonal polynomials $\{ G_{j}^{(a , b)} \}_{j = 0}^{\infty}$ form an orthogonal basis
for $L^{2}_{\rho^{(a , b)}}(\mrI)$.

%%%%%
%%%%% 
 \setcounter{equation}{0}
\setcounter{figure}{0}
\setcounter{table}{0}
\setcounter{theorem}{0}
\setcounter{lemma}{0}
\setcounter{corollary}{0}
\setcounter{definition}{0}
\section{Weighted Sobolev Space defined by interpolation}
\label{sec_Wspace1}
Following Babu\v{s}ka and Guo \cite{bab011}, and Guo and Wang \cite{guo041}, $n \in \mathbb{N}_{0}$ define the
weighted Sobolev spaces 
\begin{align} 
H^n_{\rho^{(a,b)} }(\mrI) &:= \bigg \{ v \, : 
 \sum_{j=0}^n \big \| D^j v \big \|_{\rho^{(a+j,b+j)}}^2 < \infty \bigg \} ,   \label{defHw} 
 \end{align} 
with associated norm $ \| v \|_{n , \rho^{(a,b)}}  := 
 \left( \sum_{j=0}^n \big \| D^j v \big \|_{\rho^{(a+j,b+j)}}^2  \right)^{1/2}$. 

Definition (\ref{defHw}) is extended to $s \in \mathbb{R}^{+}$ using the $K$- method of interpolation. For $s < 0$
the spaces are defined by (weighted) $L^{2}$ duality.

Babu\v{s}ka and Guo used the  $H^s_{\rho^{(a,b)} }(\mrI)$ spaces in establishing the optimal convergence properties
of the $p$- version of the finite element method. They also related the definition \eqref{defHw} to the decay property
of the coefficients of the Jacobi polynomials of the expansion of the function $v$.

In \cite{guo041} Guo and Wang derived approximation properties of Jacobi polynomials for functions in the weighted
Sobolev spaces \eqref{defHw}.

%%%%%
%%%%% 
 \setcounter{equation}{0}
\setcounter{figure}{0}
\setcounter{table}{0}
\setcounter{theorem}{0}
\setcounter{lemma}{0}
\setcounter{corollary}{0}
\setcounter{definition}{0}
\section{Weighted Sobolev Space defined by Jacobi coefficients}
\label{sec_Wspace2}
In this section we define function spaces in terms of 
the decay property of the Jacobi coefficients of their member functions. We then show that these spaces 
agree with the weighted Sobolev spaces defined using the $K$- method of interpolation.

This presentation parallels the work of Acosta, Borthagaray, Bruno and Maas in \cite{aco181}  who investigated
the regularity and approximation of the 1-d fractional Laplacian equation. The 1-d fractional Laplacian operator
they considered
is a special case of the fractional diffusion operator, $\mcL_{r}^{\alpha}(\cdot)$, for the parameter $r = 1/2$. 
In \cite{aco181} their analysis
focused on the coefficients of the Gegenbauer polynomials, which are a special case of the Jacobi polynomials
where the weight parameters are equal.

We begin by relating the decay rate of a function's Jacobi polynomial coefficient to its regularity.
We then introduce the weighted Sobolev space $H^{s}_{(a , b)}(\mrI)$ and show that it corresponds to
the space $H^{s}_{\rho^{(a , b)}}(\mrI)$. 
We conclude the section with a corollary describing the precise solvability of \eqref{DefProb1},\eqref{DefBC1},
an embedding theorem relating 
$H^{s}_{(a , b)}(\mrI)$ to $C^{k}(\mrI)$, and a 
corollary stating a suitable sufficient condition for the solution of \eqref{DefProb1},\eqref{DefBC1} to be 
continuous.

Introduce $\wtilde{G}_{j}^{(a , b)}(x) \, = \, G_{j}^{(a , b)}(x) / |\| G_{j}^{(a , b)} |\| $ as an
orthonormal basis for $L^{2}_{\rho^{(a , b)}}(\mrI)$.

Then, given $v \in L^{2}_{\rho^{(a , b)}}(\mrI)$, we have
\be
   v(x) \ = \ \sum_{j = 0}^{\infty} v_{j} \, \wtilde{G}_{j}^{(a , b)}(x) \, , 
 \label{vexp1}
\ee 
which converges in $L^{2}_{\rho^{(a , b)}}(\mrI)$, where
\be
  v_{j} \ = \ \int_{0}^{1}  \rho^{(a , b)}(x) \, v(x) \,  \wtilde{G}_{j}^{(a , b)}(x) \, dx \, .
 \label{defvj}
\ee

Using \eqref{eqC4} and \eqref{vexp1} we could conjecture (by differentiating \eqref{vexp1}) that for sufficiently smooth $v(x)$ 
\begin{align}
  v^{(k)}(x) &= \ \sum_{j = k}^{\infty} \frac{|\| G_{j - k}^{(a + k \, , b + k)} |\|}{|\| G_{j}^{(a , b)} |\|}
  \frac{\Gamma(j + k + a + b + 1)}{\Gamma(j + a + b + 1)} \, v_{j} \, 
   \wtilde{G}_{j - k}^{(a + k \, , b + k)}(x)   \nonumber \\
   &:= \  \, \sum_{j = k}^{\infty}  \, v_{j - k}^{(k)} \, 
   \wtilde{G}_{j - k}^{(a + k \, , b + k)}(x)  ,  \label{defJvk}
\end{align}
where
\be
   v_{j - k}^{(k)} \ = \  \frac{|\| G_{j - k}^{(a + k \, , b + k)} |\|}{|\| G_{j}^{(a , b)} |\|}
  \frac{\Gamma(j + k + a + b + 1)}{\Gamma(j + a + b + 1)} \, v_{j} \, .
  \label{defvjmk}
\ee  

\textbf{Remark}: The terms $v_{j}$ and $v_{j - k}^{(k)}$ denote real numbers, whereas $v^{(k)}(x)$ is used
to represent the $k^{th}$ derivative of $v(x)$. (From \eqref{defJvk}, $v_{j - k}^{(k)}$'s are the Jacobi coefficients
for $v^{(k)}(x)$.)

We investigate \eqref{defvjmk} more rigorously using the following lemma.

\begin{lemma} \label{lmaibp} (See \cite[Lemma 4.2]{aco181} Integration by parts) \\
Let $k \in \Natural$, $k \ge 2$, and let $v \in C^{k - 2}[0 , 1]$ such that for the decomposition of 
$[0 , 1] \, = \, \cup_{i = 1}^{n} [x_{i} , x_{i + 1}]$ $(0 = x_{1} < x_{2} < \ldots < x_{n+1} = 1)$, and for 
functions $\widehat{v}_{i} \in C^{k}[x_{i} , x_{i + 1}]$, we have $v(x) \, = \, \widehat{v}_{i}(x)$ for 
$x \in (x_{i} , x_{i + 1})$, and $1 \le i \le n$. Then for $j \ge k$ the $(a , b)$ Jacobi coefficient
defined by \eqref{defvj} satisfies
\begin{align}
v_{j} &=
 B_{j}^{k}  \, 
   \int_{0}^{1} \,  \rho^{(a + k \, , \, b + k)}(x) \, \wtilde{G}_{j - k}^{(a + k \, , \, b + k)}(x) \, v^{(k)}(x) \, dx  \nonumber \\
& \quad \quad \quad \quad - \     B_{j}^{k}  \, \sum_{l = 1}^{n} \, 
   \rho^{(a + k \, , \, b + k)}(x) \, \wtilde{G}_{j - k}^{(a + k \, , \, b + k)}(x) \, v_{l}^{(k - 1)}(x) \bigg|_{x = x_{l}}^{x_{l + 1}}  \, ,  \label{form43}  \\
\mbox{where } B_{j}^{k} &= \    \frac{\Gamma(j - k + 1)}{\Gamma(j + 1)} 
  \frac{|\| G_{j - k}^{(a + k , b + k)} |\|}{|\| G_{j}^{(a , b)} |\|}  \, .  \label{defBjk}
\end{align}
\end{lemma}
\textbf{Proof}:  From \eqref{eqC2}, 
\[
   \int \, \rho^{(a , b)}(x) \,  G_{n}^{(a , b)}(x) \, dx \ = \ 
   \frac{- (n - 1) !}{n !} \rho^{(a + 1 \, , \, b + 1)}(x) \, G_{n - 1}^{(a + 1 \, , \, b + 1)}(x) 
   \ = \ \frac{- 1}{n} \rho^{(a + 1 \, , \, b + 1)}(x) \, G_{n - 1}^{(a + 1 \, , \, b + 1)}(x)  \, .
\]
Hence,
\[
   \int \, \rho^{(a , b)}(x) \, \wtilde{G}_{n}^{(a , b)}(x) \, dx \ = \ 
   \frac{- 1}{n}  \, \frac{|\| G_{n - 1}^{(a + 1 , b + 1)} |\|}{|\| G_{n}^{(a , b)} |\|} 
   \rho^{(a + 1 \, , \, b + 1)}(x) \, \wtilde{G}_{n - 1}^{(a + 1 \, , \, b + 1)}(x)  \, .
\]

Beginning with \eqref{defvj} and using integration by parts,
\begin{align}
v_{j} &= \ \sum_{l = 1}^{n} \int_{x_{l}}^{x_{l + 1}} \widehat{v}_{l}(x) \,  
\rho^{(a , b)}(x) \, \wtilde{G}_{j}^{(a , b)}(x) \, dx \nonumber \\
&= \ \sum_{l = 1}^{n} \int_{x_{l}}^{x_{l + 1}}   \frac{1}{j}  \frac{|\| G_{j - 1}^{(a + 1 , b + 1)} |\|}{|\| G_{j}^{(a , b)} |\|} 
   \rho^{(a + 1 \, , \, b + 1)}(x) \, \wtilde{G}_{j - 1}^{(a + 1 \, , \, b + 1)}(x) \, \widehat{v}_{l}^{(1)}(x) \, dx  \nonumber \\
& \quad \quad \quad \quad - \     \frac{1}{j}  \frac{|\| G_{j - 1}^{(a + 1 , b + 1)} |\|}{|\| G_{j}^{(a , b)} |\|} 
   \rho^{(a + 1 \, , \, b + 1)}(x) \, \wtilde{G}_{j - 1}^{(a + 1 \, , \, b + 1)}(x) \, 
   \widehat{v}_{l}(x) \bigg|_{x \, = \, x_{l}}^{x_{l + 1}}  \nonumber \\
&= \ \sum_{l = 1}^{n} \int_{x_{l}}^{x_{l + 1}}   \frac{1}{j}  \frac{|\| G_{j - 1}^{(a + 1 , b + 1)} |\|}{|\| G_{j}^{(a , b)} |\|} 
   \rho^{(a + 1 \, , \, b + 1)}(x) \, \wtilde{G}_{j - 1}^{(a + 1 \, , \, b + 1)}(x) \, \widehat{v}_{l}^{(1)}(x) \, dx  \, . 
   \label{ibpeq1} 
\end{align}   
Repeated use of integration by parts, and using that
\[
\frac{1}{j}   \frac{1}{j - 1}  \ldots \frac{1}{j - k + 1}  \ = \ \frac{\Gamma(j - k + 1)}{\Gamma(j + 1)} \, , 
\]
\eqref{ibpeq1} becomes
\begin{align}
v_{j} &= \ \sum_{l = 1}^{n} \int_{x_{l}}^{x_{l + 1}}  \frac{\Gamma(j - k + 1)}{\Gamma(j + 1)} 
  \frac{|\| G_{j - k}^{(a + k , b + k)} |\|}{|\| G_{j}^{(a , b)} |\|} 
   \rho^{(a + k \, , \, b + k)}(x) \, \wtilde{G}_{j - k}^{(a + k \, , \, b + k)}(x) \, \widehat{v}_{l}^{(k)}(x) \, dx  \nonumber \\
& \quad \quad \quad \quad - \     \frac{\Gamma(j - k + 1)}{\Gamma(j + 1)} 
  \frac{|\| G_{j - k}^{(a + k , b + k)} |\|}{|\| G_{j}^{(a , b)} |\|} 
   \rho^{(a + k \, , \, b + k)}(x) \, \wtilde{G}_{j - k}^{(a + k \, , \, b + k)}(x) \, 
   \widehat{v}_{l}^{(k - 1)}(x) \bigg|_{x \, = \, x_{l}}^{x_{l + 1}}  \nonumber \\
 &= \  \frac{\Gamma(j - k + 1)}{\Gamma(j + 1)} 
  \frac{|\| G_{j - k}^{(a + k , b + k)} |\|}{|\| G_{j}^{(a , b)} |\|}  \, 
   \int_{0}^{1} \,  \rho^{(a + k \, , \, b + k)}(x) \, \wtilde{G}_{j - k}^{(a + k \, , \, b + k)}(x) \, v^{(k)}(x) \, dx  \nonumber \\
& \quad \quad \quad \quad - \     \frac{\Gamma(j - k + 1)}{\Gamma(j + 1)} 
  \frac{|\| G_{j - k}^{(a + k , b + k)} |\|}{|\| G_{j}^{(a , b)} |\|}  \, \sum_{l = 1}^{n} \, 
   \rho^{(a + k \, , \, b + k)}(x) \, \wtilde{G}_{j - k}^{(a + k \, , \, b + k)}(x) \, 
   \widehat{v}_{l}^{(k - 1)}(x) \bigg|_{x \, = \, x_{l}}^{x_{l + 1}}     \, , 
  \label{ibpeq2}
\end{align}
which corresponds with \eqref{form43}.
\mbox{ } \hfill \qed

Note, for $v \in C^{k}[0 , 1]$ we have that the $j^{th}$ Jacobi polynomial coefficient $v_{j}^{(k)}$ is given by
\[
  v_{j}^{(k)} \ = \ \int_{0}^{1}  \rho^{(a + k \,  , \,  b + k)} \, v^{(k)}(x) \,  \wtilde{G}_{j}^{(a + k \,  ,  \, b + k)}(x) \, dx \, ,
\]
and from \eqref{ibpeq2}
\be
 v_{j - k}^{(k)} \ = \  \frac{\Gamma(j + 1)}{\Gamma(j - k + 1)}
  \frac{|\| G_{j}^{(a , b)} |\|}{|\| G_{j - k}^{(a + k , b + k)} |\|} \, v_{j} \, , \ \ \mbox{ for } j \ge k \, .
  \label{vjmk2}
\ee 

Using \eqref{spm22}, a simple calculation establishes the equivalence of \eqref{vjmk2} and \eqref{defvjmk}. 

In order to investigate how the decay rate of the Jacobi coefficients relate to the regularity of the function we establish
the following lemma.
\begin{lemma} \label{lmaalog43}
There exists constants $C_{1}$ and $C_{2}$ such that
\be
C_{1} j^{-k} \le \frac{\Gamma(j - k + 1)}{\Gamma(j + 1)} 
  \frac{|\| G_{j - k}^{(a + k , b + k)} |\|}{|\| G_{j}^{(a , b)} |\|} \le C_{2} j^{-k} \, .
  \label{eqlma43}
\ee
\end{lemma}
\textbf{Proof}: We have that
\be
\frac{\Gamma(j - k + 1)}{\Gamma(j + 1)}  \ = \ \frac{\Gamma(j - k + 1)}{j (j - 1) \ldots (j - k + 1) \Gamma(j - k + 1)}   
\ = \ \frac{1}{j (j - 1) \ldots (j - k + 1)} \ \sim \ j^{-k} \,  .  
 \label{ewqr1}
\ee
Next consider
\begin{align}
 \frac{|\| G_{j - k}^{(a + k , b + k)} |\|^{2}}{|\| G_{j}^{(a , b)} |\|^{2}}
 &= \ 
 \frac{(2j \, + \, a \, + \, b \, + \, 1)}{1} \, 
 \frac{\Gamma(j + 1) \, \Gamma(j + a + b + 1)}{\Gamma(j + a + 1) \, \Gamma(j + b + 1)}  \nonumber \\
& \ \ 
 \cdot 
 \frac{1}{(2j \, - \, 2k \, + \, a \, + \, k \, + \, b \, + \, k \, +  \, 1)} \, 
 \frac{\Gamma(j - k + a + k + 1) \, \Gamma(j - k + b + k + 1)}%
 {\Gamma(j  - k  + 1) \, \Gamma(j -  k \, +  a +   k  +  b +  k   +  1)}  \nonumber \\
&= \ 
  \frac{\Gamma(j + 1) \, \Gamma(j + a + b + 1)}{\Gamma(j + 1 - k) \, \Gamma(j + a + b + 1 + k)}  \nonumber \\
&\sim \   (j + 1)^{k} \, (j + a + b + 1)^{-k} \ \sim \ \bigg( \frac{j + 1}{j + 1 + a + b} \bigg)^{k}  \nonumber \\
&\sim \ 1 \, . \label{ewqr2}
\end{align}
From \eqref{ewqr1} and \eqref{ewqr2}, \eqref{eqlma43} follows. \\
\mbox{ } \hfill \qed.

\begin{corollary} (See \cite[Corollary 4.3]{aco181}) \label{cor43}
Let $k \in \mathbb{N}$ and $v$ satisfy the hypothesis of Lemma \ref{lmaibp}. Then the Jacobi
coefficients in \eqref{form43} are quantities of order $O(j^{- k})$ as $j \rightarrow \infty$, i.e., 
\be
       | v_{j} | \ < \ C \, j^{- k} \, ,
\label{vjest5}
\ee
for a constant $C$ that depends on $v$ and $k$.
\end{corollary}
\textbf{Remark}: In words, this corollary states that if a function is $C^{k}(\mrI)$, with the possible exception
of a few points, then its $j^{th}$ Jacobi coefficient must decay as $C \, j^{-k}$, for some constant $C$.

\textbf{Proof}: 
In order to establish the estimate, for the first term in \eqref{form43} note that
\begin{align}
& \bigg| \int_{0}^{1} \rho^{(a + k \, , \, b + k)}(x) \, \wtilde{G}_{j - k}^{(a + k \, , \, b + k)}(x) \, v^{( k )}(x) \, dx \, \bigg| \nonumber \\
& \quad \quad \quad \le \ 
\bigg(   \int_{0}^{1} \rho^{(a + k \, , \, b + k)}(x) \, ( \wtilde{G}_{j - k}^{(a + k \, , \, b + k)}(x) )^{2} \, dx \, \bigg)^{1/2} 
\, \bigg(   \int_{0}^{1} \rho^{(a + k \, , \, b + k)}(x) \, ( v^{( k )}(x) )^{2} \, dx \, \bigg)^{1/2}  \nonumber \\
& \quad \quad \quad \le \ 1 \, 
    \| v^{(k)} \|_{L^{\infty}}  \bigg(   \int_{0}^{1} \rho^{(a + k \, , \, b + k)}(x) \, dx \, \bigg)^{1/2}  \nonumber \\
& \quad \quad \quad = \ C \, \| v^{(k)} \|_{L^{\infty}} \, .   \label{bdfh1}
\end{align}
 
To bound the second term in \eqref{form43} consider Darboux's formula \cite[pg. 196 (8.21.10)]{sze751}
\begin{align}
 P_{n}^{(a , b)}(\cos \theta) &=
 \ n^{-1/2} \, k(\theta) \, \cos( N \theta \, + \, \gamma) \ + \ O(n^{-3/2}) \, ,  \label{Dbx1}  \\
 \mbox{where } \ 
 k(\theta) &=  \ \pi^{-1/2} \, \big( \sin( \theta / 2) \big)^{-a - 1/2} \, \big( \cos( \theta / 2) \big)^{-b - 1/2} \, , 
 \ \ N \ = \ n \, + \, (a + b + 1)/2 \, ,   \label{Dbx2}  \\
 \gamma &= \ - (a + 1/2) \pi/2 \, , \ \ 0 < \theta < \pi \, .  \nonumber 
\end{align}
The bound for the error term holds uniformly in the interval $[\epsilon \, , \, \pi - \epsilon]$.
 
Using the trig. identities
\[ 
\cos(\theta/2) \ = \ \frac{1}{\sqrt{2}} \, ( 1 \, + \, \cos(\theta) )^{1/2} \, , \ \ \mbox{ and } \ 
\sin(\theta/2) \ = \ \frac{1}{\sqrt{2}} \, ( 1 \, - \, \cos(\theta) )^{1/2}  \, ,
\]
and the substitution $x = \cos(\theta)$, for $k(\theta)$ we have
\begin{align}
k(x) &= \ \pi^{-1/2} \, 2^{-\frac{1}{2} ( -a - 1/2 )} \, (1 \, - \, x)^{\frac{1}{2}(-a - 1/2)} \, 
2^{-\frac{1}{2} ( -b - 1/2 )} \, (1 \, + \, x)^{\frac{1}{2}(-b - 1/2)}  \nonumber \\
&= \ \pi^{-1/2} \, 2^{\frac{1}{2} ( a + b + 1)} \, (1 \, - \, x)^{-\frac{1}{2}(a + 1/2)} \, 
(1 \, + \, x)^{-\frac{1}{2}(b + 1/2)} \, . \nonumber
\end{align}

Therefore, for $P_{n}^{(a , b)}(x)$ we have
\be
P_{n}^{(a , b)}(x) \ = \ n^{-1/2} 
\bigg(  \pi^{-1/2} \, 2^{\frac{1}{2} ( a + b + 1)} \, (1 \, - \, x)^{-\frac{1}{2}(a + 1/2)} \, 
(1 \, + \, x)^{-\frac{1}{2}(b + 1/2)} \bigg) \, 
 \cos( N Arccos(x) \, + \, \gamma) \ + \ O(n^{-3/2}) \, .
\label{Pnrep}
\ee
Using \eqref{Pnrep}, we have for ($x \rightarrow 2t - 1$) $G_{n}^{(a , b)}(t)$ , 
\[
G_{n}^{(a , b)}(t) \ = \ n^{-1/2} 
\bigg(  \pi^{-1/2} \,  (1 \, - \, t)^{-\frac{1}{2}(a + 1/2)} \, 
t^{-\frac{1}{2}(b + 1/2)} \bigg) \, 
 \cos( N Arccos(2t - 1) \, + \, \gamma) \ + \ O(n^{-3/2}) \, .
\]

Now,
\begin{align}
\rho^{(a + k \, , \, b + k)}(x) \, G_{n}^{(a + k \, , \, b + k)}(x)
\ =   \quad \quad \quad   \quad \quad \quad & \nonumber \\
\ \ n^{-1/2} 
\bigg(  \pi^{-1/2} \,  (1 \, - \, x)^{\frac{1}{2}(a + k -1)} \, 
x^{\frac{1}{2}(b + k -1)} \bigg) \, &
 \cos( N Arccos(2x - 1) \, + \, \gamma) \ + \ O(n^{-3/2}) \,  \nonumber \\
\mbox{i.e., } \ \big| \rho^{(a + k \, , \, b + k)}(x) \, G_{n}^{(a + k \, , \, b + k)}(x) \big|
&\le \ C \, n^{-1/2} \, ,  \label{Gnrep2}
\end{align}
as $\bigg(  \pi^{-1/2} \,  (1 \, - \, x)^{\frac{1}{2}(a + k -1)} \, 
x^{\frac{1}{2}(b + k -1)} \bigg) \, 
 \cos( N Arccos(2x - 1) \, + \, \gamma)$ is bounded for $x \in [x_{2} , x_{n}]$ and on this subinterval the error estimate from 
Darboux's formula is uniform. (Note that we are not interested in the boundary pieces in the Integration by Parts formula at
$x_{1}$ and $x_{n+1}$, as at these points the boundary pieces are zero because of the 
$\rho^{(a + k \, , \, b + k)}(x)$
term). 

From  \eqref{spm22g} and \eqref{eqStrf} we have
\begin{align}
 |\| G_{j}^{(a , b)} |\|^{2} &= \
 \frac{1}{(2j \, + \, a \, + \, b \, + 1)} 
   \frac{\Gamma(j + a + 1) \, \Gamma(j + b + 1)}{\Gamma(j + 1) \, \Gamma(j + a + b + 1)}  \nonumber \\
&\sim \ j^{-1} \, (j + 1)^{a} \, (j + b + 1)^{-a} \ \sim \    j^{-1} \, .   \label{nrmGbeh}
\end{align}
 
Combining \eqref{Gnrep2} and \eqref{nrmGbeh} we obtain
\be
\big| \rho^{(a + k \, , \, b + k)}(x) \, \wtilde{G}_{n}^{(a + k \, , \, b + k)}(a_{l}) \big|
\ \le \ C \,  \mbox{ for } l = 1, 2, \ldots, n+1.
\label{bdyibp3}
\ee

Using \eqref{bdfh1} and \eqref{bdyibp3}, together with Lemma \ref{lmaalog43} estimate \eqref{vjest5} follows. \\
\mbox{ } \hfill \qed

Next we introduce the $(a , b)$-weighted Sobolev spaces. 
\begin{definition} \label{Defspace}
Let $s, a, b \in \real$, $s \ge 0$, $a, b > -1$, 
$L^{2}_{(a , b)}(\mrI) \, := \, L^{2}_{\rho^{(a , b)}}(\mrI)$, and
$v_{j}$ be given by \eqref{defvj}. Then, we define 
\be
H^{s}_{(a , b)}(\mrI) \, := \, \{v \in L^{2}_{\rho^{(a , b)}}(\mrI) \, : \,
 \sum_{j = 0}^{\infty} (1 + j^{2})^{s} \, v_{j}^{2} \, < \, \infty \}
\label{defHr}
\ee
as the $(a , b)$-weighted Sobolev space of order $s$.
\end{definition}

\begin{lemma}(See Lemma 4.7 in \cite{aco181}) \label{lmaHrHspace}  
Let $s, a, b \in \real$, $s \ge 0$, $a, b > -1$. Then the space $H^{s}_{(a , b)}(\mrI)$ endowed with the
inner product 
\be
\langle v \, , \, w \rangle^{s}_{(a , b)} \, := \, \sum_{j = 0}^{\infty} (1 + j^{2})^{s} \, v_{j} \, w_{j} \, ,
\mbox{  and associated norm  }  \| v \|_{s , (a , b)} \, := \, \bigg( \sum_{j = 0}^{\infty} (1 + j^{2})^{s} \, v_{j}^{2} \bigg)^{1/2}
\label{ipandnmdef}
\ee
is a Hilbert space.
\end{lemma}
\textbf{Proof}: The proof is completely analogously to that given in \cite[Theorem 8.2]{kre891}. \\
\mbox{ } \hfill \qed

\textbf{Remarks}. \\
1. It follows immediately from the definition that for any function $v \in H^{s}_{(a , b)}(\mrI)$ its expansion \eqref{vexp1}
is convergent in $H^{s}_{(a , b)}(\mrI)$.

2. From Parseval's identity, it follows that $H^{0}_{(a , b)}(\mrI)$
and $L^{2}_{\rho^{(a , b)}}(\mrI)$ coincide. Further, we have the dense compact embedding
$H^{t}_{(a , b)}(\mrI) \subset H^{s}_{(a , b)}(\mrI)$ whenever $t< s$. (The density of the 
embedding follows 1. and that all polynomials are contained in $H^{s}_{(a , b)}(\mrI)$ for every $s$. 
The compactness follows as in  \cite[Theorem 8.3]{kre891}).

\begin{definition} \label{DefDspace}
Let $s, a, b \in \real$, $s > 0$, $a, b > -1$. We denote by $H^{-s}_{(a , b)}(\mrI)$
the space of bounded linear functionals on $H^{s}_{(a , b)}(\mrI)$. 
\end{definition}

\begin{lemma} \label{eqDefDsp}
For $s > 0$, $a, b > -1$ the space $H^{-s}_{(a , b)}(\mrI)$ can be equivalently characterized as the set
\be
\mcH^{-s}_{(a , b)}(\mrI) \ := \ \left\{ f \, : \, f(x) \, = \, \sum_{j = 0}^{\infty} f_{j} \, \wtilde{G}_{j}^{(a , b)}(x) \, , 
\ \mbox{ where } \ \sum_{j = 0}^{\infty} (1 + j^{2})^{-s} \, f_{j}^{2} \, < \, \infty \right\} \, . 
\label{altdef1}
\ee
Additionally, for $F \in H^{-s}_{(a , b)}(\mrI)$, with representation 
$F  \, = \, \sum_{j = 0}^{\infty} f_{j} \, \wtilde{G}_{j}^{(a , b)}(x) \ := \ f(x)$ 
\be
\| F \|_{H^{-s}_{(a , b)}} \ = \ \sup_{g \in H^{s}_{(a , b)}(\mrI)} 
\frac{ | F(g) |}{ \| g \|_{s , (a , b)}} 
\ = \ \left( \sum_{j = 0}^{\infty} (1 + j^{2})^{-s} \, f_{j}^{2} \right)^{1/2}  \ := \ \| f \|_{-s , (a , b)} \, .
 \label{fdrmrm}
\ee
\end{lemma}
\textbf{Proof}: To establish the characterization we demonstrate that there is an isometry between 
$H^{-s}_{(a , b)}(\mrI)$ and $\mcH^{-s}_{(a , b)}(\mrI)$.

Let $g(x) \, = \, \sum_{j = 0}^{\infty} g_{j} \, \wtilde{G}_{j}^{(a , b)}(x) \in H^{s}_{(a , b)}(\mrI)$.

Consider $F(\cdot) \in H^{-s}_{(a , b)}(\mrI)$. Then, using the Riesz Representation Theorem,
there exists a unique $h(x) \, = \, \sum_{j = 0}^{\infty} h_{j} \, \wtilde{G}_{j}^{(a , b)}(x) \ \in H^{s}_{(a , b)}(\mrI)$, 
such that 
\begin{align}
F(g) &= \ \langle h \, ,  \, g \rangle_{(a , b)}^{s} \ = \ \sum_{j = 0}^{\infty}  (1 + j^{2})^{s} h_{j} \, g_{j}   \nonumber \\
&= \ \langle f \, ,  \, g \rangle_{(a , b)}^{0} \, ,   \nonumber
\end{align}
where $f(x) \ = \ \sum_{j = 0}^{\infty} f_{j} \, \wtilde{G}_{j}^{(a , b)}(x)$, for $f_{j} \, = \, (1 + j^{2})^{s} h_{j}$.

Note that as
\begin{align*}
 \sum_{j = 0}^{\infty}  (1 + j^{2})^{-s} f_{j}^{2} &= \  \sum_{j = 0}^{\infty} (1 + j^{2})^{-s} \, (1 + j^{2})^{2 s} h_{j}^{2}   \\
 &= \  \sum_{j = 0}^{\infty} (1 + j^{2})^{s} h_{j}^{2}   \ = \ \| h \|_{s , (a , b)}^{2} \, < \infty \, ,
\end{align*}
then $f \in \mcH^{-s}_{(a , b)}(\mrI)$, and from \eqref{fdrmrm}, $\| f \|_{-s , (a , b)} \, = \, \| h \|_{s , (a , b)}$.

Next, 
\begin{align*}
 \| F \|_{H^{-s}_{(a , b)}} &= \ \sup_{g \in H^{s}_{(a , b)}(\mrI)} 
\frac{ | F(g) |}{ \| g \|_{s , (a , b)}}  
\ = \ \sup_{g \in H^{s}_{(a , b)}(\mrI)}  \frac{ | \langle f \, ,  \, g \rangle_{(a , b)}^{0} |}{ \| g \|_{s , (a , b)}}   \\
&=  \ \sup_{g \in H^{s}_{(a , b)}(\mrI)}  \frac{ | \sum_{j = 0}^{\infty} f_{j} \, g_{j}  |}{ \| g \|_{s , (a , b)}} 
\ = \  \sup_{g \in H^{s}_{(a , b)}(\mrI)}  
\frac{ | \sum_{j = 0}^{\infty} (1 + j^{2})^{-s/2} f_{j} \, (1 + j^{2})^{s/2}  g_{j}  |}{ \| g \|_{s , (a , b)}}   \\
&\le \ \sup_{g \in H^{s}_{(a , b)}(\mrI)}  
\frac{ \left( \sum_{j = 0}^{\infty} (1 + j^{2})^{-s} f_{j}^{2} \right)^{1/2} \, 
\left( \sum_{j = 0}^{\infty} (1 + j^{2})^{s} g_{j}^{2} \right)^{1/2}}{ \| g \|_{s , (a , b)}}  
\ = \  \sup_{g \in H^{s}_{(a , b)}(\mrI)}  
\frac{ \| f \|_{-s , (a , b)} \, \| g \|_{s , (a , b)}}{ \| g \|_{s , (a , b)}}  \\
&= \ \| f \|_{-s , (a , b)} \, .
\end{align*}
Additionally, as
\[
\frac{ | F(h) |}{ \| h \|_{s , (a , b)}}  
\ = \  \frac{ | \langle f \, ,  \, h \rangle_{(a , b)}^{0} |}{ \| h \|_{s , (a , b)}}   
\ = \  \frac{ |  \sum_{j = 0}^{\infty} f_{j} \,  h_{j} |}{ \| f \|_{-s , (a , b)}}   
\ = \  \frac{ |  \sum_{j = 0}^{\infty} f_{j} \,  (1 + j^{2})^{-s} f_{j} |}{ \| f \|_{-s , (a , b)}}   
\ = \ \| f \|_{-s , (a , b)} \, ,
\]
then it follows that $ \| F \|_{H^{-s}_{(a , b)}} \ = \ \| f \|_{-s , (a , b)}$.

Finally, for $f \in \mcH^{-s}_{(a , b)}(\mrI)$ and any $g \in H^{s}_{(a , b)}(\mrI)$,
\[
\langle f \, ,  \, g \rangle_{(a , b)}^{0} \ = \  \sum_{j = 0}^{\infty} f_{j} \,  g_{j} 
\ = \  \sum_{j = 0}^{\infty} (1 + j^{2})^{-s/2} f_{j} \,  (1 + j^{2})^{s/2} g_{j} 
\ \le \ \| f \|_{-s , (a , b)} \, \| g \|_{s , (a , b)} \, .
\]
Thus, $f$ defines  a bounded  linear functional on $H^{s}_{(a , b)}(\mrI)$, i.e., $f \in H^{-s}_{(a , b)}(\mrI)$. \\
\mbox{ } \hfill \qed

With the structure of the $H^{s}_{(a , b)}(\mrI)$ spaces it is straight forward to show that $D$ is a bounded
mapping from $H^{s}_{(a , b)}(\mrI)$ onto
$H^{s - 1}_{(a + 1 \, , \, b + 1)}(\mrI)$.

\begin{lemma} \label{lmamapD}
For $s, a, b \in \mathbb{R}$, $a, b > -1$.the differential operator $D$ is a bounded mapping from $H^{s}_{(a , b)}(\mrI)$ onto
$H^{s - 1}_{(a + 1 \, , \, b + 1)}(\mrI)$.
\end{lemma}
\textbf{Proof}: Let $f \in H^{s}_{(a , b)}(\mrI)$. Then, $f(x)$ is expressible as
\[
   f(x) \ = \ \sum_{j = 0}^{\infty} c_{j} \,  \wtilde{G}_{j}^{(a , b)}(x) \, \ 
   \mbox{ where } \   \sum_{j = 0}^{\infty} (1 + j^{2})^{s} \, c_{j}^{2} < \infty \, .
\]
Using \eqref{eqC4},
\begin{align}
D f(x) &= \   \sum_{j = 1}^{\infty} c_{j} \,  \frac{1}{|\| \wtilde{G}_{j}^{(a , b)} |\|}
\, \frac{\Gamma(j + a + b + 2)}{\Gamma(j + a + b + 1)} \, G_{j - 1}^{(a + 1 \, , \, b + 1)}(x)  \nonumber  \\
&= \ \sum_{j = 1}^{\infty} c_{j} \,  \frac{|\| G_{j - 1}^{(a + 1 \, , \, b + 1)} |\|}{|\| G_{j}^{(a , b)} |\|} \, 
(j + a + b + 1) \, 
\wtilde{G}_{j - 1}^{(a + 1 \, , \, b + 1)}(x)  \, .  \label{donto1}
\end{align}
From \eqref{ewqr2} with $k = 1$, we have that there exists a constant $C > 0$ such that for $j \ge 1$
\[
 c_{j} \,  \frac{|\| G_{j - 1}^{(a + 1 \, , \, b + 1)} |\|}{|\| G_{j}^{(a , b)} |\|} \, 
(j + a + b + 1)  \ \le \ C \, j \, c_{j} \, .
\]
Thus,
\begin{align*}
\| D f \|_{s - 1 \, , \, (a + 1 \,  ,  \, b + 1)}^{2} &= \
\sum_{j = 0}^{\infty} \, (1 + j^{2})^{s - 1} \, C \, (j + 1)^{2} \, b_{j + 1}^{2} \
 \le \ C \, \sum_{j = 0}^{\infty} (1 + j^{2})^{s} \, c_{j}^{2} \ = \ C \, \| f \|_{s , (a , b)}^{2} \, .
\end{align*}
To establish that the mapping is onto, note that from \eqref{donto1} for 
$g(x) \ = \  \sum_{j = 0}^{\infty} d_{j} \,  \wtilde{G}_{j}^{(a + 1 \, , \, b + 1)}(x) 
\in H^{s - 1}_{(a + 1 \, , \, b + 1)}(\mrI)$, the function 
$f(x) \ =  \ \sum_{j = 0}^{\infty} d_{j} \, \frac{1}{j + a + b + 2} \, 
 \frac{|\| G_{j + 1}^{(a , b)} |\|}{|\| G_{j}^{(a + 1 \, , \, b + 1)} |\|} \wtilde{G}_{j + 1}^{(a  , \, b)}(x) 
\in H^{s}_{(a , b)}(\mrI)$, and $D f(x) \, = \, g(x)$.  \\
\mbox{  } \hfill \qed 

%%%%%%%%%%%
%%%%%%%%%%
\subsection{Equivalence of the spaces $H^{s}_{\rho^{(a , b)}}(\mrI)$ and $H^{s}_{(a , b)}(\mrI)$}
\label{ssec_equiv}
In this section we show that the spaces $H^{s}_{\rho^{(a , b)}}(\mrI)$ and $H^{s}_{(a , b)}(\mrI)$
are equivalent. In the first lemma we show equivalence of the space for $s \in \mathbb{N}_{0}$. Then, in the
second lemma we establish equivalence of the spaces for $s \in \mathbb{R} \backslash  \mathbb{N}_{0}$.

\begin{lemma}  \label{lmaeq1}
For $s \in \mathbb{N}_{0}$ the spaces $H^{s}_{(a , b)}(\mrI)$ and $H^{s}_{\rho^{(a , b)}}(\mrI)$
coincide, and their corresponding norms are equivalent. 
\end{lemma}
\textbf{Proof}:
For $s = m \in \mathbb{N}$, using Parseval's equality 
\begin{align*}
\| v \|_{m  , (a , b)}^{2} 
&= \  \sum_{k = 0}^{m} \| v^{(k)} \|_{L^{2}_{\rho^{(a , b)}}}^{2} 
\ = \  \sum_{k = 0}^{m} \, \sum_{j = 0}^{\infty} | v_{j}^{(k)} |^{2}   \\
&= \ \sum_{k = 0}^{m} \, \sum_{j = k}^{\infty} | v_{j - k}^{(k)} |^{2}   \\
&= \ \sum_{k = 0}^{m} \, \sum_{j = k}^{\infty} \left( \frac{|\| G_{j - k}^{(a + k \, , b + k)} |\|}{|\| G_{j}^{(a , b)} |\|}
  \frac{\Gamma(j + k + a + b + 1)}{\Gamma(j + a + b + 1)} \, v_{j} \right)^{2} 
  \ \ \ \mbox{(using \eqref{defvjmk})}  \\
&\sim   \ \sum_{k = 0}^{m} \, \sum_{j = k}^{\infty} \, j^{2 k} \, v_{j}^{2}  \ \ \ 
\mbox{(using \eqref{eqlma43}, and notation $0^{0} = 1$)}  \\
&\sim \sum_{j = k}^{\infty} \, (1 + j^{2} + \ldots +  j^{2m}) \, v_{j}^{2} \, .
\end{align*}

Noting that
\[
(1 + j^{2} + \ldots +  j^{2m}) \le (1 \, + \, j^{2})^{m} \le \left( \begin{array}{c} m  \\ \lfloor m / 2 \rfloor  \end{array} \right)
(1 + j^{2} + \ldots +  j^{2m}) \, ,
\]
where $\lfloor \cdot \rfloor$ denotes the integer part, we obtain
\[
\| v \|_{m , \rho^{(a,b)}}^{2} 
\ \sim \sum_{j = k}^{\infty} \, (1 + j^{2})^{m} \, v_{j}^{2} \ = \ \| v \|_{m , (a , b)}^{2} \, .
\]
\mbox{ } \hfill \qed

Recall that for $m, \, n \in \mathbb{N}$, $m < \, s  \, < n $, $\theta$ satisfying $s \, = \,(1 - \theta) m  \, + \, \theta \, n$,
with 
\be
  K(t , u) \, := \, \inf_{v \in  H^{n}_{\rho^{(a , b)}}} \left( \| u \, - \, v \|_{ m \, , \, \rho^{(a , b)}}
  \, + \, \| v \|_{n , \rho^{(a , b)}} \right) \, ,
  \label{def4K}
\ee
\[
   \| u \|_{s , \rho^{(a , b)}}   \, := \, 
   \| u \|_{[ H^{m}_{\rho^{(a , b)}} \, , \, H^{n}_{\rho^{(a , b)}} ]_{\theta , 2}}
   \, := \, \left( \int_{0}^{\infty} t^{-2 \theta} \, K(t , u)^{2} \, \frac{dt}{t} \right)^{1/2} \, .
\]

\begin{theorem}  \label{thmeq2}
The spaces $H^{s}_{(a , b)}(\mrI)$ and $H^{s}_{\rho^{(a , b)}}(\mrI)$
coincide, and their corresponding norms are equivalent. 
\end{theorem}
\textbf{Proof}: 
(This part of the proof is similar to an argument used in \cite{sch012}. See also \cite[Theorem 2.1]{bab011}.)
Let $m \, = \, n - 1$, $n$, $s$ and $\theta$ be as described above, and assume that 
$u \in H^{n - 1}_{\rho^{(a , b)}}(\mrI)$. 
Using 
$ \frac{1}{\sqrt{2}} (c + d) \ \le \ (c^{2} + d^{2})^{1/2} \ \le \ c + d $, and  
$ H^{j}_{\rho^{(a , b)}}(\mrI) \, = \,  H^{j}_{(a , b)}(\mrI)$ for $j \in \mathbb{N}_{0}$,
 we have that
 \[
 K(t , u) \ \sim \ \tilde{K}(t , u) \ = \ 
 \inf_{v \in H^{n}_{(a , b)}} \big( \| u - v \|_{n - 1 \, , \, (a , b)}^{2} \ + \ t^{2} \, 
 \| v \|_{n  ,  (a , b)}^{2} \big)^{1/2}  \, .
\]

For  $u \ = \ \sum_{k = 0}^{\infty} u_{k} \, \tilde{G}_{k}^{(a , b)}(x)$ 
and  $v \ = \ \sum_{k = 0}^{\infty} v_{k} \, \tilde{G}_{k}^{(a , b)}(x)$,
\[
\| u - v \|_{ n - 1 \, , \, (a , b)}^{2} \ = \ \sum_{k = 0}^{\infty} (1 + k^{2})^{n - 1} \, (u_{k} - v_{k})^{2} \, , \ \ \mbox{ and } \ \
\| v \|_{n , (a , b)}^{2} \ = \ \sum_{k = 0}^{\infty} (1 + k^{2})^{n} \, v_{k}^{2} \, .
\]
Then,
\begin{align}
\tilde{K}(t , u)^{2} &= \ 
 \inf_{ \{ (1 + k^{2})^{n/2} v_{k} \}_{k = 0}^{\infty} \in \, l^{2} }
 \sum_{k \, = \, 0}^{\infty} (1 + k^{2})^{n - 1} \, (u_{k} - v_{k})^{2} \ + \ t^{2} \, (1 + k^{2})^{n} \, v_{k}^{2}  \nonumber \\
 &\ge 
\ \sum_{k = 0}^{\infty} \inf_{v_{k} \in \mathbb{R}} 
\bigg( (1 + k^{2})^{n - 1} \, (u_{k} - v_{k})^{2} \ + \ t^{2} \, (1 + k^{2})^{n} \, v_{k}^{2}  \bigg) \, \label{msum1} .
\end{align}
Each term in the summation in \eqref{msum1} is minimized for
\begin{equation}
      v_{k} \ = \ \frac{u_{k}}{1 \, + \, t^{2} (1 + k^{2})} \, . 
\label{optvk}
\end{equation}
For these $v_{k}$ we have that
\begin{align*}
\sum_{k = 1}^{\infty} (1 \, + \, k^{2})^{n} \, v_{k}^{2} 
&= \ \sum_{k = 1}^{\infty} \frac{u_{k}^{2} \, (1 \, + \, k^{2})^{n}}{(1 \,+ \, t^{2} (1 \, + \, k^{2}) )^{2}}  \\
&= \ t^{-4} \, \sum_{k = 1}^{\infty}  \frac{u_{k}^{2} \,  \, t^{4} \,  (1 \, + \, k^{2})^{n}}{(1 \,+ \, t^{2} (1 \, + \, k^{2}) )^{2}}  
\ \le \  t^{-4} \, \sum_{k = 1}^{\infty}  u_{k}^{2} \, (1 \, + \, k^{2})^{n - 2} \\
&\le \ t^{-4} \, \| u \|_{n - 1 \, , \, (a , b)}^{2} \, ,
\end{align*}
which implies that $(1 + k^{2})^{n/2} v_{k} \in l^{2}$. Hence, 
\[
    \tilde{K}(t , u)^{2} \ = \ \sum_{k = 0}^{\infty} \frac{t^{2} (1 + k^{2})^{n}}{(1 \ + \ t^{2} (1 + k^{2}))} u_{k}^{2} \, .
\]
Then,
\begin{align*}
\| u \|_{s  , \rho^{(a , b)}} \ \sim \int_{0}^{\infty} t^{-2 \theta} \,   \tilde{K}(t , u)^{2} \, \frac{dt}{t} 
&= \ \int_{0}^{\infty}   \sum_{k = 0}^{\infty}  t^{-2 \theta} \,  \frac{t^{2} (1 + k^{2})^{n}}{(1 \ + \ t^{2} (1 + k^{2}))} u_{k}^{2} \, \frac{dt}{t}  \\
&= \  \sum_{k = 0}^{\infty} \int_{0}^{\infty}    t^{-2 \theta} \,  \frac{t^{2} (1 + k^{2})^{n}}{(1 \ + \ t^{2} (1 + k^{2}))} u_{k}^{2} \, \frac{dt}{t}  \, .
\end{align*}
Letting $\tau \ = \ (1 + k^{2})^{1/2} \, t \ \ \Leftrightarrow \ \ t \ = \ (1 + k^{2})^{- 1/2} \, \tau \ \ \Rightarrow
dt \ = \ (1 + k^{2})^{- 1/2} \, d\tau$,
\begin{align}
\| u \|_{s , \rho^{(a,b)}}^2 &\sim
 \sum_{k = 0}^{\infty} \int_{0}^{\infty}    \tau^{-2 \theta} \, (1 + k^{2})^{\theta} \,
  \frac{\tau^{2} \, (1 + k^{2})^{n - 1}}{(1 \ + \ \tau^{2})} u_{k}^{2} \, \frac{d\tau}{\tau}   \nonumber \\
&= \ \sum_{k = 0}^{\infty} (1 + k^{2})^{n - 1 + \theta} \, u_{k}^{2} \,  \int_{0}^{\infty}  
  \frac{\tau^{1 \, - \, 2 \theta}}{1 \ + \ \tau^{2}}  \, d\tau  \nonumber  \\
&= \ C_{\theta}^{2} \, \sum_{k = 0}^{\infty} (1 + k^{2})^{s} \, u_{k}^{2} 
\ = \ C_{\theta}^{2} \, \| u \|_{s , (a , b)}^{2} \, ,    \label{seq3}
\end{align}
where $C_{\theta} \ := \ \left(  \int_{0}^{\infty}  
  \frac{\tau^{1 \, - \, 2 \theta}}{1 \ + \ \tau^{2}}  \, d\tau  \right)^{1/2}$. \\
  
Thus, combining Lemma \ref{lmaeq1} and \eqref{seq3} we have for $s \ge 0$ 
the spaces $H^{s}_{(a , b)}(\mrI)$ and $H^{s}_{\rho^{(a , b)}}(\mrI)$
coincide, and their corresponding norms are equivalent. It also then follows that for $s < 0$ 
the corresponding (dual) spaces are equivalent. \\
 \mbox{ } \hfill \qed

We have the following straight forward lemma.
\begin{lemma} \label{nest1}
For $-1 < a_{1} < a$, $-1 < b_{1} < b$, $s \ge 0$, we have that $H_{(a_{1} , \, b_{1})}^{s}(\mrI) 
\subset H_{(a , \, b)}^{s}(\mrI)$. Additionally, for $s > 0$, $H_{(a , \, b)}^{-s}(\mrI)
\subset H_{(a_{1} , \, b_{1})}^{-s}(\mrI)$. 
\end{lemma}
\textbf{Proof}: \\
For $s \in \mathbb{N}_{0}$ the fact that for $x \in \mrI$, $(1 - x)^{a_{1}} \, x^{b_{1}} \ \le \
(1 - x)^{a} \, x^{b}$, together with definition \eqref{defHw} and the equivalence of
$H^{s}_{\rho^{(a , b)}}(\mrI)$ and $H^{s}_{(a , b)}(\mrI)$ establishes that 
$H_{(a_{1} , \, b_{1})}^{s}(\mrI) \subset H_{(a , \, b)}^{s}(\mrI)$.

Next, for $m \, = \, n-1$, $n$, $s$ and $\theta$ as defined in \eqref{def4K}, consider
\begin{align*}
K_{(a_{1} , \, b_{1})}(t , u) &:= \ 
\inf_{v \in H_{(a_{1} , \, b_{1})}^{n}}
\left( \| u  - v \|_{n - 1 \, , \, (a_{1} , \, b_{1})} \ + \ t \, \| v \|_{n , (a_{1} , \, b_{1})} \right)  \\
&\ge \ 
\inf_{v \in H_{(a_{1} , \, b_{1})}^{n}}
\left( \| u  - v \|_{n - 1 \, , \, (a , b)} \ + \ t \, \| v \|_{n , (a , \, b)} \right)  \\
&\ge \ 
\inf_{v \in H_{(a , \, b)}^{n}}
\left( \| u  - v \|_{n - 1 \, , \, (a , b)} \ + \ t \, \| v \|_{n , (a , \, b)} \right)  
\ \ \mbox{(as $H_{(a_{1} , \, b_{1})}^{n}(\mrI) \subset H_{(a , \, b)}^{n}(\mrI)$)} \\
&:= \ K_{(a , \, b)}(t , u) \, .
\end{align*}

Hence $\| u \|_{s , (a , b)} \le \| u \|_{s , (a_{1} , b_{1})}$, from which it follows that
$H_{(a_{1} , \, b_{1})}^{s}(\mrI) \subset H_{(a , \, b)}^{s}(\mrI)$.

That $H_{(a , \, b)}^{-s}(\mrI)
\subset H_{(a_{1} , \, b_{1})}^{-s}(\mrI)$ follows from the definition of the dual space and
that $H^{s}_{(a_{1} , b_{1})}(\mrI) \subset H^{s}_{(a , b)}(\mrI)$
for $s > 0$. \\
\mbox{  } \hfill \qed

We further illustrate the structure of the $H^{s}_{(a , b)}(\mrI)$ spaces with the following lemma,
whose proof is motivated by the proof of Theorem 2.5 in \cite{bab011}.
\begin{lemma}  \label{lmausp1}
Let $u(x) \, = \, x^{\mu}$. Then, $u \in H^{s}_{(a , b)}(\mrI)$ for $s  \, < \, 2 \mu + b + 1$. 
\end{lemma}
\textbf{Proof}: Let $\chi(x) \in C^{\infty}[0 , \infty)$ denote the cutoff function satisfying
\[
\chi(x) \ = \ \left\{ \begin{array}{rl}
1 & \mbox{ for } \ 0 < x \le 1/4  \\  0 & \mbox{ for } \   x \ge  3/4 \end{array}  \right. \, ,
\]
and let $\chi_{\delta}(x) \, := \, \chi(\frac{x}{\delta})$, for $\delta > 0$. Note that
\[
\chi_{\delta}(x) \ = \ \left\{ \begin{array}{rl}
1 & \mbox{ for } \ 0 < x \le \delta/4  \\  0 & \mbox{ for } \   x \ge  3 \delta /4 \end{array}  \right. \, ,
\quad \mbox{ and } \quad
\frac{d^{m}}{dx^{m}} \chi_{\delta}(x) \ = \ \left\{ \begin{array}{rl}
0 & \mbox{ for } \ 0 < x <  \delta/4  \\  0 & \mbox{ for } \   x  >  3 \delta /4 \end{array}  \right. \, 
\mbox{ for } m \in \mathbb{N} \, .
\]

For $\delta$ to be determined, let $u \, = \, v + w$ where $v(x) \ = \ \chi_{\delta}(x) \, u(x)$ and
$w(x) \ = \ \big(1 - \chi_{\delta}(x) \big) \, u(x)$. We have that
\[
\left| \frac{d^{m} v(x)}{dx^{m}} \right| \ \le \ C \, \sum_{j = 0}^{m} \delta^{-(m - j)} \, x^{\mu - j} \, , 
\ \mbox{ and  is zero for }  \ x \ > \ 3 \delta / 4 \, .
\]
Thus,
\begin{align*}
\int_{\mrI} (1 - x)^{a + m} \, x^{b + m} \, \left| \frac{d^{m} v(x)}{dx^{m}} \right|^{2} \, dx
&\le \  C \, \int_{0}^{3 \delta / 4} \, \sum_{j = 0}^{m} \delta^{-2 (m - j)} \, x^{2 \mu \, - \, 2 j \, + b + m} \, dx  \\
&\le \ C \, \sum_{j = 0}^{m} \delta^{-2 (m - j)} \, \delta^{2 \mu \, - \, 2 j \, + b + m + 1} \, ,
 \ \ \mbox{ provided } 2 \mu \, - \, 2 j \, + b + m > -1 \, ,   \\
 &\le \ C \, \delta^{2 \mu \,  + \, b + 1 - m} \, , 
\end{align*}
which implies that, for $m \, < \, 2 \mu \,  + \, b + 1$, $v \in H^{m}_{\rho^{(a ,  b)}}(\mrI)$ and
\be
\| v \|_{m \, , \, \rho^{(a ,  b)}}^{2} \ \le \ C \, \delta^{2 \mu \,  + \, b + 1 - m} \, .
\label{ewst21}
\ee

Next, consider $w(x)$.
\be
\left| \frac{d^{m} w(x)}{dx^{m}} \right| \ \le \ C \, \left( \big(1 - \chi_{\delta}(x) \big) x^{\mu - m} \ + \ 
 \sum_{j = 0}^{m - 1} x^{\mu - j} \,  \frac{d^{m - j}}{dx^{m - j}} \big(1 - \chi_{\delta}(x) \big) \right)  \,  . 
\label{wterm1}
\ee
The first term on the RHS of \eqref{wterm1} vanishes for $x < \delta / 4$, and the second term vanishes for
$x \, < \, \delta/4$ and $x \, > \, 3 \delta/4$. Using this,
\begin{align}
\int_{\mrI} (1 - x)^{a + m} \, x^{b + m} \, \left| \frac{d^{m} w(x)}{dx^{m}}  \right|^{2} \, dx
&\le \ C \, \left( \int_{\delta/4}^{1} x^{b + m} \, x^{2 \mu \, - \, 2 m} \, dx \ + \ 
\int_{\delta / 4}^{3 \delta / 4} \,  x^{b + m} \,  
 \sum_{j = 0}^{m - 1} x^{2 \mu \, - \, 2 j} \,  \delta^{-2 ( m - j)} \,  dx \right)   \nonumber \\
&\le \  C \, \left( \int_{\delta/4}^{1}  x^{2 \mu \, + b \, - \, m} \, dx \ + \  
 \sum_{j = 0}^{m - 1} \, \delta^{-2 m \, + \, 2 j} \, \int_{\delta / 4}^{3 \delta / 4} \, x^{2 \mu \, + \,  b \, + \, m \, - \, 2 j} \,  dx \right)
 \nonumber  \\
&
\le \ \left\{ \begin{array}{rl}
           C \, \left( 1 \, + \, \delta^{2 \mu \, + \, b \, - \, m \, + \, 1}  \right) & \mbox{ if }  2 \mu \, + \, b \, - \, m \, \neq \, -1 \\
           C \, \left( 1 \, + \, | \log \delta |  \right) & \mbox{ if }  2 \mu \, + \, b \, - \, m \, = \, -1 
           \end{array} \right. \, .  \label{ewst23}
\end{align}
Hence, for $n \, > \, 2 \mu \, + \, b \, + \, 1$, $w \in H^{n}_{\rho^{(a ,  b)}}(\mrI)$ and
\be
   \| w \|_{n \, , \, \rho{(a ,  b)}}^{2} \ \le \ C \, \delta^{2 \mu \, + \, b \, + \, 1 \, - \, n} \, .
   \label{ewst24}
\ee
(\textbf{Remark}: For $n \, > \, 2 \mu \, + \, b \, + \, 1$ the exponent of $\delta$ in \eqref{ewst24} is negative,
so the $'1'$ term in \eqref{ewst23} is bounded by the $\delta$ term.)

We have from \eqref{ewst21} and \eqref{ewst24}
\begin{align}
K(t , u) &= \ \inf_{u \, = \, u_{1} + u_{2}} \left( \| u_{1} \|_{H^{m}_{(a , b)}} \ + \ t \, \| u_{2} \|_{H^{n}_{(a , b)}} 
 \right)  \label{ewst251}  \\
 &\le \ \| v \|_{m \, , \, \rho^{(a , b)}} \ + \ t \, \| w \|_{n \, , \, \rho{(a , b)}}   \nonumber \\
 &\le \ C \, \left( \delta^{(2 \mu \,  + \, b + 1 - m) / 2} \ + \ t \, \delta^{(2 \mu \, + \, b \, + \, 1 \, - \, n) / 2} \right) \, .
 \label{ewst252}
\end{align} 

Setting $\delta \ = \ t^{2 / (n - m)}$ leads to  $K(t , u) \, \le \, C \, t^{(2 \mu \,  + \, b + 1 - m) / (n - m)}$.

Recall that
\be
\| u \|^{2}_{[H^{m}_{(a ,  b)} \, , \, H^{n}_{(a ,  b)}]_{\theta , 2}} \ = \ 
\int_{0}^{\infty} \, t^{-2 \theta} \left( K(t , u) \right)^{2} \, \frac{dt}{t} \, .
\label{ewst26}
\ee
The larger the value of $\theta$ $(0 < \theta < 1)$ in \eqref{ewst26} such that the integral is finite, 
the ``nicer'' (i.e., more regular) is the function $u$. Hence from 
\eqref{ewst26}, we are interested in the integrand about $t = 0$. We have trivially that for $u_{1} = u$, $u_{2} = 0$
in \eqref{ewst251} that $K(t , u) \, \le \, \| u \|_{H^{m}_{(a ,  b)}}  \, \le \, C$. Hence it follows that
\[
  K(t , u) \, \le \, \left\{ \begin{array}{ll}
  C \, t^{(2 \mu \,  + \, b + 1 - m) / (n - m)} & \mbox{ for } 0 < t < 1 \\
  C &  \mbox{ for } t \ge 1 
  \end{array} \right.   \, .
\]  

Using \eqref{ewst26},
\begin{align}
\| u \|^{2}_{[H^{m}_{(a ,  b)} \, , \, H^{n}_{(a ,  b)}]_{\theta , 2}} &\le \
\int_{0}^{1} C \, t^{-2 \theta \, - \, 1 \, + \, 2 (2 \mu \,  + \, b + 1 - m) / (n - m)} \, dt \ + \ 
\int_{1}^{\infty}  C \, t^{-2 \theta \, - \, 1} \, dt \ < \ \infty \, , \nonumber  \\
\mbox{ if } \ \ \theta &< \ (2 \mu \,  + \, b + 1 - m) / (n - m) \, .  \label{ewst27}
\end{align}
For $s \ = \ (1 - \theta) m \ + \ \theta n \ = \ m \ + \ \theta (n - m)$, then  
$s \ < \ 2 \mu \,  + \, b + 1$ using  \eqref{ewst27} \, .

Hence we can conclude that $u(x) \ = \ x^{\mu} \, \in \, H^{s}_{(a ,  b)}(\mrI)$ for $s < 2 \mu \,  + \, b + 1$. \\
\mbox{  } \hfill \qed

%%%%%%
%%%%%%
We now present a result which connects the weighted Sobolev spaces with the spaces of continuous functions.
\begin{theorem}(See \cite[Theorem 4.14]{aco181})  \label{cntythm}
Let $a, b > -1$, $k \in \mathbb{N}_{0}$ and 
$s >  k \, + \, 1 \, + \, \max\{a + k \,  , \, b + k \, , \, -1/2 \}$. Then we have a continuous embedding
$H^{s}_{(a , b)}(\mrI) \subset C^{k}(\mrI)$ of $H^{s}_{(a , b)}(\mrI)$ into the Banach space
$C^{k}(\mrI)$ with the norm $\| v \|_{C^{k}} \, = \, \| v \|_{\infty} \, + \, \| v^{(k)} \|_{\infty}$. 
\end{theorem}
\textbf{Proof}: 
We begin by establishing that $v(x) \in C(\mrI)$. From the representation  \eqref{vexp1},
consider the sequence of partial sums 
$\{ v_{n}(x) \}_{n = 0}^{\infty}$, where, using \eqref{spm22g}, 
\[
v_{n}(x) \ = \ \sum_{j = 0}^{n} v_{j} \, \tilde{G}_{j}^{(a , b)}(x)  \ 
= \  \sum_{j = 0}^{n} v_{j} \, \left( (2j \, + \, a \, + \, b \, + 1) 
   \frac{\Gamma(j + 1) \, \Gamma(j + a + b + 1)}{\Gamma(j + a + 1) \, \Gamma(j + b + 1)}
   \right)^{1/2} \, 
        G_{j}^{(a , b)}(x) \, .
\]
From \cite[Theorem 7.32.1]{sze751}, for $x \in \mrI$, $|G_{n}^{(a , b)}(x)| \sim n^{q}$, where $q = \max\{ a, b, -1/2\}$.
Using this, together with \eqref{nrmGbeh}, for $n > m$
\begin{align}
| v_{n}(x) \, - \, v_{m}(x) | &= \
\left| \sum_{j = m + 1}^{n}    v_{j} \, \left( (2j \, + \, a \, + \, b \, + 1) 
   \frac{\Gamma(j + 1) \, \Gamma(j + a + b + 1)}{\Gamma(j + a + 1) \, \Gamma(j + b + 1)}
   \right)^{1/2} \, 
        G_{j}^{(a , b)}(x) \right|   \nonumber \\
&\le \ C \,     \sum_{j = m + 1}^{n}   | v_{j} | \, j^{1/2} \, j^{\max\{a \, ,  \, b \, , \, -\frac{1}{2} \} } 
         \nonumber \\
&= \ C   \,    \sum_{j = m + 1}^{n}  ( 1 \, +  \, j^{2})^{- s/2} \,   
j^{\frac{1}{2} \, + \, \max\{a \, ,  \, b \, , \, -\frac{1}{2} \} } \, ( 1 \, +  \, j^{2})^{s/2} \,   | v_{j} |
          \nonumber \\
&\le \ C   \,    \left( \sum_{j = m + 1}^{n}  ( 1 \, +  \, j^{2})^{- s} \,   
j^{1 \, + \, 2 \max\{a \, ,  \, b \, , \, -\frac{1}{2} \} }  \right)^{1/2} \, 
 \left( \sum_{j = m + 1}^{n}  ( 1 \, +  \, j^{2})^{s} \,   v_{j}^{2}  \right)^{1/2}  \nonumber \\
&\le \ C   \,    \left( \sum_{j = m + 1}^{n}  
j^{1 \, - \, 2 s \, + \, 2 \max\{a \, ,  \, b \, , \, -\frac{1}{2} \} }  \right)^{1/2} \, 
 \| v \|_{s , (a , b)} \, .         \label{vinf1}         
\end{align}        
Thus, for $1 \, - \, 2 s \, + \, 2 \max\{a \, ,  \, b \, , \, -\frac{1}{2} \} \, < \, -1$, i.e., for
$s \, > \, 1 \, + \, \max\{a \, ,  \, b \, , \, -\frac{1}{2} \}$, $\{ v_{n}(x) \}_{n = 0}^{\infty}$ converges
uniformly on $\mrI$. As the limit function of a sequence of uniformly convergent, continuous functions is also
continuous, it follows that $v(x) \in C(\mrI)$.

To obtain the continuity of $v^{(l)}(x)$, for $1 \le l \le k$, we proceed in a similar manner. Using \eqref{eqC4}, 
$m \ge l$
\begin{align*}
   &     | v_{n}^{(l)}(x) \, - \, v_{m}^{(l)}(x) |  \ = \\ 
   & \quad  \    \left| \sum_{j = m+1}^{n} v_{j} \, 
        \left( (2j \, + \, a \, + \, b \, + 1) 
   \frac{\Gamma(j + 1) \, \Gamma(j + a + b + 1)}{\Gamma(j + a + 1) \, \Gamma(j + b + 1)}
   \right)^{1/2} \, 
   \frac{\Gamma(j + l + a + b + 1)}{  \Gamma(j + a + b + 1)} 
        G_{j - l}^{(a + l \,  , \, b + l)}(x) \right| \, .
\end{align*}    

Using
\[ \frac{\Gamma(j + l + a + b + 1)}{  \Gamma(j + a + b + 1)} \ \sim \ j^{l} \, \quad 
\mbox{ for $j$ large }  \, ,
\]

\begin{align}
| v_{n}^{(l)}(x) \, - \, v_{m}^{(l)}(x) |&\le \  C  \sum_{j = l}^{n} | v_{j}  | \, 
        j^{1/2} \,  j^{l} \, 
        (j - l)^{\max\{a + l \,  , \, b + l \, , \, -1/2 \}}  \nonumber \\
&= \ C   \,    \left| \sum_{j = m + 1}^{n}  ( 1 \, +  \, j^{2})^{- s/2} \,   
j^{l \, + \, \frac{1}{2} \, + \, \max\{a + l \,  , \, b + l \, , \, -1/2 \} } \, ( 1 \, +  \, j^{2})^{s/2} \,   | v_{j} |
        \right|  \nonumber \\
&\le \ C   \,    \left( \sum_{j = m + 1}^{n}  
j^{2 l \, + \, 1 \, - \, 2 s \, + \, 2 \max\{a + l \,  , \, b + l \, , \, -1/2 \} }  \right)^{1/2} \, 
 \| v \|_{s , (a , b)} \, .       \label{vinf2}       
\end{align}          

Thus for $2 l \, + \, 1 \, - \, 2 s \, + \, 2 \max\{a + l \,  , \, b + l \, , \, -1/2 \} \, < \, -1$, i.e., for
$s \, > \,  l \, + \, 1 \, + \, \max\{a + l \,  , \, b + l \, , \, -1/2 \}$, $\{ v^{(l)}_{n}(x) \}_{n = 0}^{\infty}$ converges
uniformly on $\mrI$. Furthermore, as $\{ v_{n}(x) \}_{n = 0}^{\infty}$ converges
to $v(x)$, then $\{ v^{(l)}_{n}(x) \}_{n = 0}^{\infty}$ converges to $v^{(l)}(x) \in C(\mrI)$ (see \cite[Theorem 7.17]{rud761}).

A straight forward modification of the argument used in \eqref{vinf1} and \eqref{vinf2} can be used to
establish that, for $0 \le l \le k$, $\| v^{(l)} \|_{\infty} \, \le \, C \,  \| v \|_{H^{s}_{(a , b)}}$, where the 
precise value for $C$ depends upon $l$.  \\
\mbox{  } \hfill \qed

%%%%%
%%%%% 
 \setcounter{equation}{0}
\setcounter{figure}{0}
\setcounter{table}{0}
\setcounter{theorem}{0}
\setcounter{lemma}{0}
\setcounter{corollary}{0}
\setcounter{definition}{0}
\section{Regularity of the solution to the fractional diffusion equation}
\label{secRFD}
Now we can give sharp regularity results for the fractional diffusion equation  \eqref{DefProb1},\eqref{DefBC1}.
\begin{theorem} \label{regmapFDiff}
For $\mcL^{\alpha}_{r} (\cdot)$ defined for $1 < \alpha < 2$, $0 \le r \le 1$, let $\beta$ be determined by 
\textbf{Condition A}. Then the mapping $\mcL^{\alpha}_{r} (\cdot)  \, : \, \rho^{(\alpha - \beta , \beta)}(x)
 \otimes H^{s + \alpha}_{(\alpha - \beta \, , \beta)}(\mrI) \rightarrow H^{s}_{(\beta \, , \alpha - \beta)}(\mrI)$
is bijective, continuous,  and has a continuous inverse.
\end{theorem}

\textbf{Proof}: From \cite{erv162, jia181},
\[
   \mcL^{\alpha}_{r} \left( \rho^{(\alpha - \beta \, , \, \beta)} \, \tilde{G}_{k}^{(\alpha - \beta \, , \, \beta)} \right)(x)
   \ = \ \lambda_{k} \, \tilde{G}_{k}^{(\beta \, , \, \alpha - \beta)}(x) \, , \ \ 
   \mbox{ where } \ \ 
   \lambda_{k} \ = \ - c_{*}^{*} \frac{\Gamma(k + 1 + \alpha)}{\Gamma(k + 1)} \, , \ k = 0, 1, 2, \ldots,
\]
and $c_{*}^{*} $ given by \eqref{defcss}. 

Using Stirling's formula,
\begin{align}
\lambda_{k}^{2} \ = \ (c_{*}^{*})^{2} \frac{\Gamma(k + 1 + \alpha)^{2}}{\Gamma(k + 1)^{2}} 
&\sim  \ (k + 1)^{2 \alpha} \, \ \mbox{ as } k \rightarrow \infty \, , \nonumber   \\
\mbox{i.e.}, \ \ c \, (1 + k^{2})^{\alpha} &\le \ \lambda_{k}^{2} \ \le  \ C \, (1 + k^{2})^{\alpha} \, .  \label{hhh1}
\end{align}

Let $\phi(x) \ = \ \sum_{k = 0}^{\infty} \phi_{k} \tilde{G}_{k}^{(\alpha - \beta \, , \, \beta)}(x) 
\in H_{(\alpha - \beta \, , \, \beta)}^{s + \alpha}(\mrI)$, and 
 $\phi_{N}(x) \ = \ \sum_{k = 0}^{N} \phi_{k} \tilde{G}_{k}^{(\alpha - \beta \, , \, \beta)}(x) $. Then 
 $\{ \phi_{n} \}_{n = 1}^{\infty}$ is a Cauchy sequence in $H_{(\alpha - \beta \, , \, \beta)}^{s + \alpha}(\mrI)$,
 with $\lim_{n \rightarrow \infty} \phi_{n} \ = \ \phi$.
 
 Consider, 
 \[
 f_{N}(x) \ = \ \mcL^{r}_{\alpha} \left( \rho^{(\alpha - \beta \, , \, \beta)} \phi_{N} \right)(x) 
 \ = \ \sum_{k = 0}^{N} \lambda_{k} \, \phi_{k} \, \tilde{G}_{k}^{(\beta \, , \, \alpha - \beta)}(x)  \
 \in  H_{(\beta \, , \, \alpha - \beta)}^{s}(\mrI) \, .
\]
Then using \eqref{hhh1},
\begin{align*}
\|f_{N} - f_{M} \|_{s , (\beta \, , \, \alpha - \beta)}^{2} 
&= \ \sum_{k = M}^{N}  (1 + k^{2})^{s} \, \left( \lambda_{k} \, \phi_{k} \right)^{2} 
\ \le \ C \, \sum_{k = M}^{N}  (1 + k^{2})^{s} \,  (1 + k^{2})^{\alpha}  \, \phi_{k}^{2}   \\
&= \ C \, \|\phi_{N} - \phi_{M} \|_{s + \alpha \, ,  \, (\alpha - \beta \, , \, \beta)}^{2} \, .
\end{align*}
Thus $\{ f_{n} \}_{n = 1}^{\infty}$ is a Cauchy sequence in $H_{(\beta \, , \, \alpha - \beta)}^{s}(\mrI)$. It then follows
that $f \, := \, \lim_{n \rightarrow \infty} f_{n}$ satisfies 
\begin{align*}
f &= \ \mcL^{\alpha}_{r} \left( \rho^{(\alpha - \beta \, , \, \beta)} \, \phi \right)(x) \ = \ 
\mcL^{\alpha}_{r} \left( \rho^{(\alpha - \beta \, , \, \beta)} \, \sum_{k = 0}^{\infty} \phi_{k} \tilde{G}_{k}^{(\alpha - \beta \, , \, \beta)} \right)(x)  \\
&= \ \sum_{k = 0}^{\infty} \phi_{k} \,  \mcL^{\alpha}_{r} \left( \rho^{(\alpha - \beta \, , \, \beta)}
  \tilde{G}_{k}^{(\alpha - \beta \, , \, \beta)} \right)(x)  \\
&= \ \sum_{k = 0}^{\infty} \lambda_{k} \, \phi_{k} \, \tilde{G}_{k}^{(\beta \, , \, \alpha - \beta)}(x) \, .  
\end{align*}
%%Then,
%%\begin{equation}
%%\| f \|_{s , (\beta \, , \, \alpha - \beta)}^{2} \ = \ \sum_{k = 0}^{\infty} (1 + k^{2})^{s} \, \lambda_{k}^{2} \, \phi_{k}^{2} \, .
%%\label{hest0}
%%\end{equation}
%%As
%%\begin{align*}
%%\lambda_{k}^{2} &= \ (c_{*}^{*})^{2} \frac{\Gamma(k + 1 + \alpha)^{2}}{\Gamma(k + 1)^{2}} 
%%\ \sim  \ (k + 1)^{2 \alpha} \, \ \mbox{ as } k \rightarrow \infty
%%\mbox{ (using Stirling's formula)}  \\
%%%
%%&\le \ C \, (1 + k^{2})^{\alpha} \, ,
%%\end{align*}
%
Also, using \eqref{hhh1},
\begin{align}
\| f \|_{s , (\beta \, , \, \alpha - \beta)}^{2} &\le \
 \sum_{k = 0}^{\infty} (1 + k^{2})^{s}  \, \lambda_{k}^{2} \, \phi_{k}^{2}  %\nonumber \\
\ \le \ C \,  \sum_{k = 0}^{\infty} (1 + k^{2})^{s + \alpha} \,  \phi_{k}^{2} 
\ = \ C \, \| \phi \|_{s + \alpha \, , \, (\alpha - \beta \, , \beta)}^{2} \, .   \label{hest1}
\end{align}
Hence $\mcL^{\alpha}_{r}(\cdot)$ is a one-to-one, continuous mapping from 
$\rho^{(\alpha - \beta , \beta)}(x)
 \otimes H^{s + \alpha}_{(\alpha - \beta \, , \beta)}(\mrI)$ into $H^{s}_{(\beta \, , \alpha - \beta)}(\mrI)$.

Next, for $f(x) \ = \ \sum_{k = 0}^{\infty} f_{k} \tilde{G}_{k}^{(\beta \, , \, \alpha - \beta)}(x) 
\in H_{(\beta \, , \, \alpha - \beta)}^{s}(\mrI)$, let 
$\phi(x) \ = \ \sum_{k = 0}^{\infty} \frac{1}{\lambda_{k}} f_{k} \tilde{G}_{k}^{(\alpha - \beta \, , \, \beta)}(x)$.
Note that $\mcL^{\alpha}_{r} \left( \rho^{(\alpha - \beta , \beta)} \, \phi \right)(x) \ = \ f(x)$, 
and using a similar argument to that in \eqref{hest1},
$\| \phi \|_{s + \alpha \, , \, (\alpha - \beta \, , \, \beta)} \, \le \, C \, \| f \|_{s , (\beta \, , \, \alpha - \beta)}$, from which 
the stated result then follows. \\
\mbox{  } \hfill \qed

Using this theorem we obtain the following sharp regularity result for \eqref{DefProb1},\eqref{DefBC1}.

\begin{corollary}  \label{exuncor}
For $f \in H^{s}_{(\beta \, , \, \alpha - \beta)}(\mrI)$ there exists a unique solution $u$ to \eqref{DefProb1},\eqref{DefBC1},
which can be
expressed as $u(x) \ = \ \rho^{(\alpha - \beta , \beta)}(x) \phi(x)$, where 
$\phi(x) \in H^{s + \alpha}_{(\alpha - \beta \, , \, \beta)}(\mrI)$, with 
$\| \phi \|_{s + \alpha \,  , \, (\alpha - \beta \, , \, \beta)} \, \le \, C \, \| f \|_{s , (\beta \, , \, \alpha - \beta)}$ \, .
\end{corollary}
\mbox{ } \hfill \qed

Combining Theorem \ref{cntythm} and Corollary \ref{exuncor} we obtain a sufficient condition for the solution
of \eqref{DefProb1},\eqref{DefBC1} to be continuous on $\mrI$.

\begin{corollary} \label{corctssol}
Let $1 < \alpha < 2$, $0 \le r \le 1$, and $\beta$ be determined
by \textbf{Condition A}. Then for $f \in H^{s}_{(\beta \, ,  \, \alpha - \beta)}(\mrI)$ where
\begin{equation}
  s \ > \ (1 - \alpha) \, + \, \max\{ \alpha - \beta \, , \, \beta \} \, ,     \label{defscts}
\end{equation}
the solution of \eqref{DefProb1},\eqref{DefBC1} is 
continuous on $\mrI$.
\end{corollary}
\mbox{  } \hfill \qed

\textbf{Note}: For $r = 0$ ($r = 1$) Corollary \ref{corctssol} implies for $f \in H^{s}_{(1 \, ,  \, \alpha - 1)}(\mrI)$
\big($f \in H^{s}_{(\alpha - 1 \, ,  \, 1)}(\mrI) \big)$, where $s \, > \, 2 - \alpha$, that $u \in C(\mrI)$. \\
For $r = 1/2$ Corollary \ref{corctssol} implies for $f \in H^{s}_{(\alpha/2 \, ,  \, \alpha/2)}(\mrI)$,
where $s \, > \, 1 - \alpha/2$, that $u \in C(\mrI)$. \\

%%%%%
%%%%% 
 \setcounter{equation}{0}
\setcounter{figure}{0}
\setcounter{table}{0}
\setcounter{theorem}{0}
\setcounter{lemma}{0}
\setcounter{corollary}{0}
\setcounter{definition}{0}
\section{Weighted Sobolev Space defined using a Sobolev-Slobodeckij type semi-norm}
\label{sec_Wspace3}
To obtain the regularity results for \eqref{DefProb2},\eqref{DefBC2} in the following section we use
the regularity result from  \eqref{DefProb1},\eqref{DefBC1} and a boot-strapping argument. To
apply the boot-strapping argument we need a precise estimate for $t$ such that if
$\phi(x) \in H^{r}_{(\alpha - \beta \, , \, \beta)}(\mrI)$, then
$(1 - x)^{\alpha - \beta} \, x^{\beta} \, \phi(x) \in H^{t}_{(\beta \, , \, \alpha - \beta)}(\mrI)$.

In Sections \ref{sec_Wspace1} and \ref{sec_Wspace2} we gave two characterizations of the
weighted Sobolev space $H^{s}_{(a ,  b)}(\mrI)$. Neither of these definitions are easily
to apply in order to determine the value for $t$ in order that
$(1 - x)^{\alpha - \beta} \, x^{\beta} \, \phi(x) \in H^{t}_{(\beta \, , \, \alpha - \beta)}(\mrI)$.
Therefore, following the work of Bernardi, Dauge, and Maday \cite{ber921}, 
and Fdez-Manin, and Munoz-Sola \cite{fde951}
we introduce
a further characterization of $H^{s}_{(a  , \, b)}(\mrI)$.

\textbf{Remark}: The work in \cite{ber921, fde951} is presented for a single weighted space, i.e., 
$H^{s}_{(a  , \, a)}(\mrI)$. However,  it extends in a straight forward manner to 
$H^{s}_{(a  , \, b)}(\mrI)$ spaces.

\subsection{Equivalent definition for $H^{s}_{(a , b)}(\mrI)$}
\label{ssec_eqdf3}

\begin{figure}[!ht]
%\begin{figure}[t]
\begin{minipage}{.46\linewidth}

%\begin{figure}[!ht]
\begin{center}
 \includegraphics[height=2.5in]{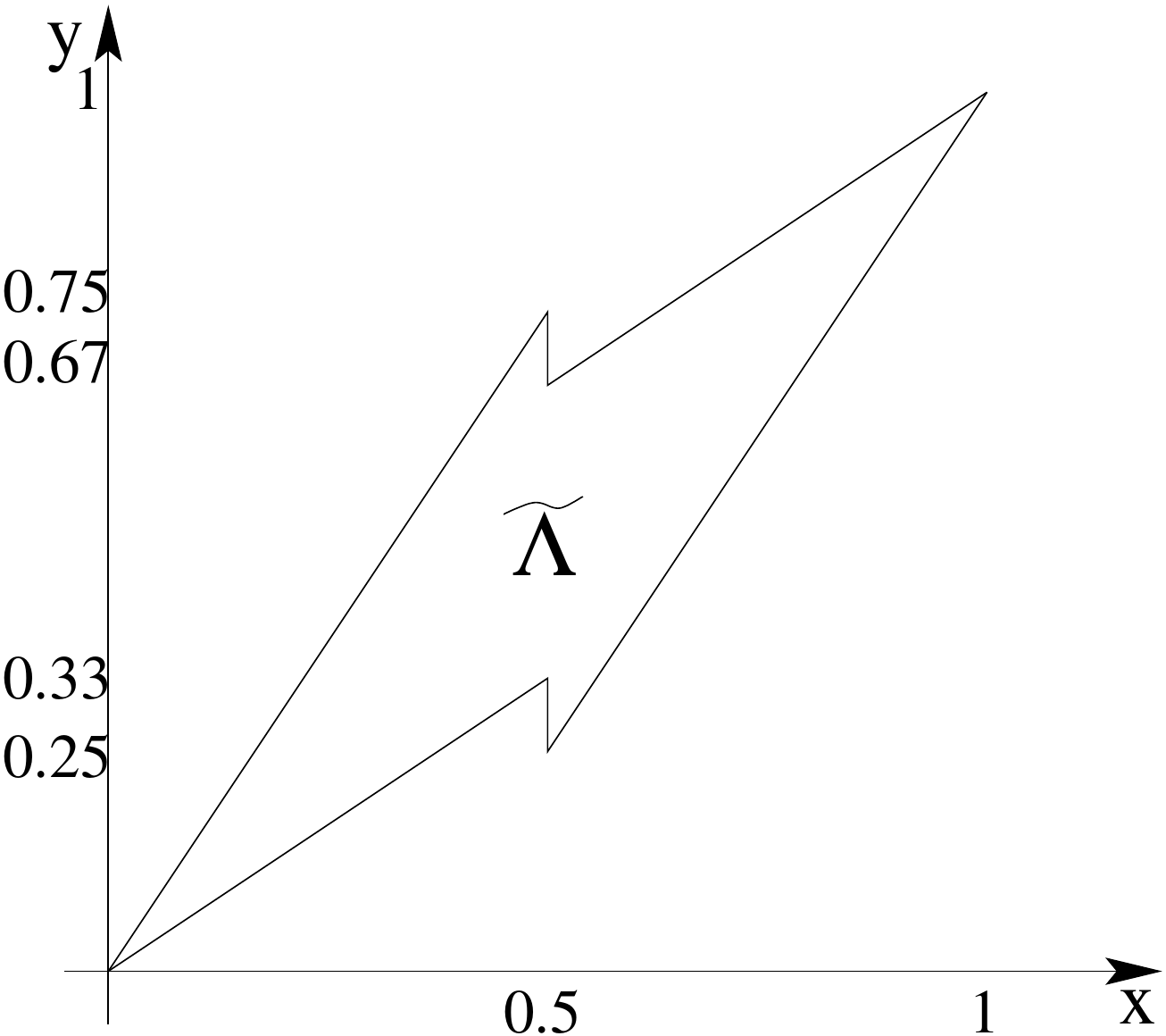}
   \caption{Domain $\widetilde{\Lambda}$.}
   \label{figdomlam}
\end{center}
%\end{figure}

\end{minipage} \hfill
\begin{minipage}{.46\linewidth}
 
\begin{center}
 \includegraphics[height=2.5in]{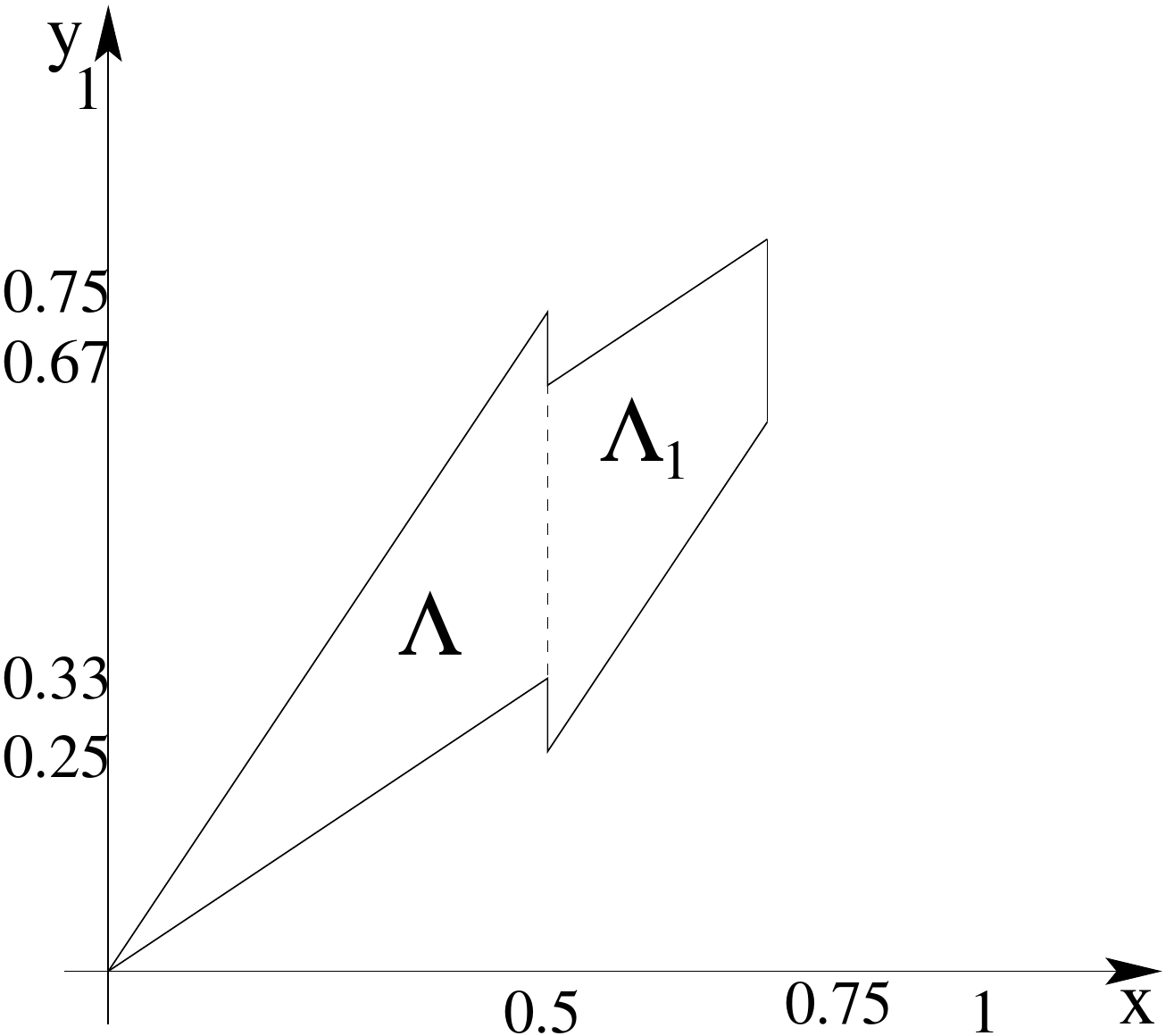}
   \caption{Domain $\Lambda^{*} \, = \, \Lambda \cup \Lambda_{1}$.}
   \label{figdomlams}
\end{center}
\end{minipage} 
\end{figure}

For $s > 0$, let $s \ = \ \floor{s} \, + \, r$, where $0 < r < 1$. Let the
semi-norm, $| \cdot |_{H^{s}_{(a , b)}(\mrI)}$, and norm, $\| \cdot \|_{H^{s}_{(a , b)}(\mrI)}$ 
 be defined as 
 \begin{align*}
  | f |_{H^{s}_{(a , b)}(\mrI)}^{2} &:= \ 
  \iint_{\widetilde{\Lambda}} (1 - x)^{a + s} \, x^{b + s} \, 
  \frac{ | D^{\floor{s}} f(x) \, - \, D^{\floor{s}} f(y) |^{2}}{ | x \, - \, y |^{1 \, + \, 2 (s - \floor{s})}} dy \,  dx \,   \\
 &  \\
\mbox{and } \ 
  \| f \|_{H^{s}_{(a , b)}(\mrI)}^{2} &:= \left\{ \begin{array}{rl}  
   \sum_{j = 0}^{s} \| D^{j} f \|_{L^{2}_{(a + j \, , \, b + j)}(\mrI)}^{2} \, , & \mbox{ for } s \in \mathbb{N}_{0}  \\
 \sum_{j = 0}^{\floor{s}} \| D^{j} f \|_{L^{2}_{(a + j \, , \, b + j)}(\mrI)}^{2} \ + \ 
  | f |_{H^{s}_{(a , b)}(\mrI)}^{2} \, , & \mbox{ for } s \in \mathbb{R}^{+} \backslash \mathbb{N}_{0} 
 \end{array} \right. \, ,  
\end{align*}
where (see Figure \ref{figdomlam})
\be
\widetilde{\Lambda} \ := \ 
\bigg\{ (x , y) \, : \, \frac{2}{3} x < y < \frac{3}{2} x , \, 0 < x < \frac{1}{2} \bigg\} \cup
\left\{ (x , y) \, : \, \frac{3}{2} x - \frac{1}{2} \,  < y <  \, \frac{2}{3} x + \frac{1}{3} , \, 1/2 \le x < 1 \right\} \, .
\label{defLtda}
\ee

%%\begin{figure}[!ht]
%%\begin{center}
%% %\includegraphics[height=2.5in]{SolEx1.eps}
%% \includegraphics[height=2.5in]{Dom2.pdf}
%%   \caption{Domains $\widetilde{\Lambda}$ and $\Lambda^{*}$.}
%%   \label{domfig1}
%%\end{center}
%%\end{figure}

\begin{lemma} \cite{fde951} Let $s \ge 0$ with $s \neq 1 + a$ if $a \in (-1 , 0)$ and $s \neq 1 + b$ if
$b \in (-1 , 0)$. Then, 
$\{ f \, : \, f \mbox{ is measurable and } \| f \|_{H^{s}_{(a , b)}(\mrI)} \, < \, \infty \} \, = \, H^{s}_{(a , b)}(\mrI)$.
\end{lemma}

\textbf{Proof}: From \cite{fde951} it follows that 
$\{ f \, : \, f \mbox{ is measurable and } \| f \|_{H^{s}_{(a , b)}(\mrI)} \, < \, \infty \} \, = \, H^{s}_{\rho^{(a , b)}}(\mrI)$.
The equality of the spaces  $H^{s}_{\rho^{(a , b)}}(\mrI)$ and $H^{s}_{(a , b)}(\mrI)$ was established in
Section \ref{sec_Wspace2}.  \\
\mbox{ } \hfill \qed

The $H^{s}_{(a , b)}(\mrI)$ space a function $f$ lies in is determined by its behavior at: (i) the left endpoint ($x = 0$), 
(ii) the right endpoint ($x = 1$), and (iii) away from the endpoints. In order to separate the consideration of the endpoint 
behaviors, following \cite{ber921} we introduce the following function space $H^{s}_{(\gamma)}(\mrJ)$. 
Let $\mrJ \, := \, (0 , \, 3/4)$, and
\begin{align*}
\Lambda^{*} &:= \ 
\left\{ (x , y) \, : \, \frac{2}{3} x < y < \frac{3}{2} x , \, 0 < x < \frac{1}{2} \right\} \cup
\left\{ (x , y) \, : \, \frac{3}{2} x - \frac{1}{2} \,  < y <  \, \frac{2}{3} x + \frac{1}{3} , \, 1/2 \le x < 3/4 \right\} \,  \\
&:= \ \Lambda \cup \Lambda_{1}  \ \ \mbox{(see Figure \ref{figdomlams})} \, .
\end{align*}

Introduce the semi-norm and norm
\begin{align}
 | f |_{H^{s}_{(\gamma)}(\mrJ)}^{2} &:= \ 
  \iint_{\Lambda}  x^{\gamma + s} \, 
  \frac{ | D^{\floor{s}} f(x) \, - \, D^{\floor{s}} f(y) |^{2}}{ | x \, - \, y |^{1 \, + \, 2 (s - \floor{s})}} dy \,  dx  \ + \ 
  \iint_{\Lambda_{1}}  x^{\gamma + s} \, 
  \frac{ | D^{\floor{s}} f(x) \, - \, D^{\floor{s}} f(y) |^{2}}{ | x \, - \, y |^{1 \, + \, 2 (s - \floor{s})}} dy \,  dx \, \\
 &:=   | f |_{H^{s}_{(\gamma)}(\Lambda)}^{2}  \ + \  | f |_{H^{s}_{(\gamma)}(\Lambda_{1})}^{2}  \, , \\
 &  \\
\mbox{and } \ 
  \| f \|_{H^{s}_{(\gamma)}(\mrJ)}^{2} &:= \  \left\{ \begin{array}{rl}  
  \sum_{j = 0}^{s} \| D^{j} f \|_{L^{2}_{(\gamma + j)}(\mrJ)}^{2} \, , & \mbox{ for } s \in \mathbb{N}_{0}  \\
 \sum_{j = 0}^{\floor{s}} \| D^{j} f \|_{L^{2}_{(\gamma + j)}(\mrJ)}^{2} \ + \  | f |_{H^{s}_{(\gamma)}(\mrJ)}^{2} \, , & \mbox{ for } s \in \mathbb{R}^{+} \backslash \mathbb{N}_{0} 
 \end{array} \right. \, ,   \label{deffnJ} \\
  &  \\
\mbox{where } \ 
\| g \|_{L^{2}_{(\gamma)}(\mrJ)}^{2} &:= \ 
  \int_{\mrJ}  x^{\gamma} \, g^{2}(x) \, dx \, .  \nonumber
\end{align}

Then, $H^{s}_{(\gamma)}(\mrJ) \, := \, \{ f \, : \, f \mbox{ is measurable and }  \| f \|_{H^{s}_{(\gamma)}(\mrJ)} \, < \, \infty \}$. 

\textbf{Note}: A function $f(x)$ is in $H^{s}_{(a , b)}(\mrI)$ if and only if $f(\frac{3}{4} x) \in H^{s}_{(b)}(\mrJ)$ and
$f(\frac{3}{4}(1 - x)) \in H^{s}_{(a)}(\mrJ)$.

The following lemma and discussion allows us to further focus our analysis for determining the regularity of 
$(1 - x)^{\alpha - \beta} \, x^{\beta} \phi(x)$.

\begin{lemma} \label{lmaoutnrm}
Let $s \ge 0$, $\psi \in H^{s}_{(\gamma)}(\mrJ)$, and $g \in C^{\ceil{s}}(\mrJ)$. Then
\[
 \| g \, \psi \|_{H^{s}_{(\gamma)}(\mrJ)}^{2} \ \ \lesssim \ \| g \|_{C^{\ceil{s}}(\mrJ)}^{2} \, 
 \| \psi \|_{H^{s}_{(\gamma)}(\mrJ)}^{2} \, .
\]
\end{lemma}
\textbf{Proof}: For $s = 0$,
\be
\| g \, \psi \|_{H^{s}_{(\gamma)}(\mrJ)}^{2} \ = \ \| g \, \psi \|_{L^{2}_{(\gamma)}(\mrJ)}^{2}
\ \le \  \| g \|_{L^{\infty}(\mrJ)}^{2} \,  \| \psi \|_{L^{2}_{(\gamma)}(\mrJ)}^{2} \, .
\label{po1}
\ee
Then, for $s = 1$,
\begin{align}
\| g \, \psi \|_{H^{s}_{(\gamma)}(\mrJ)}^{2} 
&\le \ \| g \, \psi \|_{L^{2}_{(\gamma)}(\mrJ)}^{2}  \ + \ \| D(g \, \psi) \|_{L^{2}_{(\gamma + 1)}(\mrJ)}^{2}  \nonumber \\
&\lesssim \ \| g \|_{L^{\infty}(\mrJ)}^{2}  \, \| \psi \|_{L^{2}_{(\gamma)}(\mrJ)}^{2} 
 \ + \ \| \psi \, Dg \|_{L^{2}_{(\gamma + 1)}(\mrJ)}^{2} 
  \ + \ \| g \, D\psi \|_{L^{2}_{(\gamma + 1)}(\mrJ)}^{2}  \nonumber \\
&\lesssim \ \| g \|_{L^{\infty}(\mrJ)}^{2}  \, \| \psi \|_{L^{2}_{(\gamma)}(\mrJ)}^{2} 
 \ + \ \| Dg \|_{L^{\infty}(\mrJ)}^{2}  \, \| \psi \|_{L^{2}_{(\gamma)}(\mrJ)}^{2} 
  \ + \ \| g \|_{L^{\infty}(\mrJ)}^{2}  \, \| D\psi \|_{L^{2}_{(\gamma + 1)}(\mrJ)}^{2}  \nonumber \\
&\lesssim \ \| g \|_{C^{1}(\mrJ)}^{2}  \, \| \psi \|_{H^{1}_{(\gamma)}(\mrJ)}^{2} \, .  \label{po2}
\end{align}
Next, for $0 \, <  \, s \, < \, 1$, consider the mapping 
$\mcF \, : \, H^{s}_{(\gamma)}(\mrJ) \ \longrightarrow \ H^{s}_{(\gamma)}(\mrJ)$, defined by
$\mcF(\psi) \ := \ g \, \psi$. From \eqref{po1} and \eqref{po2} we have that $\mcF$ is a bounded mapping
for $s = 0$ and $s = 1$. As $H^{s}_{(\gamma)}(\mrJ)$ is a family of interpolation spaces, it follows that
$\mcF$ is a bounded mapping for $0 \, <  \, s \, < \, 1$, with $\| \mcF \| \, \lesssim \, \| g \|_{C^{1}(\mrJ) }$.

For $s > 1$ the above argument extends in a straight forward manner. \\
\mbox{ } \hfill \qed

Consider $h(x) \, = \, g(x) \psi(x)$, where $\psi(x) \in H^{s}_{(\mu)}(\mrJ)$ and 
$g(x) \, = \, \left\{ \begin{array}{rl}
     x^{\alpha}, & \ \frac{1}{2} \le x \le \frac{3}{4}   \\
     g_{ext}(x), & \ 0 \le x \le \frac{1}{2} 
     \end{array} \, , \right.$
for  $ g_{ext}(x)$ a $C^{\ceil{s}}(0 \, , \, \frac{1}{2})$ extension of $x^{\alpha}$ satisfying
$\| g \|_{C^{\ceil{s}}(\mrJ)} \ \le \ \| x^{\alpha} \|_{C^{\ceil{s}}[\frac{1}{2}  , \frac{3}{4} ]}$.

Note that for $t \le s$, $t \not \in \mathbb{N}$,
\begin{align*}
\| x^{\alpha} \, \psi(x) \|^{2}_{H^{t}_{(\sigma)}(\mrJ)} 
&\lesssim \ \sum_{j = 0}^{\floor{t}} \| D^{j} (x^{\alpha} \, \psi(x) ) \|_{L^{2}_{(\sigma + j)}(\mrJ)}^{2} \ + \ 
    \iint_{\Lambda}  x^{\sigma + t} \, 
  \frac{ | D^{\floor{t}} (x^{\alpha} \, \psi(x) ) \, - \, D^{\floor{t}} (y^{\alpha} \, \psi(y) ) |^{2}}{ | x \, - \, y |^{1 \, + \, 2 (t - \floor{t})}} dy \,  dx
     \\
&\quad  \quad \quad +     
     \iint_{\Lambda_{1}}  x^{\mu + t} \, 
  \frac{ | D^{\floor{t}} (x^{\alpha} \, \psi(x) ) \, - \, D^{\floor{t}} (y^{\alpha} \, \psi(y) ) |^{2}}{ | x \, - \, y |^{1 \, + \, 2 (t - \floor{t})}} dy \,  dx  \\
&\lesssim \ \sum_{j = 0}^{\floor{t}} \| D^{j} (x^{\alpha} \, \psi(x) ) \|_{L^{2}_{(\sigma + j)}(\mrJ)}^{2} \ + \ 
    \iint_{\Lambda}  x^{\sigma + t} \, 
  \frac{ | D^{\floor{t}} (x^{\alpha} \, \psi(x) ) \, - \, D^{\floor{t}} (y^{\alpha} \, \psi(y) ) |^{2}}{ | x \, - \, y |^{1 \, + \, 2 (t - \floor{t})}} dy \,  dx
     \\
&\quad  \quad \quad + \ \| h \|^{2}_{H^{t}_{(\mu)}(\mrJ)}  \\
&\lesssim 
\sum_{j = 0}^{\floor{t}} \| D^{j} (x^{\alpha} \, \psi(x) ) \|_{L^{2}_{(\sigma + j)}(\mrJ)}^{2} \ + \ 
    \iint_{\Lambda}  x^{\sigma + t} \, 
  \frac{ | D^{\floor{t}} (x^{\alpha} \, \psi(x) ) \, - \, D^{\floor{t}} (y^{\alpha} \, \psi(y) ) |^{2}}{ | x \, - \, y |^{1 \, + \, 2 (t - \floor{t})}} dy \,  dx
   \\
&\quad  \quad \quad + \  \| \psi \|^{2}_{H^{s}_{(\mu)}(\mrJ)} \quad  \mbox{(using Lemma \ref{lmaoutnrm})}.
\end{align*}

Hence for the analysis of the regularity of $f(x) \, = \, (1 - x)^{\alpha - \beta} x^{\beta} \phi(x)$ 
we can restrict our attention on the analysis
of the semi-norm to $ | f |_{H^{t}_{(\sigma)}(\Lambda)}$.

\subsection{Regularity of $x^{p} \psi(x)$}
\label{ssec_reg3}
In this section we establish the value of $t$ and $\sigma$ such that $x^{p} \psi(x) \in H^{t}_{(\sigma)}(\mrJ)$
for $\psi(x) \in H^{s}_{(\mu)}(\mrJ)$.
The general result is given in Theorem \ref{thrmall1}. There are two key terms which arise in the proof of
Theorem \ref{thrmall1}. The analysis for one of these terms follows similarly to a term which occurs for the
case of $s$ between $0$ and $1$, discussed in Theorem \ref{thrm1}. The other of these terms arises
for the case of $s$ between $1$ and $2$, discussed in Theorem \ref{thrm2}. We begin this section with
an embedding theorem which is used in the proofs of the subsequent theorems in this section.

\begin{theorem} \label{imbedL2}
For $s \ge 0$, $\gamma - s \, > \, -1$, then $H_{(\gamma)}^{s}(\mrJ) \subset L^{2}_{(\gamma - s)}(\mrJ)$.
\end{theorem}
\textbf{Proof}: Firstly we consider the case for $0 \le s \le 1$. \\
For $s = 0$ we have that 
\be
H_{(\gamma)}^{s}(\mrJ) \, = \, H_{(\gamma)}^{0}(\mrJ) \, = \, L^{2}_{(\gamma)}(\mrJ). 
\label{gr1}
\ee
For $s = 1$, consider $\phi \in L^{2}_{(\gamma - s)}(\mrJ) \, = \, L^{2}_{(\gamma - 1)}(\mrJ)$. Then,
using Hardy's inequality \cite[Lemma 3.2]{ber071},
\begin{align}
\| \phi \|_{L^{2}_{(\gamma - 1)}}^{2} &= \ \int_{\mrJ} x^{\gamma - 1} \big( \phi(x) \big)^{2} \, dx  \nonumber \\
&\lesssim \ \int_{\mrJ} x^{\gamma + 1} \big( \phi'(x) \big)^{2} \, dx \ + \ \int_{\mrJ} x^{\gamma + 1} \big( \phi(x) \big)^{2} \, dx \nonumber \\
&\le \int_{\mrJ} x^{\gamma + 1} \big( \phi'(x) \big)^{2} \, dx \ + \ \int_{\mrJ} x^{\gamma} \big( \phi(x) \big)^{2} \, dx \nonumber \\
&= \ \| \phi \|_{H^{1}_{(\gamma)}}^{2}  \, .   \label{gr2}
\end{align}

From \eqref{gr1} and \eqref{gr2} the identity operator $\mathbb{I}$ mapping from 
$H_{(\gamma)}^{s}(\mrJ) \longrightarrow L^{2}_{(\gamma - s)}(\mrJ)$, $s \, = \, 0, 1$ is a bounded operator.

Additionally, the spaces $L^{2}_{\sigma}(\mrJ)$ are interpolation spaces \cite[Lemma 23.1]{tar071}.

Hence, for $\theta \ = \  (1 - \theta) 0 \, + \, \theta 1$, using the fact that $H_{(\gamma)}^{\theta}(\mrJ)$ and
$L^{2}_{(\theta)}(\mrJ)$ are interpolation spaces, it follows that
\[
   H_{(\gamma)}^{\theta}(\mrJ) \stackrel{\mathbb{I}}{\longrightarrow} L^{2}_{( (1 - \theta) \gamma \, + \, \theta (\gamma - 1) )}(\mrJ)
   \ = \ L^{2}_{(\gamma - \theta)}(\mrJ) 
\]
is bounded. Thus, if $u \in  H_{(\gamma)}^{\theta}(\mrJ)$ then $u \in L^{2}_{(\gamma - \theta)}(\mrJ)$ with 
$\| u \|_{L^{2}_{(\gamma - \theta)}} \, \le \, C \, \| u \|_{H_{(\gamma)}^{\theta}}$.

Next, for $1 \le s \le 2$, consider $s = 2$. For $\phi \in L^{2}_{(\gamma - s)}(\mrJ) \, = \, L^{2}_{(\gamma - 2)}(\mrJ)$, again
using Hardy's inequality (and that $\gamma - s > -1$),
\begin{align}
\| \phi \|_{L^{2}_{(\gamma - 2)}}^{2} &= \ \int_{\mrJ} x^{\gamma - 2} \big( \phi(x) \big)^{2} \, dx  \nonumber   \\
&\lesssim \ \int_{\mrJ} x^{\gamma} \big( \phi'(x) \big)^{2} \, dx \ + 
\ \int_{\mrJ} x^{\gamma} \big( \phi(x) \big)^{2} \, dx \nonumber   \\
&\lesssim \int_{\mrJ} x^{\gamma + 2} \big( \phi''(x) \big)^{2} \, dx \ + \ \int_{\mrJ} x^{\gamma + 2} \big( \phi'(x) \big)^{2} \, dx \ + \ 
\int_{\mrJ} x^{\gamma} \big( \phi(x) \big)^{2} \, dx \nonumber   \\
&\le \int_{\mrJ} x^{\gamma + 2} \big( \phi''(x) \big)^{2} \, dx \ + \ \int_{\mrJ} x^{\gamma + 1} \big( \phi'(x) \big)^{2} \, dx \ + \ 
\int_{\mrJ} x^{\gamma} \big( \phi(x) \big)^{2} \, dx \nonumber   \\
&= \ \| \phi \|_{H^{2}_{(\gamma)}}^{2}  \, .   \label{gr3}
\end{align}

Again, using the fact that $H_{(\gamma)}^{\theta}(\mrJ)$ and
$L^{2}_{(\theta)}(\mrJ)$ are interpolation spaces, it now follows that for $0 \le \theta \le 2$
if $u \in  H_{(\gamma)}^{\theta}(\mrJ)$ then $u \in L^{2}_{(\gamma - \theta)}(\mrJ)$ with 
$\| u \|_{L^{2}_{(\gamma - \theta)}} \, \le \, C \, \| u \|_{H_{(\gamma)}^{\theta}}$.

The argument extends in an obvious manner to arbitrary $s > 2$. \\
\mbox{ } \hfill \qed

%%%%%
%%%%%
\begin{theorem} \label{thrm1}
Let $0 \le s < 1$, $\mu > -1$, and $\psi \in H^{s}_{(\mu)}(\mrJ)$. Then $x^{p} \, \psi \in H^{t}_{(\sigma)}(\mrJ)$
provided
\be
 0 \le  t \le s \, , \ \ \sigma \, + \, 2 p \, \ge  \, \mu  \, , \  \ \sigma \, + \, 2 p \, - t \, >  \,  - 1 \, , \ 
   \ \mbox{ and } \ \ \sigma \, + \, 2 p \, + t \, \ge \, \mu \, + \, s   \, .
 \label{piu1}
\ee
Additionally, when \eqref{piu1} is satisfied, there exists $C > 0$ (independent of $\psi$) such that
$\| x^{p} \, \psi \|_{H^{t}_{(\sigma)}(\mrJ)} \, \le \, C \, \| \psi \|_{H^{s}_{(\mu)}(\mrJ)}$.
\end{theorem} 
\textbf{Proof}: Firstly, for $s \, = \, t \, = \, 0$,
\begin{align}
\| x^{p} \, \psi \|_{H^{t}_{(\sigma)}}^{2} \ = \ \| x^{p} \, \psi \|_{H^{0}_{(\sigma)}}^{2} 
&= \ \int_{\mrJ} \, x^{\sigma} \, \left( x^{p} \, \psi(x) \right)^{2} \, dx
 \ = \ \int_{\mrJ} \, x^{\sigma \, + \, 2 p} \, \left( \psi(x) \right)^{2} \, dx  \label{piu2} \\
&\le \  \int_{\mrJ} \, x^{\mu} \, \left( \psi(x) \right)^{2} \, dx  \, , \mbox{ provided } \ \sigma \, + \, 2 p \ge \mu \, , \nonumber \\
&= \ \| \psi \|_{H^{0}_{(\mu)}}^{2} \, = \,  \| \psi \|_{H^{s}_{(\mu)}}^{2} \, . \nonumber
\end{align}

For $0 < s < 1$, in addition to \eqref{piu2} we must also consider the semi-norm  $|  x^{p} \, \psi |_{H^{t}_{(\sigma)}(\Lambda)}$.
\begin{align}
|  x^{p} \, \psi |_{H^{t}_{(\sigma)}(\Lambda)}^{2} 
&= \ \iint_{\Lambda} \, x^{\sigma + t} \, 
\frac{ | x^{p} \, \psi(x) \ - \ y^{p} \, \psi(y) |^{2}} { | x \ - \ y |^{1 \, + \, 2 t}} \, dy \, dx    \label{piu3}    \\
&\lesssim \  \iint_{\Lambda} \, x^{\sigma + t} \, y^{2 p} \, 
\frac{ | \psi(x) \ - \ \psi(y)  |^{2}} { | x \ - \ y |^{1 \, + \, 2 t}} \, dy \, dx  \ 
+ \ \iint_{\Lambda} \, x^{\sigma + t} \, \psi^{2}(x) \, 
\frac{ |  x^{p}  \ - \ y^{p}  |^{2}} { | x \ - \ y |^{1 \, + \, 2 t}} \, dy \, dx    \nonumber \\
&:= \ I_{1} \ + \ I_{2} \, . \nonumber 
\end{align}

Noting that in $\Lambda$, $y  \, <  \, \frac{3}{2} x$, for $I_{1}$ we have
\begin{align}
I_{1} &\lesssim \  \iint_{\Lambda} \, x^{\sigma + t} \, x^{2 p} \, 
\frac{ | \psi(x) \ - \ \psi(y)  |^{2}} { | x \ - \ y |^{1 \, + \, 2 t}} \, dy \, dx   \nonumber  \\
&\le \  \iint_{\Lambda} \, x^{\mu + s} \, 
\frac{ | \psi(x) \ - \ \psi(y)  |^{2}} { | x \ - \ y |^{1 \, + \, 2 s}} \, dy \, dx \, ,
\ \  \mbox{ provided } \ \sigma  + t \, + \, 2 p \ge \mu + s \, ,  \mbox{ and } \ t \le s \, , \nonumber \\
&= \ | \psi |^{2}_{H^{s}_{(\mu)}(\Lambda)} \, .    \label{piu4}
\end{align}

If $p = 0$ then $I_{2} = 0$.
To bound $I_{2}$ for $p \ne 0$
we introduce the change of variable: $y \, = \, \frac{1}{z} x$, where $\frac{2}{3} < z < \frac{3}{2}$.
With this change of variable, we have
\begin{align}
dy &= \frac{-1}{z^{2}} x \, dz \, , \quad  |x - y | \ = \  x \, \frac{1}{z} |z - 1| \, , \   \mbox{ and for } I_{2} \nonumber \\ 
I_{2} &= \ \int_{0}^{1/2} \, x^{\sigma \, + \, 2 p \, - t } \psi^{2}(x) \, dx \
   \int_{2/3}^{3/2} \, z^{-1 \, + \, 2 t \, - \, 2 p} \, | 1 - z |^{-1 \, - \, 2 t} |z^{p}  - 1|^{2} \, dz  \label{piu4p5} \\
&\lesssim \    \| \psi \|^{2}_{L^{2}_{(\sigma  \, + \, 2 p \, - t)}} \ \cdot \ 
 \int_{2/3}^{3/2} \,  | 1 - z |^{-1 \, - \, 2 t} |z^{p}  - 1|^{2} \, dz\, .    \label{piu5}  
\end{align}

Using Hardy's inequality \cite[Lemma 3.2]{ber071}, we bound the integral in \eqref{piu5} as follows.
\begin{align}
 \int_{2/3}^{3/2} \,  | 1 - z |^{-1 \, - \, 2 t} |z^{p}  - 1|^{2} \, dz
 &= \  \int_{2/3}^{1} \,  ( 1 - z )^{-1 \, - \, 2 t} (z^{p}  - 1)^{2} \, dz \ + \ 
  \int_{1}^{3/2} \,  ( z - 1 )^{-1 \, - \, 2 t} (z^{p}  - 1)^{2} \, dz   \nonumber \\
&\lesssim \  \int_{2/3}^{1} \,  ( 1 - z )^{-1 \, - \, 2 t + 2} ( z^{p - 1} )^{2} \, dz \ + \ 
  \int_{1}^{3/2} \,  ( z - 1 )^{-1 \, - \, 2 t + 2} ( z^{p - 1} )^{2} \, dz  \nonumber  \\
&\lesssim 1 \, , \quad \ \mbox{ provided } t < 1 \, .    \label{piu6}
\end{align}

From Theorem \ref{imbedL2}, we have
\begin{align}
\| \psi \|_{L^{2}_{(\sigma  \, + \, 2 p \, - t)}}  &\lesssim \ \| \psi \|_{H^{t}_{(\sigma \, + \, 2 p)}}  , 
\quad \ \mbox{ provided } \sigma \, + \, 2 p  \, - t  \, >  \, -1 \, ,  \nonumber \\
&\le \ \| \psi \|_{H^{s}_{(\mu)}}  , 
\quad \ \mbox{ provided } \sigma \, + \, 2 p   \ge  \mu \, ,  \ \mbox{ and } t \le s \, .  \label{piu7}
\end{align}
Finally, combining \eqref{piu2}-\eqref{piu7} we obtain the stated results. \\
\mbox{ } \hfill \qed

%%%%
%%%%
\begin{theorem} \label{thrm2}
Let $1 \le s < 2$, $\mu > -1$, and $\psi \in H^{s}_{(\mu)}(\mrJ)$. 
Then $x^{p} \, \psi \in H^{t}_{(\sigma)}(\mrJ)$
provided
\be
  0 \le t \le s \, , \ \ \sigma \, + \, 2 p \, \ge  \, \mu  \, , \  \ \sigma \, + \, 2 p \, - t \, >  \,  - 1 \, , \ 
   \ \mbox{ and } \ \ \sigma \, + \, 2 p \, + t \, \ge \, \mu \, + \, s   \, .
 \label{piu11}
\ee
Additionally, when \eqref{piu11} is satisfied, there exists $C > 0$ (independent of $\psi$) such that
$\| x^{p} \, \psi \|_{H^{t}_{(\sigma)}(\mrJ)} \, \le \, C \, \| \psi \|_{H^{s}_{(\mu)}(\mrJ)}$.
\end{theorem} 
\textbf{Proof}: For $0 \le t < 1$ Theorem \ref{thrm1} applies. We assume that $t \ge 1$. Hence, of interest
is
\[
  \| x^{p} \, \psi \|_{L^{2}_{(\sigma)}} , \quad \ \| D \left(  x^{p} \, \psi \right) \|_{L^{2}_{(\sigma + 1)}} 
  , \quad \ | D \left(  x^{p} \, \psi \right) |_{H^{t - 1}_{(\sigma + 1)}(\Lambda)}  \, .
\]
From Theorem \ref{thrm1}, 
\be
\| x^{p} \, \psi \|_{L^{2}_{(\sigma)}} \lesssim \| \psi \|_{H_{(\mu)}^{0}} \, ,  \ \ \mbox{ provided } \ \
\sigma \, + \, 2 p  \,  \ge  \, \mu \, .
\label{piu12}
\ee

For $D \left(  x^{p} \, \psi \right)$ we have
\begin{align*}
D \left(  x^{p} \, \psi \right) &= \ x^{p} \, \psi' \ + \   p \, x^{p - 1} \, \psi \ \mbox{ and } \\
\left( D \left(  x^{p} \, \psi \right) \right)^{2} &\lesssim \ \left( x^{p} \, \psi' \right)^{2} \ + \ \left( x^{p - 1} \, \psi  \right)^{2} \, .
\end{align*}

Then, using Hardy's inequality \cite[Lemma 3.2]{ber071} (using $\sigma \, + \, 2 p \, - t > -1$, i.e., $\sigma  \, + \, 2 p  > \, 0$),
\begin{align}
\| D \left(  x^{p} \, \psi \right) \|_{L^{2}_{(\sigma + 1)}}^{2} 
&\lesssim \  \int_{0}^{3/4} \, x^{\sigma + 1} \,  \left( x^{p} \, \psi'(x) \right)^{2} \, dx \ + \ 
\int_{0}^{3/4} \, x^{\sigma + 1} \,  \left( x^{p - 1} \, \psi(x)  \right)^{2} \, dx  \nonumber \\
&= \  \int_{0}^{3/4} \, x^{\sigma  \, + \, 2 p \, + 1} \,  \left( \psi'(x) \right)^{2} \, dx \ + \ 
\int_{0}^{3/4} \, x^{\sigma  \, + \, 2 p \, - 1 } \,  \left( \psi(x) \right)^{2} \, dx  \nonumber \\
&\lesssim \  \int_{0}^{3/4} \, x^{\sigma  \, + \, 2 p \,  + 1 } \,  \left( \psi'(x) \right)^{2} \, dx   \nonumber  \\
 & \quad \quad \quad  \ + \ 
\int_{0}^{3/4} \, x^{\sigma  \, + \, 2 p \, + 1} \,  \left( \psi'(x) \right)^{2} \, dx \ + \ 
\int_{0}^{3/4} \, x^{\sigma  \, + \, 2 p \, + 1 } \,  \left( \psi(x) \right)^{2} \, dx  \nonumber \\
&\lesssim \ \| \psi \|_{H^{1}_{(\mu)}}^{2} \, , \ \ \mbox{ provided } \ \ \sigma \, + \, 2 p \, \ge  \, \mu  \, . \label{piu13}
\end{align}

Equation \eqref{piu13}, together with \eqref{piu12} establishes the stated result for $t = s = 1$.

Next, for $| D \left(  x^{p} \, \psi \right) |_{H^{t - 1}_{(\sigma + 1)}(\Lambda)} $ we have
\begin{align}
& | D \left(  x^{p} \, \psi \right) |_{H^{t - 1}_{(\sigma + 1)}(\Lambda)}^{2}
\ = \ \iint_{\Lambda} \, x^{\sigma + t} \, 
\frac{ | D \left(  x^{p} \, \psi(x) \right) \, - \, D \left(  y^{p} \, \psi(y) \right) |^{2}}%
{ | x \, - \, y |^{1 \, + \, 2 (t - 1)} } \ dy \, dx  \nonumber \\
& \quad \quad \quad \lesssim \  \iint_{\Lambda} \, x^{\sigma + t} \, 
\frac{ | x^{p - 1} \, \psi(x)  \ - \  y^{p - 1} \, \psi(y) |^{2}}%
{ | x \, - \, y |^{1 \, + \, 2 (t - 1)} } \ dy \, dx \
+ \   \iint_{\Lambda} \, x^{\sigma + t} \, 
\frac{ | x^{p} \, \psi'(x)  \ - \  y^{p} \, \psi'(y) |^{2}}%
{ | x \, - \, y |^{1 \, + \, 2 (t - 1)} } \ dy \, dx \  \nonumber \\
&  \quad \quad \quad := \ I_{3}  \ + \ I_{4} \, .    \label{piu14}
\end{align}

For $I_{3}$, proceeding as in \eqref{piu3},
\begin{align}
I_{3} &\lesssim \  \iint_{\Lambda} \, x^{\sigma + t} \, y^{2 (p - 1)} \, 
\frac{ | \psi(x) \ - \ \psi(y)  |^{2}} { | x \ - \ y |^{1 \, + \, 2 (t - 1)}} \, dy \, dx  \ 
+ \ \iint_{\Lambda} \, x^{\sigma + t} \, \psi^{2}(x) \, 
\frac{ |  x^{p - 1}  \ - \ y^{p - 1}  |^{2}} { | x \ - \ y |^{1 \, + \, 2 (t - 1)}} \, dy \, dx    \nonumber \\
&:= \ I_{3 , 1} \ + \ I_{3 , 2} \, . \nonumber 
\end{align}

For $I_{3 , 1}$, using $y \, < \, \frac{3}{2} x$, and the introducing the change of variable 
$y \, = \, \frac{1}{z} x$, where $\frac{2}{3} < z < \frac{3}{2}$ we obtain
\begin{align}
I_{3 , 1} &\lesssim \int_{x = 0}^{1/2} \, x^{\sigma \, + \, t \, + \, 2 ( p - 1)} \, \int_{z = 2/3}^{3/2} \,
\frac{z^{1 \, + \, 2 (t - 1)  \, - \, 2 (p -1)}}{ ( x \, | z - 1| )^{1 \, + \, 2 (t - 1)}}   \, 
| \psi(x) \, - \, \psi(\frac{1}{z} x) |^{2} \, z^{-2} \, x \, dz \, dx  \nonumber  \\
&\lesssim \  \int_{x = 0}^{1/2} \, x^{\sigma \, + \, 2 p \, - \, t} 
\bigg( \int_{z = 2/3}^{1} \, (1  - z)^{1 \, - \, 2 t}  \, ( \psi(x) \, - \, \psi(\frac{1}{z} x) )^{2} \, dz  \nonumber \\
 & \quad \quad \quad \quad \quad  \quad \quad \quad \quad \quad \quad \quad \quad  \ + 
  \int_{z = 1}^{3/2} \, (z  - 1)^{1 \, - \, 2 t}  \, ( \psi(x) \, - \, \psi(\frac{1}{z} x) )^{2} \, dz \bigg) dx  \nonumber \\
&\lesssim \   \int_{x = 0}^{1/2} \, x^{\sigma \, + \, 2 p \, - \, t} 
\bigg( \int_{z = 2/3}^{1} \, (1  - z)^{3 \, - \, 2 t}  \, z^{-4} \, x^{2} \, ( \psi'(\frac{1}{z} x) )^{2} \, dz  \nonumber \\
 & \quad \quad \quad \quad \quad  \quad \quad \quad \quad \quad \quad \quad \quad  \ + 
  \int_{z = 1}^{3/2} \, (z  - 1)^{3 \, - \, 2 t}  \, z^{-4} \, x^{2} \, ( \psi'(\frac{1}{z} x) )^{2}  \, dz \bigg) dx  \nonumber \\
& \quad \quad \quad \quad \mbox{(using Hardy's inequality)}   \nonumber \\
&\lesssim \  \int_{x = 0}^{1/2} \, x^{\sigma \, + \, 2 p \, - \, t \, + \, 2} \, 
\int_{z = 2/3}^{3/2} | 1 - z |^{3 \, - \, 2 t} \, ( \psi'(\frac{1}{z} x) )^{2}  \, dz \, dx  \, . \label{piu141}
\end{align}

Next, letting $w \, = \, \frac{1}{z} x$, then $dw \, = \,  \frac{1}{z} \, dx$ and
\begin{align}
I_{3 , 1} &\lesssim \ 
 \int_{z = 2/3}^{3/2} | 1 - z |^{3 \, - \, 2 t} \,  \int_{w = 0}^{1 \, / \, 2 z} \, 
 (z \, w)^{\sigma \, + \, 2 p \, - \, t \, + \, 2} \, ( \psi'(w) )^{2}  \, z \,  dw \, dz  \nonumber \\
&\lesssim \ 
 \int_{2/3}^{3/2} | 1 - z |^{3 \, - \, 2 t} \, dz \,  \int_{0}^{3/4} \, 
 w^{\sigma \, + \, 2 p \, - \, t \, + \, 2} \, ( \psi'(w) )^{2}  \,  dw   \nonumber \\
&\lesssim \ \| \psi' \|^{2}_{L^{2}_{(\sigma \, + \, 2 p \, - \, t \, + \, 2)}}  \, , \mbox{ as } t < 2 \, ,  \nonumber \\
&\lesssim \ \| \psi' \|^{2}_{H^{t - 1}_{(\sigma \, + \, 2 p \, + \, 1)}}  \, , \ \mbox{(using Theorem \ref{imbedL2})}
 \ \ \mbox{provided } \ \sigma + \, 2 p \, - \, t \, \, + \, 2 \, > \, -1 \, , \nonumber \\
&\lesssim \ \| \psi' \|^{2}_{H^{s - 1}_{(\mu \, + \, 1)}} \, , \ \ 
\mbox{provided } \ \sigma \, + \, 2 p \, \ge \, \mu \, , \mbox{ and } t \le s \, ,  \nonumber \\
&\lesssim \ \| \psi \|^{2}_{H^{s}_{(\mu)}} \, .     \label{piu142}
\end{align}

If $p = 1$ then $I_{3 , 2} = 0$.
For $I_{3 , 2}$ with $p \neq 1$, proceeding as in the approach used to obtain the bound for $I_{2}$, 
\eqref{piu4p5} - \eqref{piu6},
\begin{align}
I_{3 , 2} &\lesssim \
\int_{0}^{1/2} x^{(\sigma + 1) \, - \, (t - 1) \, + \, 2 (p - 1)} ( \psi(x) )^{2} \, dx \
\int_{2/3}^{3/2} \, z^{-1 \, + \, 2 (t - 1) \, - \, 2 (p - 1)} \, |1 - z|^{-1 \, - \, 2 ( t - 1)} \, |z^{p - 1} \, - \, 1|^{2} \, dz \nonumber \\
&\lesssim \ \int_{x = 0}^{1/2} x^{\sigma \, + \, 2 p \, - \, t} ( \psi(x) )^{2} \, dx \, , \
\mbox{ provided } (t  - 1) \, <  \, 1 \, ,  \mbox{ i.e., } t < 2 \, ,  \nonumber \\
&\lesssim \ \int_{0}^{1/2} x^{\sigma \, + \, 2 p \, - \, t \, + \, 2} ( \psi'(x) )^{2} \, dx \ + \ 
 \int_{0}^{1/2} x^{\sigma \, + \, 2 p \, - \, t \, + \, 2} ( \psi(x) )^{2} \, dx   \, , \nonumber \\
 & \quad \quad 
 \ \ \mbox{(using Hardy's inequality)}  \mbox{ provided } \sigma \, + \, 2 p \, - \, t   \, > \, -1 \, ,  \nonumber \\
&\lesssim \ \| \psi' \|^{2}_{L^{2}_{(\sigma \, + \, 2 p \, - \, t \, + \, 2)}} \ 
+ \ \| \psi \|^{2}_{L^{2}_{(\sigma \, + \, 2 p \, - \, t \, + \, 2)}}  \nonumber \\
&\lesssim \ \| \psi' \|^{2}_{H^{t - 1}_{(\sigma \, + \, 2 p \, + \, 1)}} \ 
+ \ \| \psi \|^{2}_{H^{t}_{(\sigma \, + \, 2 p \, + \, 2)}}  \, , \ \mbox{(using Theorem \ref{imbedL2})}
\mbox{ provided } \sigma \, + \, 2 p  \, - \, t   \, + \, 2 \, > \, -1 \, ,   \nonumber \\
&\lesssim \ \| \psi' \|^{2}_{H^{s - 1}_{(\mu \, + \, 1)}} \ 
+ \ \| \psi \|^{2}_{H^{s}_{(\mu \, + \, 2)}}   \, , \ \mbox{ provided } \sigma \, + \, 2 p \, \ge \, \mu \, 
\mbox{ and } t \le s,   \nonumber \\
&\lesssim  \| \psi \|^{2}_{H^{s}_{(\mu)}}  \, .   \label{piu144}
\end{align}

With $I_{4}$, proceeding as in \eqref{piu3},
\begin{align}
I_{4} &\lesssim  \  \iint_{\Lambda} \, x^{\sigma + t} \, y^{2 p } \, 
\frac{ | \psi'(x) \ - \ \psi'(y)  |^{2}} { | x \ - \ y |^{1 \, + \, 2 (t - 1)}} \, dy \, dx  \ 
+ \ \iint_{\Lambda} \, x^{\sigma + t} \, ( \psi'(x) )^{2}  \, 
\frac{ |  x^{p }  \ - \ y^{p}  |^{2}} { | x \ - \ y |^{1 \, + \, 2 (t - 1)}} \, dy \, dx    \nonumber \\
&\lesssim \ | \psi' |^{2}_{H^{s - 1}_{(\mu + 1)}(\Lambda)} 
\,  \  \mbox{ (provided } \ \sigma  + \, 2 p \, + t  \ge \mu + s \, ,  \mbox{ and } \ t \le s \,)   \nonumber \\
& \quad \ + \ \| \psi' \|^{2}_{H^{s - 1}_{(\mu + 1)}} 
\  \ ( \mbox{ provided  } \sigma \, + \, 2 p  \, - t  + 2 \, >  \, -1 \, ,
  \sigma \, + \, 2 p   \ge  \mu \, ,  \ \mbox{ and } 1 < t \le s \, )    \nonumber \\
&\lesssim \  \| \psi \|^{2}_{H^{s}_{(\mu)}}  \, . \label{piu16}
\end{align}

Combining \eqref{piu12} - \eqref{piu16} we obtain the stated result. \\
\mbox{ } \hfill \qed

%%%%
%%%%
We are now in a position to state the general result.
\begin{theorem} \label{thrmall1}
Let $s \ge 0$, $\mu > -1$, and $\psi \in H^{s}_{(\mu)}(\mrJ)$. 
Then $x^{p} \, \psi \in H^{t}_{(\sigma)}(\mrJ)$
provided
\be
 0 \le  t \le s \, , \ \ \sigma \, + \, 2 p \, \ge  \, \mu  \, , \  \ \sigma \, + \, 2 p \, - t \, >  \,  - 1 \, , \ 
   \ \mbox{ and } \ \ \sigma \, + \, 2 p \, + t \, \ge \, \mu \, + \, s   \, .
 \label{piu101}
\ee
Additionally, when \eqref{piu101} is satisfied, there exists $C > 0$ (independent of $\psi$) such that
$\| x^{p} \, \psi \|_{H^{t}_{(\sigma)}(\mrJ)} \, \le \, C \, \| \psi \|_{H^{s}_{(\mu)}(\mrJ)}$.
\end{theorem}
\textbf{Proof}: The proof is an induction argument, using Theorem \ref{thrm2} as the initial step. 

For $t = s = n$ we have,
\begin{align}
\| D^{t} ( x^{p} \, \psi(x) ) \|^{2}_{L^{2}_{(\sigma + t)}} 
&\lesssim \ \sum_{j = 0}^{n} \| x^{p - j} \, D^{n - j} \psi(x) \|^{2}_{L^{2}_{(\sigma + n)}} \nonumber \\
&= \ \sum_{j = 0}^{n} \, \int_{0}^{3/4} \, x^{\sigma + n} \, \left( x^{p - j} \, D^{n - j} \psi(x) \right)^{2} \, dx  \nonumber \\
&= \ \sum_{j = 0}^{n} \, \int_{0}^{3/4} \, x^{\sigma + n \, + \, 2 (p - j)} \, \left(  D^{n - j} \psi(x) \right)^{2} \, dx  \nonumber \\
&  \mbox{ (applying Hardy's inequality $j$ times, using $\sigma \, + \, 2 p \, - \, n \, > \, -1$) }  \nonumber  \\
&\lesssim \ \sum_{j = 0}^{n}  \, \sum_{k = 0}^{j} \, 
 \int_{0}^{3/4} \, x^{\sigma \, + \, 2 p  \, + n} \, \left(  D^{n - j + k} \psi(x) \right)^{2} \, dx  \nonumber \\
 &\lesssim \ \sum_{j = 0}^{n} \| D^{j} \psi \|^{2}_{L^{2}_{(\sigma \, + \, 2 p  \, + n)}}  \nonumber \\
 &\lesssim \ \sum_{j = 0}^{n} \| D^{j} \psi \|^{2}_{L^{2}_{(\mu + n)}}  \ \ 
 \mbox{ (using $\sigma \, + \, 2 p \, \ge \, \mu$) } \nonumber \\
&\lesssim \ \sum_{j = 0}^{n} \| D^{j} \psi \|^{2}_{L^{2}_{(\mu + j)}} 
\  \lesssim \ \| \psi \|^{2}_{H^{n}_{(\mu)}} \ = \ \| \psi \|^{2}_{H^{s}_{(\mu)}} \, .  \label{piu102}
\end{align}
Equation \eqref{piu102}, together with the induction assumption establishes the result for  $t = s = n$ .

For $n \, < \, t , s \, < \, n + 1$ we also need to consider $| D^{n}  ( x^{p} \, \psi(x) ) |_{H^{t - n}_{(\sigma + n)}(\Lambda)}$.
\begin{align}
| D^{n}  ( x^{p} \, \psi(x) ) |^{2}_{H^{t - n}_{(\sigma + t)}(\Lambda)}
&= \ \sum_{j = 0}^{n} \iint_{\Lambda} x^{\sigma + t} 
\frac{| x^{p - j} \, D^{n - j} \psi(x) \ - \ y^{p - j} \, D^{n - j} \psi(y) |^{2}}%
{ | x \, - \, y |^{1 \, + \, 2 (t - n)}} \, dy \, dx  \nonumber \\
&\lesssim \  \sum_{j = 0}^{n} \, 
\left(  \iint_{\Lambda} x^{\sigma + t} \, y^{2 (p - j)} \,
\frac{|  D^{n - j} \psi(x) \ - \ D^{n - j} \psi(y) |^{2}}%
{ | x \, - \, y |^{1 \, + \, 2 (t - n)}} \, dy \, dx  \right.  \nonumber \\
& \quad \quad \quad \left. \ + \ 
 \iint_{\Lambda} x^{\sigma + t} \, \left( D^{n - j} \psi(x) \right)^{2} \, 
\frac{| x^{p - j}  \ - \ y^{p - j}  |^{2}}%
{ | x \, - \, y |^{1 \, + \, 2 (t - n)}} \, dy \, dx    \right)  \nonumber \\
&= \ \iint_{\Lambda} x^{\sigma + t} \, y^{2 p} \,
\frac{|  D^{n} \psi(x) \ - \ D^{n} \psi(y) |^{2}}%
{ | x \, - \, y |^{1 \, + \, 2 (t - n)}} \, dy \, dx   \nonumber  \\
& \quad \ + \ 
\sum_{j = 1}^{n} \, 
 \iint_{\Lambda} x^{\sigma + t} \, y^{2 (p - j)} \,
\frac{|  D^{n - j} \psi(x) \ - \ D^{n - j} \psi(y) |^{2}}%
{ | x \, - \, y |^{1 \, + \, 2 (t - n)}} \, dy \, dx  \nonumber \\
& \quad \ + \  \sum_{j = 0}^{n} \, 
 \iint_{\Lambda} x^{\sigma + t} \, \left( D^{n - j} \psi(x) \right)^{2} \, 
\frac{| x^{p - j}  \ - \ y^{p - j}  |^{2}}%
{ | x \, - \, y |^{1 \, + \, 2 (t - n)}} \, dy \, dx    \, .  \label{piu103}
\end{align}

The first term in \eqref{piu103} is bounded in a similar manner to $I_{1}$ in Theorem \ref{thrm1}, with 
$\sigma \, \rightarrow \, \sigma + n$, $t \, \rightarrow \, t - n$, to obtain
\be
\iint_{\Lambda} x^{\sigma + t} \, y^{2 p} \,
\frac{|  D^{n} \psi(x) \ - \ D^{n} \psi(y) |^{2}}%
{ | x \, - \, y |^{1 \, + \, 2 (t - n)}} \, dy \, dx  
\ \lesssim \ | \psi |^{2}_{H^{s}_{(\mu)}(\Lambda)} \, ,
 \ \mbox{ provided } \ \sigma  \, + \, 2 p \, + t  \ge \mu + s \, ,  \mbox{ and } \ t \le s  \, .
 \label{piu104} 
\ee

For the second term in \eqref{piu103} the terms in the summation are bounded in a similar manner to
$I_{3 , 1}$ in Theorem \ref{thrm2}, 
$\sigma \, \rightarrow \, \sigma + n - 1$, $t \, \rightarrow \, t - n + 1$, $p \, \rightarrow \, p + 1 - j$, to obtain
\begin{align}
& \iint_{\Lambda} x^{\sigma + t} \, y^{2 (p - j)} \,
\frac{|  D^{n - j} \psi(x) \ - \ D^{n - j} \psi(y) |^{2}}%
{ | x \, - \, y |^{1 \, + \, 2 (t - n)}} \, dy \, dx   \nonumber \\
& \quad \lesssim \ \| D^{n - j + 1} \psi \|^{2}_{L^{2}_{(\sigma  \, + \, 2 p \,  + n - j + 1 - (t - n + j - 1))}} \nonumber  \\
& \quad \lesssim \  \| D^{n - j + 1} \psi \|^{2}_{H^{t - (n - j + 1)}_{(\sigma  \, + \, 2 p \, + n - j + 1)}} \, , \
\mbox{ provided }  \sigma \, + \, 2 p \, + n - j + 1 - (t - n + j - 1)  \ > \ -1  \nonumber \\
& \quad \lesssim \  \| D^{n - j + 1} \psi \|^{2}_{H^{s - (n - j + 1)}_{(\mu + n - j + 1)}} \, ,  \
\mbox{ provided } \sigma \, + \, 2 p \, \ge \, \mu \, 
\mbox{ and } t \le s,  \nonumber \\
& \quad \lesssim \ \| \psi \|^{2}_{H^{s}_{(\mu)}} \, .  \label{piu105}
\end{align}

For the third term in \eqref{piu103} the terms in the summation are bounded in a similar manner to
$I_{2}$ in Theorem \ref{thrm1}, 
$\sigma \, \rightarrow \, \sigma + n$, $t \, \rightarrow \, t - n$, $p \, \rightarrow \, p - j$, to obtain
\begin{align}
& \iint_{\Lambda} x^{\sigma + t} \, \left( D^{n - j} \psi(x) \right)^{2} \, 
\frac{| x^{p - j}  \ - \ y^{p - j}  |^{2}}%
{ | x \, - \, y |^{1 \, + \, 2 (t - n)}} \, dy \, dx   \nonumber   \\
& \quad \lesssim \ \| D^{n - j} \psi \|^{2}_{L^{2}_{(\sigma  \, + \, 2 (p - j) \, + n - (t - n) )}} 
\ = \ \| D^{n - j} \psi \|^{2}_{L^{2}_{(\sigma  \, + \, 2 p  \, + n - j - (t - n + j))}}  \, , \
\ \mbox{ provided } t \,  < \,  n + 1 ,  \nonumber \\
&\quad \lesssim \ \| D^{n - j} \psi \|^{2}_{H^{t - (n - j)}_{(\sigma \, + \, 2 p  \, + n - j)}}   \, , \ 
\mbox{ provided }  \sigma \, + \, 2 p \, + n - j - (t - n + j)  \ > \ -1  \nonumber \\
&\quad \lesssim \ \| D^{n - j} \psi \|^{2}_{H^{s- (n - j)}_{(\mu + n - j )}}   \, , \
\mbox{ provided } \sigma \, + \, 2 p \, \ge \, \mu \, 
\mbox{ and } t \le s,  \nonumber \\
& \quad \lesssim \ \| \psi \|^{2}_{H^{s}_{(\mu)}} \, .  \label{piu106}
\end{align}

Combining \eqref{piu102} - \eqref{piu106} the stated result follows. \\
\mbox{ } \hfill \qed

 \setcounter{equation}{0}
\setcounter{figure}{0}
\setcounter{table}{0}
\setcounter{theorem}{0}
\setcounter{lemma}{0}
\setcounter{corollary}{0}
\setcounter{definition}{0}
\section{Regularity of the solution to the fractional diffusion, advection, reaction equation \eqref{DefProb2}}
\label{secRFDAR2}

In this section we present the analysis for the regularity of the solution to \eqref{DefProb2}. The solution's regularity
can be influenced by the regularity of the coefficients $b(x)$ and $c(x)$. We begin with two lemmas which enables
us to insulate the influence of these terms.

%%%%%
%%%%%

Introduce the space $W^{k , \infty}_{w}(\mrI)$ and its associated norm, defined for $k \in \mathbb{N}_{0}$, as
\begin{align}
W^{k , \infty}_{w}(\mrI) &:= \ \left\{ f \, : \ (1 - x)^{j/2} x^{j/2} D^{j}f(x) \in L^{\infty}(\mrI) , \ j = 0, 1, \ldots, k \right\} , 
\label{defWinftw}  \\
\| f \|_{W^{k , \infty}_{w}} &:= \ \max_{0 \le j \le k} \| (1 - x)^{j/2} x^{j/2} D^{j}f(x) \|_{L^{\infty}(\mrI)} \, .
\label{defWinftwnorm}
\end{align}
The subscript $w$ denotes the fact that $W^{k , \infty}_{w}(\mrI)$ is a weaker space than $W^{k , \infty}(\mrI)$ in that
the derivative of functions in $W^{k , \infty}_{w}(\mrI)$ may be unbounded at the endpoints of the interval. 
\begin{lemma} \label{lmaprodsp}
Let $\alpha, \, \beta > -1$, $0 \le s \le k \in \mathbb{N}_{0}$, and $f \in W^{k , \infty}_{w}(\mrI)$. Then, for
\begin{align}
g &\in H^{s}_{(\alpha , \beta)}(\mrI) \mbox{ we have that  } \ f g \in H^{s}_{(\alpha , \beta)}(\mrI) \, .
   \label{prodr1}  
%   \\
%(ii)  \ g &\in H^{-s}_{(\alpha , \beta)}(\mrI) \mbox{ then  } \ f g \in H^{-s}_{(\alpha , \beta)}(\mrI) . 
%   \label{prodr2} 
\end{align}
\end{lemma}
\textbf{Proof}: We establish Lemma \ref{lmaprodsp} for $s = 0$, $s = 1$, and $0 < s < 1$. The proof extends in an
obvious manner for $s > 1$. \\
If $s = 0$ then $k \ge 0$, and
\begin{align}
\| f g \|_{H^{s}_{(\alpha , \beta)}(\mrI)}^{2} &= \ \| f g \|_{H^{0}_{(\alpha , \beta)}(\mrI)}^{2}
  \  = \ \int_{0}^{1} (1 - x)^{\alpha} x^{\beta} \, (f(x) g(x))^{2} \, dx   \nonumber  \\
&\lesssim \| f \|_{L^{\infty}(\mrI)}^{2} \, \| g \|_{H^{0}_{(\alpha , \beta)}(\mrI)}^{2} 
\ = \  \| f \|_{W^{0 , \infty}_{w}(\mrI)}^{2} \, \| g \|_{H^{0}_{(\alpha , \beta)}(\mrI)}^{2}  \, .  \label{erwy1}
\end{align}
Hence \eqref{prodr1} is established for $s = 0$.

For $s = 1$ then $k \ge 1$, and
\begin{align}
\| f g \|_{H^{1}_{(\alpha , \beta)}(\mrI)}^{2} &= \ \| f g \|^{2}_{L^{2}_{(\alpha , \beta)}(\mrI)} \ + \ 
 \| D( f g ) \|^{2}_{L^{2}_{(\alpha + 1 \,  , \, \beta + 1)}(\mrI)}   \nonumber  \\
&= \ 
\int_{0}^{1} (1 - x)^{\alpha} x^{\beta} \, \left( f(x) g(x) \right)^{2} \, dx \ + \ 
\int_{0}^{1} (1 - x)^{\alpha + 1} x^{\beta + 1} \, \left( g(x) \, Df(x)  \ + \ f(x) \, Dg(x) \right)^{2} \, dx    \nonumber  \\
&\le \ 
\int_{0}^{1}  f(x)^{2} \ (1 - x)^{\alpha} x^{\beta} \, \left( g(x) \right)^{2} \, dx \ + \ 
2 \, \int_{0}^{1} (1 - x)^{1} x^{1} \, \left( Df(x) \right)^{2} \
 (1 - x)^{\alpha} x^{\beta} \, \left( g(x) \right)^{2}  dx     \nonumber  \\
& \quad \quad  \quad \quad  \quad \quad 
+ \ 2 \, \int_{0}^{1} f(x)^{2} \
 (1 - x)^{\alpha + 1} x^{\beta + 1} \, \left( Dg(x) \right)^{2}  dx     \nonumber  \\
&\lesssim \ 
 \left( \| f \|^{2}_{L^{\infty}(\mrI)} \ + \ \| (1 - x)^{1/2} x^{1/2} \,  Df(x) \|^{2}_{L^{\infty}(\mrI)} \right) \, 
\left( \| g \|^{2} _{L^{2}_{(\alpha , \beta)}(\mrI)} \ + \ \| Dg \|^{2} _{L^{2}_{(\alpha + 1 \, , \, \beta + 1)}(\mrI)} \right)    \nonumber  \\
&= \  \| f \|^{2}_{W^{1 , \infty}_{w}(\mrI)} \, \| g \|^{2}_{H^{1}_{(\alpha , \beta)}(\mrI)} \, .  \label{erwy2}
\end{align}
Hence \eqref{prodr1} is established for $s = 1$.

Note that for $f \in W^{1 , \infty}_{w}(\mrI)$, from \eqref{erwy1}, 
$\| f g \|_{H^{0}_{(\alpha , \beta)}(\mrI)} \ \le \ \| f \|_{W^{1 , \infty}_{w}(\mrI)} \, \| g \|_{H^{0}_{(\alpha , \beta)}(\mrI)}$.
Combining this with \eqref{erwy2} and the fact that $H^{s}_{(\alpha , \beta)}(\mrI)$ are interpolation spaces, it follows
that for $0 < s < 1$,  
$\| f g \|_{H^{s}_{(\alpha , \beta)}(\mrI)} \ \le \ \| f \|_{W^{1 , \infty}_{w}(\mrI)} \, \| g \|_{H^{s}_{(\alpha , \beta)}(\mrI)}$. \\
\mbox{  } \hfill \qed

%%For (ii), let $g \in H^{-s}_{(\alpha , \beta)}(\mrI)$. Using $L^{2}$ duality, there exists $\hat{g} \in L^{2}_{(\alpha , \beta)}(\mrI)$
%%such that for $h \in H^{s}_{(\alpha , \beta)}(\mrI)$, $g(h) \ = \ ( \hat{g} , h )_{L^{2}_{(\alpha , \beta)}(\mrI)}$. Then, 
%%\[
%%f g(h) \ = \ ( f \hat{g}  \, ,  \, h )_{L^{2}_{(\alpha , \beta)}(\mrI)} \ = \ ( \hat{g}  \, ,  \,  f h )_{L^{2}_{(\alpha , \beta)}(\mrI)}
%%\ = \ g( f h )  < \infty \, ,
%%\]
%%as from (i) $f h \, \in \, H^{s}_{(\alpha , \beta)}(\mrI)$, from which it follows that $f g \, \in \, H^{-s}_{(\alpha , \beta)}(\mrI)$. \\
%%\mbox{  } \hfill \qed
 
Arising in the analysis is the product of a function and a functional (e.g., $b(x) \, D u(x)$). To define such a product
it is convenient to consider $L^{2}_{(\alpha , \beta)}(\mrI) = H^{0}_{(\alpha , \beta)}(\mrI)$ as the pivot space
for $H^{s}_{(\alpha , \beta)}(\mrI)$ and $H^{-s}_{(\alpha , \beta)}(\mrI)$, with $H^{-s}_{(\alpha , \beta)}(\mrI)$
characterized as the closure of $L^{2}_{(\alpha , \beta)}(\mrI)$ with respect to the operator norm
\[
    |\| v \|| \ := \ \sup_{h \in H^{s}_{(\alpha , \beta)}(\mrI)} 
    \frac{ \langle v \, , \, h \rangle_{\rho^{(\alpha , \beta)}}}{ \| h \|_{H^{s}_{(\alpha , \beta)}(\mrI)}} \, .
\]

\textbf{Definition}: Product of a function and a functional. \\
Let  $g \in H^{-s}_{(\alpha , \beta)}(\mrI)$. Then there exists 
$ \{ g_{i} \}_{i = 1}^{\infty} \subset  L^{2}_{(\alpha , \beta)}(\mrI)$, such that
$\lim_{i \rightarrow \infty} |\| g - g_{i} \|| \, = \, 0$. Thus, for any $h \in H^{s}_{(\alpha , \beta)}(\mrI)$,
$g(h) \ = \ \lim_{i \rightarrow \infty} \langle g_{i} \, , \, h \rangle_{\rho^{(\alpha , \beta)}}$.

For $f \in W^{k , \infty}_{w}(\mrI)$, let $f g$ be defined as
\be
f  g (h) \ := \ \lim_{i \rightarrow \infty} \langle f g_{i} \, , \, h \rangle_{\rho^{(\alpha , \beta)}} \, .
\label{deffg1}
\ee
\begin{lemma} \label{lmaprodsp2}
Let $\alpha, \, \beta > -1$, $0 \le s \le k \in \mathbb{N}_{0}$, and $f \in W^{k , \infty}_{w}(\mrI)$. Then, for 
\begin{align}
g &\in H^{-s}_{(\alpha , \beta)}(\mrI) \mbox{ we have that  } \ f g \in H^{-s}_{(\alpha , \beta)}(\mrI) . 
   \label{prodr2} 
\end{align}
\end{lemma}
\textbf{Proof}: 
To establish that $f g$ is well defined and contained in $H^{-s}_{(\alpha , \beta)}(\mrI)$, we show that
$\{ f g_{i} \}_{i = 1}^{\infty}$ is a Cauchy sequence with respect to the norm $|\| \cdot \||$. 

As \eqref{prodr2} trivially holds for $f = 0$, assume $0 \neq f \in W^{k , \infty}_{w}(\mrI)$. Then, 
\begin{align*}
| \| f g_{i} \, - \, f g_{j} \| | &= \ \sup_{h \in H^{s}_{(\alpha , \beta)}(\mrI)} 
    \frac{ \langle f g_{i} \, - \, f g_{j} \, , \, h \rangle_{\rho^{(\alpha , \beta)}}}{ \| h \|_{H^{s}_{(\alpha , \beta)}(\mrI)}}  
    \ = \  \| f \|_{W^{k , \infty}_{w}(\mrI)} \, \sup_{h \in H^{s}_{(\alpha , \beta)}(\mrI)} 
    \frac{ \langle g_{i}  - g_{j} \, , \, f h \rangle_{\rho^{(\alpha , \beta)}}}%
    {\| f \|_{W^{k , \infty}_{w}(\mrI)} \,  \| h \|_{H^{s}_{(\alpha , \beta)}(\mrI)}}  
\\
    &\le \   \| f \|_{W^{k , \infty}_{w}(\mrI)} \,  \sup_{h \in H^{s}_{(\alpha , \beta)}(\mrI)} 
    \frac{ \langle g_{i} -  g_{j} \, , \, f h \rangle_{\rho^{(\alpha , \beta)}}}{ \| f  h \|_{H^{s}_{(\alpha , \beta)}(\mrI)}} 
    \ \mbox{ (using Lemma \ref{lmaprodsp})},  \\
    &\le \
       \| f \|_{W^{k , \infty}_{w}(\mrI)} \,   \sup_{\tilde{h} \in H^{s}_{(\alpha , \beta)}(\mrI)} 
    \frac{ \langle g_{i} -  g_{j} \, , \, \tilde{h} \rangle_{\rho^{(\alpha , \beta)}}}{ \| \tilde{h} \|_{H^{s}_{(\alpha , \beta)}(\mrI)}} 
     \\
  &\le \   \| f \|_{W^{k , \infty}_{w}(\mrI)} \,  | \| g_{i} -  g_{j}  \| | \, .
\end{align*}
As $\{ g_{i} \}_{i = 1}^{\infty}$ is a Cauchy sequence in $H^{-s}_{(\alpha , \beta)}(\mrI)$, then it follows that
$\{ f g_{i} \}_{i = 1}^{\infty}$ is also a Cauchy sequence with limit in $H^{-s}_{(\alpha , \beta)}(\mrI)$. Hence,
\eqref{deffg1} defines a linear functional in $H^{-s}_{(\alpha , \beta)}(\mrI)$. \\
\mbox{ } \hfill \qed

\begin{theorem} \label{thmreg11}
Let $s > -1$, $\beta$ be determined by \textbf{Condition A}, 
$c \in W^{\ceil{ \min \{s \, , \, \alpha \, + \, (\alpha - \beta) \, + \, 1 \, , \,  \alpha \, + \, \beta \, + \, 1\}}  , \infty}_{w}(\mrI)$ 
satisfying $c(x) \ge 0$
and
\begin{equation}
f \in H^{-\alpha/2}(\mrI) \cap H^{s}_{(\beta \, , \, \alpha - \beta)}(\mrI).
\label{def4f1}
\end{equation}  
Then there exists a unique solution $u(x) \ = \ (1 - x)^{\alpha - \beta} \, x^{\beta} \, \phi(x)$, with  \linebreak[4]
$\phi(x) \in 
H^{\alpha \, + \, \min \{s \, , \, \alpha \, + \, (\alpha - \beta) \, + \, 1 - \eps \, , \,  \alpha \, + \, \beta \, + \, 1 - \eps\}}_{(\alpha - \beta \, , \, \beta)}(\mrI)$ for arbitrary $\eps > 0$,
 to
\begin{equation}
\mcL_{r}^{\alpha} u(x) \ + \ c(x) \, u(x) \ = \ f(x) \, , \ x \in \mrI , \ 
\mbox{ subject to } u(0) = u(1) = 0 \, .
\label{theq1}
\end{equation}
\end{theorem}
\textbf{Proof}:  
The stated result is established using two steps. In Step 1 existence of a solution $u \in H_{0}^{\alpha/2}(\mrI)$
to \eqref{theq1} is shown. Then, in Step 2 a boot strapping argument is applied 
to  improve the regularity of $u$.

\underline{Step 1}: 
For $f$ satisfying \eqref{def4f1},
from \cite{erv061}, there exists a unique solution $u \in H^{\alpha/2}_{0}(\mrI)$ to \eqref{theq1}.

\underline{Step 2}: 
For $u \in H_{0}^{\alpha/2}(\mrI)$, from \cite[Theorem 1.2.16]{gri921}, $u(x) \ = \ (1 - x)^{\alpha/2} \, x^{\alpha/2} \, g(x)$,
where $g \in L^{2}(\mrI)$.    \\
Then, as for $0 < \epsilon \le \min\{ \alpha - \beta \, , \, \beta \}$
\[
\int_{0}^{1} (1 - x)^{-1 + \epsilon} \, x^{-1 + \epsilon} \, \left( u(x) \right)^{2} \, dx 
\ = \ \int_{0}^{1} (1 - x)^{-1 + \epsilon + \alpha} \, x^{-1 + \epsilon + \alpha} \, \left( g(x) \right)^{2} \, dx 
\ < \ \int_{0}^{1} \left( g(x) \right)^{2} \, dx \ < \ \infty \, ,
\]
it follows that 
$u \in H^{0}_{(-1 + \eps \, , \, -1 + \eps)}(\mrI) \subset H^{0}_{(\beta - 1\, , \, \alpha - \beta - 1)}(\mrI)
 \subset H^{0}_{(\beta \, , \, \alpha - \beta)}(\mrI)$,
and using Lemma \ref{lmaprodsp} (with the association 
$g(x) = u(x) \in H^{0}_{(\beta \, , \, \alpha - \beta)}(\mrI)$, $f(x) = c(x) \in W^{0  , \infty}_{w}(\mrI)$)
\be
c(x) \, u(x)  \in H^{0}_{(\beta \, , \, \alpha - \beta)}(\mrI) \, .
\label{regu1}
\ee

Using \eqref{regu1}, the solution $u$ of \eqref{theq1} satisfies
\[
 \mcL_{r}^{\alpha} u(x)  \ = \ f(x)  \ - \ c(x) \, u(x) \ := \
  f_{1}(x) \in H_{(\beta  \, , \, \alpha - \beta)}^{\min\{s , 0\}}(\mrI) \, .
\]
From Corollary \ref{exuncor} it follows that
\[
 u(x) \ = \ (1 - x)^{\alpha - \beta} \, x^{\beta} \, \phi_{1}(x) \, , \mbox{  where }
  \phi_{1} \in H_{(\alpha - \beta  \, , \, \beta)}^{\alpha \, + \, \min\{s \,  , \, 0\}}(\mrI) \, .
\]  

Using Theorem \ref{thrmall1}  
(with its parameters $s$, $\mu$, $p$, $\sigma$ replaced by
 $\min\{s + \alpha \,  , \,  \alpha\}, \  \alpha - \beta  ,  \ \alpha - \beta , \  \beta$, respectively; and in
 the second instance,  with its parameters $s$, $\mu$, $p$, $\sigma$ replaced by
 $\min\{s + \alpha \,  , \,  \alpha\}, \ \beta ,  \ \beta , \ \alpha - \beta$, respectively)
%$(s \, = \, \min\{s + \alpha \,  , \,  \alpha\}, \ \mu \, = \, \alpha - \beta  ,  \, p \, = \, \alpha - \beta , \, \sigma \, = \, \beta \ ; \
%  s \, = \, \min\{s + \alpha \,  , \,  \alpha\}, \ \mu \, = \, \beta  ,  \, p \, = \, \beta , \, \sigma \, = \, \alpha - \beta)$, 
  we have that, for arbitrary $\eps > 0$,
$u \in  H_{(\beta\, , \, \alpha - \beta)}^{\min\{s + \alpha \,  , \, \alpha \, , \, \alpha \, + \, (\alpha - \beta) \, + 1 - \eps \, , \,
\alpha \, + \, \beta \, + \, 1 - \eps \}}(\mrI)$
and using Lemma \ref{lmaprodsp},
\be
c(x) \, u(x)  \in H_{(\beta\, , \, \alpha - \beta)}^{\min\{s  \,  , \, \alpha \, , \, \alpha \, + \, (\alpha - \beta) \, + 1 - \eps \, , \,
\alpha \, + \, \beta \, + \, 1 - \eps \}}(\mrI)
\label{regu2}
\ee

Again, using that the solution $u$ of \eqref{theq1} satisfies
\[
 \mcL_{r}^{\alpha} u(x)  \ = \ f(x)  \ - \ c(x) \, u(x) \ := \
  f_{2}(x) \in H_{(\beta \, , \, \alpha - \beta)}^{\min\{s \,  , \, \alpha \, , \, \alpha \, + \, (\alpha - \beta) \, + 1 - \eps \, , \,
\alpha \, + \, \beta \, + \, 1 - \eps \}}(\mrI) \, , 
\]
and from Corollary \ref{exuncor},
\[
 u(x) \ = \ (1 - x)^{\alpha - \beta} \, x^{\beta} \, \phi_{2}(x) \, , \mbox{  where }
  \phi_{2} \in H_{(\alpha - \beta  \, , \, \beta)}^{\alpha \, + \, \min\{s \, , \,  \alpha \, , \, \alpha \, + \, (\alpha - \beta) \, + 1 - \eps \, , \,
\alpha \, + \, \beta \, + \, 1 - \eps \}}(\mrI)  \, .
\]  

Using Theorem \ref{thrmall1}  
(with its parameters $s$, $\mu$, $p$, $\sigma$ replaced by
$\min\{s + \alpha \,  , \,  2 \alpha\}, \  \alpha - \beta ,  \ \alpha - \beta , \ \beta$, respectively; and in
 the second instance,  with its parameters $s$, $\mu$, $p$, $\sigma$ replaced by
$\min\{s + \alpha \,  , \,  2 \alpha\}, \  \beta  ,  \ \beta , \ \alpha - \beta$, respectively),
%%s \, = \, \min\{s + \alpha \,  , \,  2 \alpha\}, \ \mu \, = \, \alpha - \beta  ,  \, p \, = \, \alpha - \beta , \, \sigma \, = \, \beta \ ; \
%%  s \, = \, \min\{s + \alpha \,  , \,  2 \alpha\}, \ \mu \, = \, \beta  ,  \, p \, = \, \beta , \, \sigma \, = \, \alpha - \beta)$, 
 and  Lemma \ref{lmaprodsp},
\[
c(x) \, u(x)  \in H_{(\beta\, , \, \alpha - \beta)}^{\min\{s \,  , \, 2 \alpha \, , \, \alpha \, + \, (\alpha - \beta) \, + 1 - \eps \, , \,
\alpha \, + \, \beta \, + \, 1 - \eps \}}(\mrI) \, ,
\]
from which it then follows that 
\[
 u(x) \ = \ (1 - x)^{\alpha - \beta} \, x^{\beta} \, \phi_{2}(x) \, , \mbox{  where }
  \phi_{2} \in H_{(\alpha - \beta  \, , \, \beta)}^{\alpha \, + \, \min\{s \, , \, 2 \alpha \, , \, \alpha \, + \, (\alpha - \beta) \, + 1 - \eps \, , \,
\alpha \, + \, \beta \, + \, 1 - \eps \}}(\mrI)  \, .
\]  
Noting that $4 \alpha \ge \min\{ \alpha \, + \, (\alpha - \beta) \, + 1 \, , \,
\alpha \, + \, \beta \, + \, 1 \}$ for $1 < \alpha < 2$, repeating the boot strapping argument two more times establishes the
stated result. \\
\mbox{ } \hfill \qed

%%%%
%%%%

The inclusion of an advection term can significantly reduced the regularity of the solution.

\begin{theorem} \label{thmreg13}
Let $s > -1$, $\beta$ be determined by \textbf{Condition A},  \\
$b \in W^{\max\{1 , \, \ceil{\min \{s \, , \, \alpha \, + \, (\alpha - \beta) \, - \, 1 \, , \,  \alpha \, + \, \beta \, - \, 1\} } \}, \infty}_{w}(\mrI)$, 
\, $c \in W^{\ceil{\min \{s \, , \, \alpha \, + \, (\alpha - \beta) \, - \, 1 \, , \,  \alpha \, + \, \beta \, - \, 1\}} , \infty}_{w}(\mrI)$ satisfying 
$c(x)  \, - \, 1/2 D b(x) \ \ge 0$,
and
\begin{equation}
f \in H^{-\alpha/2}(\mrI) \cap H^{s}_{(\beta \, , \, \alpha - \beta)}(\mrI).
\label{def4f2}
\end{equation}  
Then there exists a unique solution $u(x) \ = \ (1 - x)^{\alpha - \beta} \, x^{\beta} \, \phi(x)$, with  \linebreak
$\phi(x) \in 
H^{\alpha \, + \, \min \{s \, , \, \alpha \, + \, (\alpha - \beta) \, - \, 1 - \eps \, , \,  \alpha \, + \, \beta \, - \, 1 - \eps\}}_{(\alpha - \beta \, , \, \beta)}(\mrI)$ for arbitrary $\eps > 0$,
to
\begin{equation}
\mcL_{r}^{\alpha} u(x) \ + \ b(x) \, D u(x) \ + \ c(x) \, u(x) \ = \ f(x) \, , \ x \in \mrI , \ 
\mbox{ subject to } u(0) = u(1) = 0 \, .
\label{theq2}
\end{equation}
\end{theorem}
\textbf{Proof}:  
The proof follows the same two steps as in Theorem \label{thmreg2}. Step 1, establishing the
existence of a solution is exactly the same. In Step 2 the boot strapping argument is applied 
$m$ times, where $m$ is the least integer such that $m (\alpha - 1) \, \ge \, 
\min\{ \alpha + (\alpha - \beta) \, , \, \alpha + \beta \}$, to obtain the
stated result.

\underline{Step 2}: 
For $u \in H_{0}^{\alpha/2}(\mrI)$ and $0 < \epsilon \le \min\{ \alpha - \beta \, , \, \beta \}$, 
$u \in H^{0}_{(-1 + \eps \, , \, -1 + \eps)}(\mrI) 
\subset H^{0}_{(\beta - 1 \, , \, \alpha - \beta - 1)}(\mrI)$. Then using Lemma \ref{lmamapD}, 
$D u \in H^{-1}_{(\beta \, , \, \alpha - \beta)}(\mrI)$. Hence we have using \eqref{regu1} 
and \eqref{prodr2},
\[
c(x) \, u(x)  \in H^{0}_{(\beta \, , \, \alpha - \beta)}(\mrI)  \ \mbox{ and } \ 
b(x) \, D u(x)  \in H^{-1}_{(\beta \, , \, \alpha - \beta)}(\mrI) \, .
\]
This leads to the solution of \eqref{theq2} satisfying
\[
 \mcL_{r}^{\alpha} u(x)  \ = \ f(x)  \ - \ b(x) \, D u(x) \ - \ c(x) \, u(x) \ := \
  f_{1}(x) \in H_{(\beta  \, , \, \alpha - \beta)}^{\min\{s , -1\}}(\mrI) \, .
\]
From Corollary \ref{exuncor} it follows that
\[
 u(x) \ = \ (1 - x)^{\alpha - \beta} \, x^{\beta} \, \phi_{1}(x) \, , \mbox{  where }
  \phi_{1} \in H_{(\alpha - \beta  \, , \, \beta)}^{\alpha \, + \, \min\{s \,  , \,   - 1\}}(\mrI) \, .
\]  

Using Theorem \ref{thrmall1} 
(with its parameters $s$, $\mu$, $p$, $\sigma$ replaced by
$ \min\{s + \alpha \,  , \,  \alpha - 1\}, \  \alpha - \beta  ,  \  \alpha - \beta , \ \, \beta - 1$, respectively; and in
 the second instance,  with its parameters $s$, $\mu$, $p$, $\sigma$ replaced by
$\min\{s + \alpha \,  , \,  \alpha - 1\}, \  \beta  ,  \ \beta , \ \alpha - \beta - 1$, respectively)
%%$(s \, = \, \min\{s + \alpha \,  , \,  \alpha - 1\}, \ \mu \, = \, \alpha - \beta  ,  \, p \, = \, \alpha - \beta , \, \sigma \, = \, \beta - 1\ ; \
%%  s \, = \, \min\{s + \alpha \,  , \,  \alpha - 1\}, \ \mu \, = \, \beta  ,  \, p \, = \, \beta , \, \sigma \, = \, \alpha - \beta - 1)$, 
we have that, for arbitrary $\eps > 0$,
$u \in  H_{(\beta - 1 \, , \, \alpha - \beta - 1)}^{\min\{s + \alpha \,  , \, \alpha - 1 \, , \, \alpha \, + \, (\alpha - \beta) - \eps \, , \,
\alpha \, + \, \beta - \eps\}}(\mrI) \ = \  
H_{(\beta - 1 \, , \, \alpha - \beta - 1)}^{\min\{s + \alpha \,  , \, \alpha - 1 \}}(\mrI)$
and  using Lemmas \ref{lmaprodsp} and \ref{lmamapD}
\be
 D u(x)  \in 
H^{-1 + \min\{s + \alpha \,  , \, \alpha - 1  \}}_{(\beta \, , \, \alpha - \beta)}(\mrI) , \ 
c(x) \, u(x)  \in H^{\min\{s \,  , \, \alpha - 1 \}}_{(\beta - 1 \, , \, \alpha - \beta - 1)}(\mrI)  \ \mbox{ and } \ 
b(x) \, D u(x)  \in 
H^{ \min\{s \,  , \, \alpha - 2  \}}_{(\beta \, , \, \alpha - \beta)}(\mrI) \, .
\label{regu22}
\ee

The solution $u$ of \eqref{theq1} then must satisfies
\[
 \mcL_{r}^{\alpha} u(x)  \ = \ f(x)  \ - \ b(x) \, D u(x) \ - \ c(x) \, u(x) \ := \
  f_{2}(x) \in H_{(\beta \, , \, \alpha - \beta)}^{\min\{s \, , \,  \alpha - 2 \}}(\mrI) \, , 
\]
and from Corollary \ref{exuncor},
\[
 u(x) \ = \ (1 - x)^{\alpha - \beta} \, x^{\beta} \, \phi_{2}(x) \, , \mbox{  where }
  \phi_{2} \in H_{(\alpha - \beta  \, , \, \beta)}^{\alpha \, + \, 
  \min\{s  \, , \,  \alpha - 2  \}}(\mrI)  \, .
\]  

Using Theorem \ref{thrmall1}  
(with its parameters $s$, $\mu$, $p$, $\sigma$ replaced by
$min\{s + \alpha \,  , \,  2\alpha - 2\}, \  \alpha - \beta  ,  \ \alpha - \beta , \ \beta - 1$, respectively; and in
 the second instance,  with its parameters $s$, $\mu$, $p$, $\sigma$ replaced by
$\min\{s + \alpha \,  , \,  2\alpha - 2\}, \ \beta  ,  \ \beta , \ \alpha - \beta - 1$, respectively),
%%$(s \, = \, \min\{s + \alpha \,  , \,  2\alpha - 2\}, \ \mu \, = \, \alpha - \beta  ,  \, p \, = \, \alpha - \beta , \, \sigma \, = \, \beta - 1\ ; \
%%  s \, = \, \min\{s + \alpha \,  , \,  2\alpha - 2\}, \ \mu \, = \, \beta  ,  \, p \, = \, \beta , \, \sigma \, = \, \alpha - \beta - 1)$, 
 we have that
$u \in  H_{(\beta - 1 \, , \, \alpha - \beta - 1)}^{\min\{s + \alpha \,  , \, 2\alpha - 2 \, , \, \alpha \, + \, (\alpha - \beta) - \eps \, , \,
\alpha \, + \, \beta - \eps \}}(\mrI)$,
and   $D u(x)  \in 
H^{-1 + \min\{s + \alpha \,  , \, 2\alpha - 2  \, , \, \alpha \, + \, (\alpha - \beta) - \eps \, , \,
\alpha \, + \, \beta - \eps \}}_{(\beta \, , \, \alpha - \beta)}(\mrI)$,
\be
c(x) \, u(x)  \in H^{\min\{s \,  , \, 2\alpha - 2  \, , \, \alpha \, + \, (\alpha - \beta) - \eps \, , \,
\alpha \, + \, \beta - \eps \}}_{(\beta - 1 \, , \, \alpha - \beta - 1)}(\mrI)  \ \mbox{ and } \ 
b(x) \, D u(x)  \in 
H^{ \min\{s \,  , \, 2\alpha - 3  \, , \, \alpha \, + \, (\alpha - \beta)  - 1 - \eps \, , \,
\alpha \, + \, \beta  - 1 - \eps \}}_{(\beta \, , \, \alpha - \beta) - 1}(\mrI) \, .
\label{regu25}
\ee

The solution $u$ of \eqref{theq1} then must satisfies
\[
 \mcL_{r}^{\alpha} u(x)  \ = \ f(x)  \ - \ b(x) \, D u(x) \ - \ c(x) \, u(x) \ := \
  f_{2}(x) \in H_{(\beta \, , \, \alpha - \beta)}^{\min\{s \,  , \, 2\alpha - 3  \, , \, \alpha \, + \, (\alpha - \beta)  - 1 - \eps \, , \,
\alpha \, + \, \beta  - 1 - \eps \}}(\mrI) \, , 
\]
and from Corollary \ref{exuncor},
\[
 u(x) \ = \ (1 - x)^{\alpha - \beta} \, x^{\beta} \, \phi_{3}(x) \, , \mbox{  where }
  \phi_{3} \in H_{(\alpha - \beta  \, , \, \beta)}^{\alpha \, + \, 
  \min\{s \,  , \, 2\alpha - 3  \, , \, \alpha \, + \, (\alpha - \beta)  - 1 - \eps \, , \,
\alpha \, + \, \beta  - 1 - \eps \}}(\mrI)  \, .
\]  

Repeatedly applying this boot stepping procedure we obtain after $(m - 2)$ additional steps
 $u(x) \ = \ (1 - x)^{\alpha - \beta} \, x^{\beta} \, \phi_{m + 1}(x) $,  where 
\[
  \phi_{m + 1} \in  
  H_{(\alpha - \beta  \, , \, \beta)}^{\alpha \, + \, \min \{s \, , \, m \alpha - (m + 1) \, , \, 
  \alpha \, + \, 
  (\alpha - \beta) \, - \, 1 - \eps  \, , \,  \alpha \, + \, \beta \, - \, 1 - \eps \}}(\mrI) \ = \ 
 H_{(\alpha - \beta  \, , \, \beta)}^{\alpha \, + \, \min \{s \, , \, 
  \alpha \, + \, 
  (\alpha - \beta) \, - \, 1 - \eps \, , \,  \alpha \, + \, \beta \, - \, 1 - \eps \}}(\mrI)   \, .
\]  
\mbox{ } \hfill \qed

%%%%%
%%%%% 
 \setcounter{equation}{0}
\setcounter{figure}{0}
\setcounter{table}{0}
\setcounter{theorem}{0}
\setcounter{lemma}{0}
\setcounter{corollary}{0}
\setcounter{definition}{0}
\section{Regularity of the solution to the fractional diffusion, advection, reaction equation in 
unweighted Hilbert spaces}
\label{secRFDAR3}

To connect the regularity results for the solution of \eqref{DefProb2}, given in Theorems 
\ref{thmreg11} and \ref{thmreg13} to the usual (unweighted) Hilbert spaces we use four steps. 
Step 1 uses Theorem \ref{thrmall1} to determine $H^{t}_{(\sigma)}(\mrJ)$ such that
$u_{0}(x) \, = \, x^{\beta} \psi_{0}(x)  \in  H^{t}_{(\sigma)}(\mrJ)$, for 
$\psi_{0}(x)  \in  H^{s^{*}}_{(\beta)}(\mrJ)$. Step 2 applies an embedding theorem
(Corollary \ref{embErv}) to then obtain $u_{0} \in  H^{v_{0}}(\mrJ)$.
Step 3 repeats Steps 1 and 2 for $u_{1}(x) \, = \, x^{\alpha - \beta} \psi_{1}(x)$, for 
$\psi_{1}(x)  \in  H^{s^{*}}_{(\alpha - \beta)}(\mrJ)$ to obtain
$u_{1} \in  H^{v_{1}}(\mrJ)$. The final step combines Steps 2 and 3 to conclude
that $u \in H^{\min\{v_{0} , v_{1}\}}(\mrI)$.

To begin we introduce the space $W^{s , 2}_{(c , d)}(\mrI) \, := \, 
\{ f : f \mbox{ is  measurable and } \| f \|_{W^{s , 2}_{(c , d)}(\mrI)} < \infty \}$, where
\begin{align}
  \| f \|_{W^{s , 2}_{(c , d)}(\mrI)}^{2} &:= \left\{ \begin{array}{rl}  
   \sum_{j = 0}^{s} \| D^{j} f \|_{L^{2}_{(c \, , \, d)}(\mrI)}^{2} \, , & \mbox{ for } s \in \mathbb{N}_{0}  \\
 \sum_{j = 0}^{\floor{s}} \| D^{j} f \|_{L^{2}_{(c \, , \, d)}(\mrI)}^{2} \ + \ 
  | f |_{W^{s , 2}_{(c , d)}(\mrI)}^{2} \, , & \mbox{ for } s \in \mathbb{R}^{+} \backslash \mathbb{N}_{0} 
 \end{array} \right. \, ,     \nonumber \\
 \mbox{for } \
  | f |_{W^{s , 2}_{(c , d)}(\mrI)}^{2} &:= \
   \iint_{\widetilde{\Lambda}} (1 - x)^{c} \, x^{d} \, 
  \frac{ | D^{\floor{s}} f(x) \, - \, D^{\floor{s}} f(y) |^{2}}{ | x \, - \, y |^{1 \, + \, 2 (s - \floor{s})}} dy \,  dx  \, , \nonumber
\end{align}
and $\widetilde{\Lambda}$  as defined in \eqref{defLtda}.

Following \eqref{deffnJ}, also introduce $W^{s , 2}_{(\delta)}(\mrJ) \, := \, 
\{ f : f \mbox{ is  measurable and } \| f \|_{W^{s , 2}_{(\delta)}(\mrJ)} < \infty \}$, where
\begin{align*}
  \| f \|_{W^{s , 2}_{(\delta)}(\mrJ)}^{2} &:= \  \left\{ \begin{array}{rl}  
  \sum_{j = 0}^{s} \| D^{j} f \|_{L^{2}_{(\delta)}(\mrJ)}^{2} \, , & \mbox{ for } s \in \mathbb{N}_{0}  \\
 \sum_{j = 0}^{\floor{s}} \| D^{j} f \|_{L^{2}_{(\delta)}(\mrJ)}^{2} \ + \  | f |_{W^{s , 2}_{(\delta)}(\mrJ)}^{2} \, , 
 & \mbox{ for } s \in \mathbb{R}^{+} \backslash \mathbb{N}_{0} 
 \end{array} \right. \, ,    \nonumber \\
\mbox{and } \
 | f |_{W^{s , 2}_{(\delta)}(\mrJ)}^{2} &:= \ 
  \iint_{\Lambda}  x^{\delta} \, 
  \frac{ | D^{\floor{s}} f(x) \, - \, D^{\floor{s}} f(y) |^{2}}{ | x \, - \, y |^{1 \, + \, 2 (s - \floor{s})}} dy \,  dx  \ + \ 
  \iint_{\Lambda_{1}}  x^{\delta} \, 
  \frac{ | D^{\floor{s}} f(x) \, - \, D^{\floor{s}} f(y) |^{2}}{ | x \, - \, y |^{1 \, + \, 2 (s - \floor{s})}} dy \,  dx \, .
\end{align*}

From \cite{ber921} we have the following embedding result.
\begin{theorem}\cite[See Theorem 1.d.2]{ber921}    \label{embBDY}
Let $\mu$, $\sigma \, > -1$, and $v$, $w$ be two real numbers such that $0 \le v \le w$. Then, if
\begin{equation}
\left\{ \begin{array}{lcl}
v - \frac{\sigma}{2} & < & w - \frac{\mu}{2}   \\
\mbox{ or } & & \\
v - \frac{\sigma}{2} & = & w - \frac{\mu}{2} \ \mbox{ with } \ w - \frac{\mu}{2} - \frac{1}{2} \not \in \mathbb{N} \, 
\end{array} \, ,  \right.
\end{equation} 
we have  $W^{w , 2}_{(\mu)}(\mrJ) \subset W^{v , 2}_{(\sigma)}(\mrJ)$.
\end{theorem}

\begin{corollary}  \label{embErv}
Let $\gamma \, > -1$, and $v$, $w$ be two real numbers such that $0 \le v \le w$. Then, if
\begin{equation} 
\left\{ \begin{array}{lcl}
v  & < & \frac{w \, - \, \gamma}{2}   \\
\mbox{ or } & & \\
v  & = & \frac{w \, - \, \gamma}{2} \ \mbox{ with } \  \frac{w \, -  \, \gamma}{2} - \frac{1}{2} \not \in \mathbb{N} \, 
\end{array} \, ,  \right.
\end{equation} 
we have  $H^{w}_{(\gamma)}(\mrJ) \subset H^{v}(\mrJ)$.
\end{corollary}
\textbf{Proof}: From \cite[Theorem 3.3]{nic001}, it follows that $H^{s}_{(a , b)}(\mrI)$ and 
$W^{s , 2}_{(a + s , b + s)}(\mrI)$ are equivalent spaces, as are $H^{s}_{(\gamma)}(\mrJ)$ and 
$W^{s , 2}_{(\gamma + s)}(\mrJ)$. Using that $H^{r}(\mrI)$ and $W^{r , 2}_{(0 , 0)}(\mrI)$ are equivalent spaces,
the stated result follows from Theorem \ref{embBDY} for $\sigma = 0$ and $\mu \, = \, w + \gamma$. \\
\mbox{ } \hfill \qed

\begin{figure}[!ht]
%\begin{figure}[t]
\begin{minipage}{.46\linewidth}

%\begin{figure}[!ht]
\begin{center}
 \includegraphics[height=2.5in]{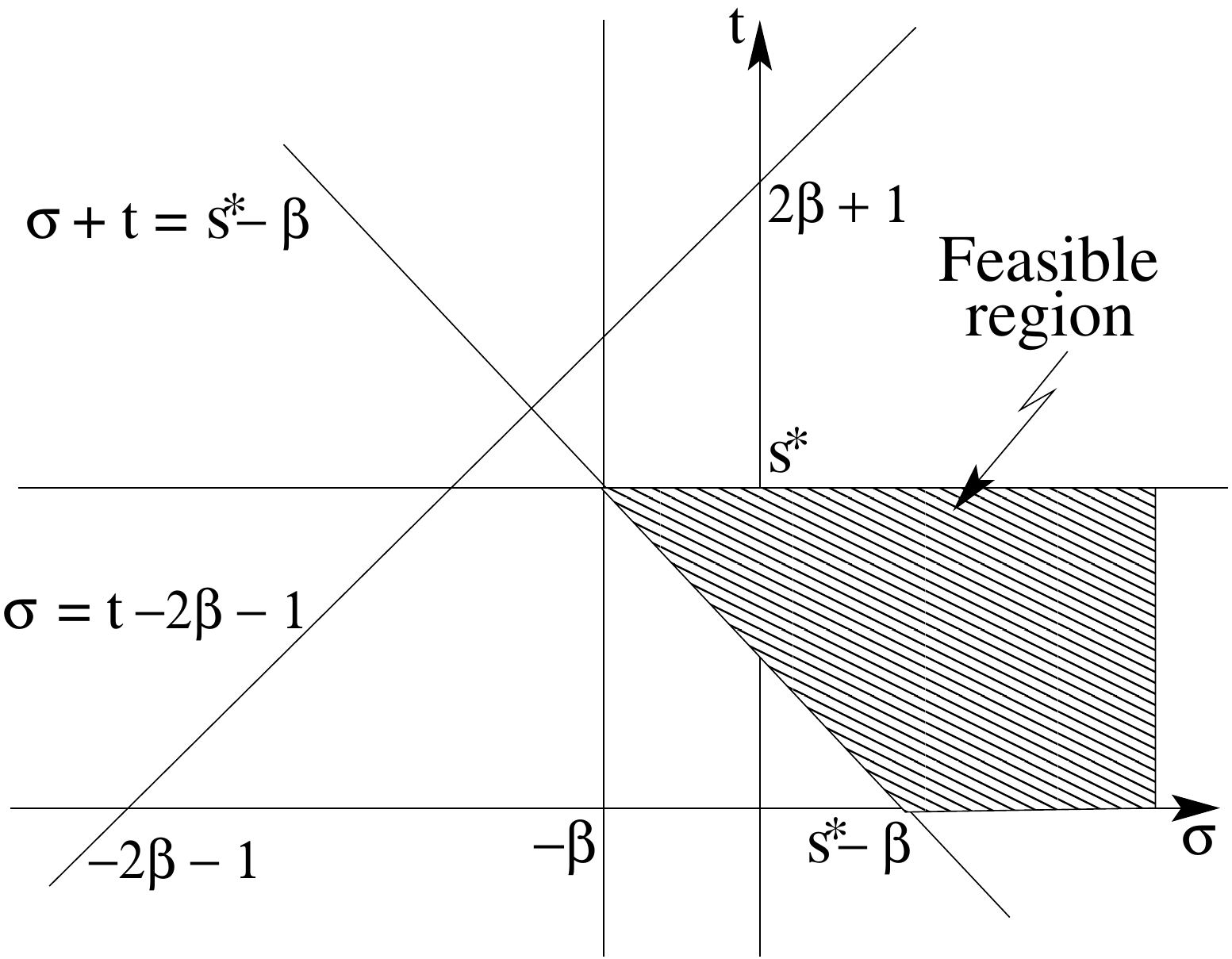}
   \caption{Illustration of Case 1 for constraints \eqref{yute1}.}
   \label{figcase1}
\end{center}
%\end{figure}

\end{minipage} \hfill
\begin{minipage}{.46\linewidth}
 
\begin{center}
 \includegraphics[height=2.5in]{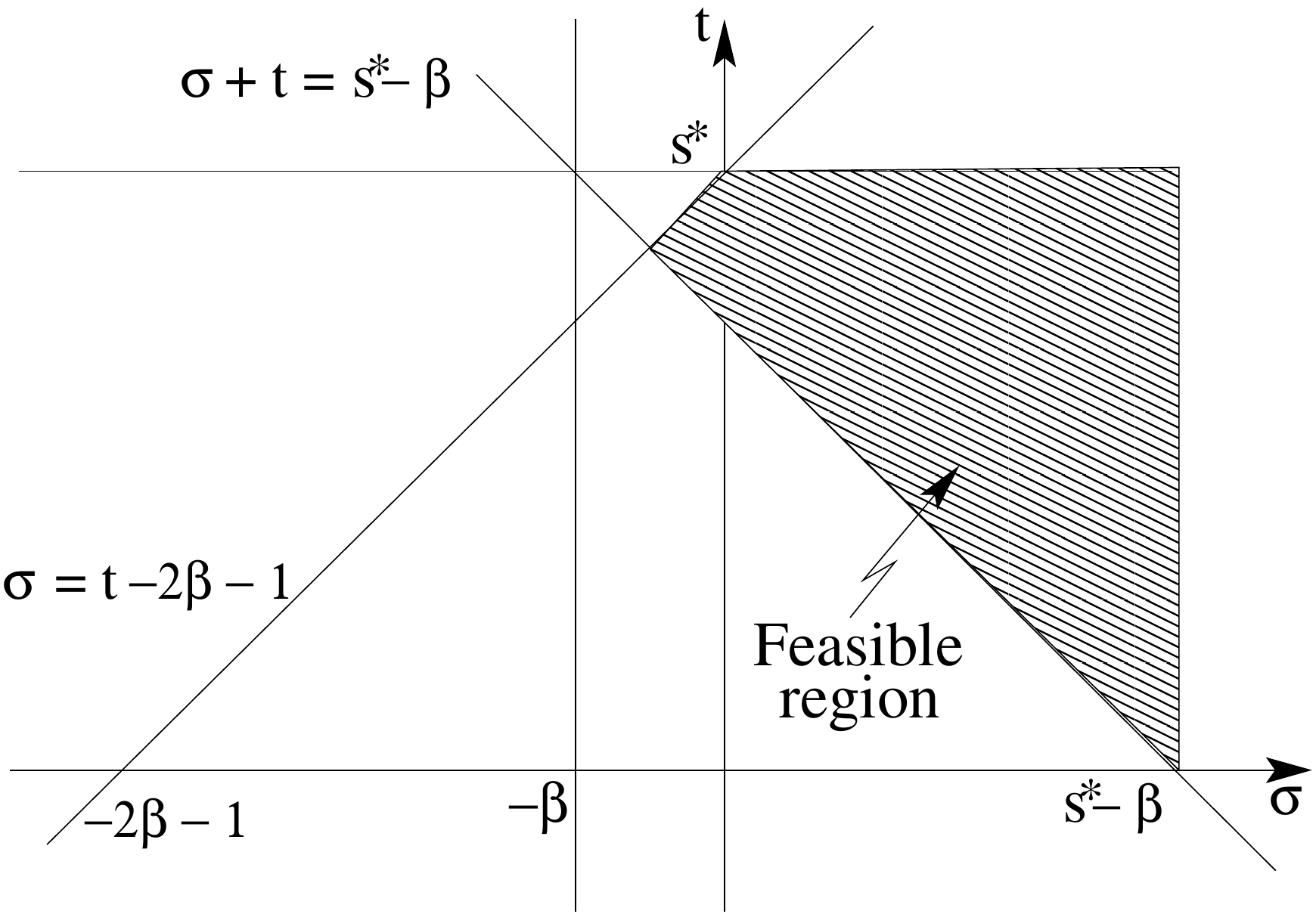}
   \caption{Illustration of Case 2 for constraints \eqref{yute1}.}
   \label{figcase2}
\end{center}
\end{minipage} 
\end{figure}

\begin{corollary} [See Corollary \ref{exuncor}.]  \label{corHreg1}
Let $s \ge \, - \alpha$, $f \in H^{s}_{(\beta \, , \, \alpha - \beta)}(\mrI)$, and $s^{*} \, := \, s \, + \, \alpha$.
Then the unique solution to \eqref{DefProb1},\eqref{DefBC1}, satisfies for any $\epsilon > 0$
\be
 u \in H^{\min\{ \frac{s^{*} + (\alpha - \beta)}{2} , \, \frac{s^{*} +  \beta}{2} , 
 \, (\alpha - \beta) + \frac{1}{2} - \epsilon , \, \beta + \frac{1}{2} - \epsilon \} }(\mrI) \, .
\label{rregu1}
\ee
In particular, for $s \, > \, - \alpha + 1 + \min\{(\alpha - \beta) , \, \beta\}$,
\be
 u \in H^{\min\{(\alpha - \beta)  , \, \beta\} \, + \, \frac{1}{2} - \epsilon}(\mrI) \, .
\label{rregu2}
\ee
\end{corollary}
\textbf{Proof}:
Proceeding as described at the beginning of this section, consider
$u_{0}(x) \, = \, x^{\beta} \psi_{0}(x)$, for 
$\psi_{0}(x)  \in  H^{s^{*}}_{(\beta)}(\mrJ)$. Using Theorem \ref{thrmall1}, the most regular (i.e., ``nicest'') weighted
Sobolev space that $u_{0}$ lies in is  given by the largest value for $t$ and the smallest value
for $\sigma$ such that the conditions stated in \eqref{piu101} are satisfied. To apply Theorem \ref{thrmall1} in
this case we have: $s \rightarrow s^{*}$, $\mu \rightarrow \beta$, $p \rightarrow \beta$. 
Equation \eqref{piu101} then require
that $\sigma$ and $t$ satisfy
\be
0 \le t \le s^{*} , \quad \sigma \ge - \beta , \quad \sigma \, >  \, t \, - \, 2 \beta - 1 , 
\quad \sigma \, \ge \, s^{*} \, - \, t \, - \, \beta \, .
\label{yute1}
\ee
Two cases arise for consideration.\\
\underline{Case 1}. If $s^{*} < \beta + 1$ then $t$ and $\sigma$ satisfying \eqref{yute1} are determined by:
$0 \le t \le s^{*}$, and $\sigma \ge -\beta$ (see Figure \ref{figcase1}). \\
With the choices $t = s^{*}$, $\sigma = -\beta$, using Corollary \ref{embErv} we obtain
\be
u_{0} \in H^{s^{*}}_{(-\beta)}(\mrJ) \subset H^{\frac{s^{*} + \beta}{2}}(\mrJ) \, .
\label{yhgf1}
\ee
\underline{Case 2}. If $s^{*} \ge \beta + 1$ then $t$ and $\sigma$ satisfying \eqref{yute1} are determined by:
$0 \le t \le s^{*}$, $\sigma \, > \, t \, - \, 2 \beta \, - 1$, and $\sigma \, \ge \, s^{*} - t - \beta$ (see Figure \ref{figcase2}).
With $t \le s^{*}$ and $t - \sigma \, < \, 2 \beta \, + \, 1$, using Corollary \ref{embErv} we obtain, for $\epsilon > 0$,
\be
u_{0} \in H^{t}_{(\sigma)}(\mrJ) \subset H^{\beta + \frac{1}{2} - \epsilon}(\mrJ) \, .
\label{yhgf2}
\ee

Combining \eqref{yhgf1} and \eqref{yhgf2} yields
\be
u_{0} \in  H^{\min\{ \frac{s^{*} + \beta}{2} \, , \,  \beta + \frac{1}{2} - \epsilon \}}(\mrJ) \, .
\label{yhgf3}
\ee
For $u_{1}(x) \, = \, x^{\alpha - \beta} \, \psi_{1}(x)$, , for 
$\psi_{1}(x)  \in  H^{s^{*}}_{(\alpha - \beta)}(\mrJ)$, a similar analysis
leads to  
\be
u_{1} \in  H^{\min\{ \frac{s^{*} + (\alpha - \beta)}{2} \, , \,  (\alpha - \beta) + \frac{1}{2} - \epsilon \}}(\mrJ) \, .
\label{yhgf4}
\ee
Combining \eqref{yhgf3} and \eqref{yhgf4} we obtain
\be
u \in  H^{\min\{ \frac{s^{*} + (\alpha - \beta)}{2} \, , \,   \frac{s^{*} + \beta}{2} \, , 
(\alpha - \beta) + \frac{1}{2} - \epsilon \, , \,  \beta + \frac{1}{2} - \epsilon \}}(\mrI) \, .
\label{yhgf5}
\ee
Noting that for $s \, > \, -\alpha + 1 + \, \min\{(\alpha - \beta) \, , \, \beta\}$ that 
$s^{*} \, > \, \min\{(\alpha - \beta) \, , \, \beta\} \, + \, 1$, from \eqref{yhgf5},
\[
  u \in  H^{\min\{(\alpha - \beta)  \, , \, \beta \} \, + \, \frac{1}{2} - \epsilon}(\mrI) \, .
\]
\mbox{ } \hfill \qed

Corresponding to Theorem \ref{thmreg11} we have the following.
\begin{corollary} [See Theorem \ref{thmreg11}.]
Assuming the hypothesis of Theorem \ref{thmreg11} are satisfied, and let 
$s^{*} \, := \, \alpha \, + \, \min \{s \, , \, \alpha \, + \, (\alpha - \beta) \, + \, 1 \, , \,  \alpha \, + \, \beta \, + \, 1\}$.
Then
the unique solution of \eqref{theq1} satisfies for any $\epsilon > 0$
\be
 u \in H^{\min\{ \frac{s^{*} + (\alpha - \beta)}{2} , \, \frac{s^{*} +  \beta}{2} , 
 \, (\alpha - \beta) + \frac{1}{2} - \epsilon , \, \beta + \frac{1}{2} - \epsilon \} }(\mrI) \, .
\label{rregu11}
\ee
In particular, for $s \, > \, - \alpha + 1 + \min\{(\alpha - \beta) , \, \beta\}$,
\be
 u \in H^{\min\{(\alpha - \beta)  , \, \beta\} \, + \, \frac{1}{2} - \epsilon}(\mrI) \, .
\label{rregu12}
\ee
\end{corollary}
\textbf{Proof}: Proof follows exactly as that for Corollary \ref{corHreg1}. \\
\mbox{ } \hfill \qed

%%%%
%%%%%
Corresponding to Theorem \ref{thmreg13} we have the following.
\begin{corollary} [See Theorem \ref{thmreg13}.]
Assuming the hypothesis of Theorem \ref{thmreg13} are satisfied, and let 
$s^{*} \, := \, \alpha \, + \, \min \{s \, , \, \alpha \, + \, (\alpha - \beta) \, - \, 1 \, , \,  \alpha \, + \, \beta \, - \, 1\}$.
Then
the unique solution of \eqref{theq2} satisfies for any $\epsilon > 0$
\be
 u \in H^{\min\{ \frac{s^{*} + (\alpha - \beta)}{2} , \, \frac{s^{*} +  \beta}{2} , 
 \, (\alpha - \beta) + \frac{1}{2} - \epsilon , \, \beta + \frac{1}{2} - \epsilon \} }(\mrI) \, .
\label{rregu13}
\ee
In particular, for $s \, > \, - \alpha + 1 + \min\{(\alpha - \beta) , \, \beta\}$,
\be
 u \in H^{\min\{(\alpha - \beta)  , \, \beta\} \, + \, \frac{1}{2} - \epsilon}(\mrI) \, .
\label{rregu14}
\ee
\end{corollary}
\textbf{Proof}: Proof follows exactly as that for Corollary \ref{corHreg1}. \\
\mbox{ } \hfill \qed

\textbf{Remark}: If the regularity of the right hand side function $f$ is further restricted then
the regularity of the solution may be improved. For example, if $r = 1/2$ the operator
$\mcL_{1/2}^{\alpha}(\cdot)$ corresponds to the integral fractional Laplacian operator. For this 
operator Hao and Zhang in \cite{hao181} showed that for $f \in H^{s}(\mrI)$ the solution
of \eqref{theq2} satisfied $u \in H^{\min\{s + \alpha \, , \, \alpha/2 \, + \, 1/2 \, - \, \epsilon\}}(\mrI)$
for any $\epsilon > 0$.

%%%%

%\bibliographystyle{plain}
%\bibliography{FADEbib}

\end{document}